\documentclass[10pt,a4paper]{article}
\usepackage{graphics}
\usepackage{amsmath,amssymb,amsthm,graphics,graphicx,epic,eepic,multicol,ascmac}
\usepackage[mathscr]{euscript}
\usepackage[all]{xy}
\usepackage{fancybox}
\usepackage{url}
\usepackage{ytableau}
\usepackage{here}

\ytableausetup{centertableaux,mathmode,boxsize=1em}
\def\young(#1){\ytableaushort{#1}}
\def\yng(#1){\ydiagram{#1}}

\numberwithin{equation}{section}
\theoremstyle{theorem}
\newtheorem{thm}{Theorem}[section]
\newtheorem{prop}[thm]{Proposition}
\newtheorem{lem}[thm]{Lemma}
\newtheorem{rem}[thm]{Remark}

\theoremstyle{definition}
\newtheorem{defn}[thm]{Definition}
\newtheorem{ex}[thm]{Example}

\def\al{\alpha}

\def\wht(#1){\widehat{\ #1\ }}

\newcommand{\ch}{\mathrm{ch}}

\newcommand{\lbr}{\begin{bmatrix}}
\newcommand{\rbr}{\end{bmatrix}}

\newcommand{\cd}{commutative diagram }

\def\al{\alpha}

\def\beneme{\begin{enumerate}}
\def\beq{\begin{equation}}
\def\beqn{\begin{eqnarray}}
\def\beqnn{\begin{eqnarray*}}

\def\bfii0{{\bf i_0}}

\def\bbra#1,#2,#3{\left\{\begin{array}{c}\hspace{-5pt}
#1;#2\\ \hspace{-5pt}#3\end{array}\hspace{-5pt}\right\}}
\def\cd{\cdots}
\def\ci(#1,#2){c_{#1}^{(#2)}}
\def\Ci(#1,#2){C_{#1}^{(#2)}}
\def\mpp(#1,#2,#3){#1^{(#2)}_{#3}}
\def\bCi(#1,#2){\ovl C_{#1}^{(#2)}}
\def\ch(#1,#2){c_{#2,#1}^{-h_{#1}}}
\def\cc(#1,#2){c_{#2,#1}}

\def\di(#1,#2){D_{#1}^{(#2)}}
\def\dbi(#1,#2){\ovl D_{#1}^{(#2)}}

\def\eit{\tilde{e}_i}
\def\eneme{\end{enumerate}}

\def\eeq{\end{equation}}
\def\eeqn{\end{eqnarray}}
\def\eeqnn{\end{eqnarray*}}
\def\fit{\tilde{f}_i}

\def\gau#1,#2{\left[\begin{array}{c}\hspace{-5pt}#1\\
\hspace{-5pt}#2\end{array}\hspace{-5pt}\right]}

\def\ify{\infty}
\def\io{\iota}
\def\ji(#1,#2){j_{#1}^{(#2)}}
\def\kp{k^{(+)}}
\def\km{k^{(-)}}

\def\lan{\langle}

\def\nd{\noindent}

\def\ovl{\overline}

\def\qed{\hfill\framebox[2mm]{}}
\def\QQ{\mathbb Q}

\def\ran{\rangle}

\def\TY(#1,#2,#3){#1^{(#2)}_{#3}}

\def\vp{\varphi}

\def\xxi(#1,#2,#3){\displaystyle {}^{#1}\Xi^{(#2)}_{#3}}
\def\xsi(#1,#2,#3){\displaystyle {}^{#1}\Sigma^{(#2)}_{#3}}
\def\xE(#1,#2,#3){\displaystyle {}^{#1}E_{#2}[#3]}
\def\xF(#1,#2){\displaystyle {}^{#1}F_{#2}}
\def\xx(#1,#2){\displaystyle {}^{#1}\Xi_{#2}}
\def\W1{W(\varpi_1)}

\def\ZZ{\mathbb Z}

\def\m@th{\mathsurround=0pt}
\def\fsquare(#1,#2){
\hbox{\vrule$\hskip-0.4pt\vcenter to #1{\normalbaselines\m@th
\hrule\vfil\hbox to #1{\hfill$\scriptstyle #2$\hfill}\vfil\hrule}$\hskip-0.4pt
\vrule}}

\newcommand{\ba}{\begin{array}}
\newcommand{\ea}{\end{array}}

\newcommand{\eq}{\begin{eqnarray}}
\newcommand{\eneq}{\end{eqnarray}}

\title{\textbf{\large{Polyhedral realizations for $B(\infty)$
and extended Young diagrams, Young walls of type
${\rm A}^{(1)}_{n-1}$, ${\rm C}^{(1)}_{n-1}$, ${\rm A}^{(2)}_{2n-2}$, ${\rm D}^{(2)}_{n}$
}}}
\author{\normalsize{YUKI KANAKUBO\thanks{Faculty of Pure and Applied Sciences, University of Tsukuba,
1-1-1 Tennodai, Tsukuba, Ibaraki 305-8577,
Japan: {y-kanakubo@math.tsukuba.ac.jp}.}}
}
\date{}

\begin{document}

\maketitle
\vspace{-10pt}

\begin{abstract}
The crystal bases are quite useful combinatorial tools to study the
representations of quantized universal enveloping algebras $U_q(\mathfrak{g})$. The polyhedral realization for $B(\infty)$
is a combinatorial description of the crystal base, which is defined as an image of embedding $\Psi_{\iota}:B(\infty)\hookrightarrow \mathbb{Z}^{\infty}_{\iota}$,
where $\iota$ is an infinite sequence of indices and
$\mathbb{Z}^{\infty}_{\iota}$ is an infinite $\mathbb{Z}$-lattice with a crystal structure associated with $\iota$.
It is a natural problem to find an explicit form of the polyhedral realization ${\rm Im}(\Psi_{\iota})$.
In this article, supposing that $\mathfrak{g}$ is of affine type ${\rm A}^{(1)}_{n-1}$, ${\rm C}^{(1)}_{n-1}$, ${\rm A}^{(2)}_{2n-2}$ or ${\rm D}^{(2)}_{n}$
and $\iota$ satisfies the condition of `adaptedness', we describe ${\rm Im}(\Psi_{\iota})$ by using several combinatorial objects
such as extended Young diagrams and Young walls.
\end{abstract}

\section{Introduction}

The combinatorics in representation theory of quantized universal enveloping algebra
$U_q(\mathfrak{g})$ has been developed by a lot of authors.
The crystal bases invented in \cite{K0,L} are significant combinatorial tools to know
the skeleton structures of representations of $U_q(\mathfrak{g})$.
It is well-known that the crystal bases are described by
a bunch of objects, like as Young tableaux, LS paths, monomials and so on.

In \cite{Ha, JMMO, KMM},
for an affine Lie algebra $\mathfrak{g}$ of type ${\rm A}^{(1)}_n$, ${\rm C}^{(1)}_n$, ${\rm A}^{(2)}_{2n}$ or ${\rm D}^{(2)}_{n+1}$
and almost all fundamental weights $\Lambda$,
the irreducible integrable highest weight modules $V(\Lambda)$
are realized as Fock space representations
by using `extended Young diagrams'. 
For ${\rm A}^{(1)}_n$ and ${\rm C}^{(1)}_n$ cases,
the crystal base of $B(\Lambda)$ is expressed by extended Young diagrams \cite{JMMO, MM, P}.
In \cite{Kang}, for several affine Lie algebras $\mathfrak{g}$ and level $1$ representations $V(\lambda)$ of $U_q(\mathfrak{g})$, the crystal bases $B(\lambda)$
are realized as sets of reduced proper Young walls.

In \cite{NZ}, the polyhedral realization of crystal base $B(\infty)$ for the negative part 
$U^-_q(\mathfrak{g})$ was introduced as an image of  `Kashiwara embedding'
$\Psi_{\iota}:B(\infty)\hookrightarrow \mathbb{Z}^{\infty}_{\iota}$ associated with an infinite sequence $\iota$ of indices $I$.
Here, $\mathbb{Z}^{\infty}_{\iota}=\{(\cd,a_k,\cd,a_2,a_1)| a_k\in\ZZ
\,\,{\rm and}\,\,a_k=0\,\,{\rm for}\,\,k\gg 0\}$ has a crystal structure associated with $\iota$.
If $\iota$ satisfies the `positivity condition', an algorithm computing an explicit form of the image ${\rm Im}(\Psi_{\iota})$
is given. It is a natural problem to express ${\rm Im}(\Psi_{\iota})$ explicitly. 
In \cite{H1, H2, KS, NZ}, explicit forms of inequalities defining ${\rm Im}(\Psi_{\iota})$ are given in the case $\mathfrak{g}$ is a finite dimensional
simple Lie algebra or classical affine Lie algebra and $\iota=(\cdots,n,\cdots,2,1,n,\cdots,2,1)$.
In the case $\mathfrak{g}$ is a finite dimensional simple Lie algebra and $\iota=(\cdots,i_{N+1},i_N,\cdots,i_2,i_1)$ is a sequence
such that
$(i_N,\cdots,i_2,i_1)$ is a reduced word of the longest element in the Weyl group $W$,
the polyhedral realization ${\rm Im}(\Psi_{\iota})$ coincides with the set of integer points in the string cone of \cite{Lit} associated to the reduced word 
$(i_1,i_2\cdots,i_N)$, which is a polyhedral convex cone. In \cite{GKS16}, a combinatorial expression of string cone via rhombus tiling tools is given.

In \cite{KaN}, assuming that
the sequences $\iota$ satisfies a condition called `adaptedness' (Definition \ref{adapt}), 
we found explicit forms of inequalities defining polyhedral realizations
${\rm Im}(\Psi_{\iota})$ in terms of column tableaux in the case $\mathfrak{g}$ is a finite
dimensional simple Lie algebra of type ${\rm A}_n$, ${\rm B}_n$, ${\rm C}_n$ or ${\rm D}_n$.
It is well known that the crystal base $B(\Lambda_i)$ of fundamental representation $V(\Lambda_i)$ for a classical Lie algebra is
described in terms of column tableaux \cite{KN}.
Therefore, we can expect inequalities defining ${\rm Im}(\Psi_{\iota})$ are expressed as some
combinatorial objects which describe fundamental representations or their crystal bases $B(\Lambda_i)$ in the case $\mathfrak{g}$
is a Kac-Moody algebra other than classical Lie algebras.

In this article, assuming $\iota$ is adapted and $\mathfrak{g}$ is of type ${\rm A}^{(1)}_{n-1}$, ${\rm C}^{(1)}_{n-1}$, ${\rm A}^{(2)}_{2n-2}$ or ${\rm D}^{(2)}_{n}$,
we describe inequalities defining polyhedral realizations
in terms of extended Young diagrams and Young walls.
More precisely, we will define the set of inequalities $\Xi'_{\iota}$ in (\ref{xiiodef}), that is,
it holds ${\rm Im}(\Psi_{\iota})=\{\textbf{a}\in\mathbb{Z}^{\infty} | \varphi(\textbf{a})\geq0\ \text{for any }\varphi\in\Xi'_{\iota}\}$.
The set $\Xi'_{\iota}$ is naturally decomposed as $\Xi'_{\iota}=\bigcup_{k\in I, s\in\mathbb{Z}_{\geq1}} \Xi'_{s,k,\iota}$ with certain subsets $\Xi'_{s,k,\iota}$
(Sect.6). For any $s\in\mathbb{Z}_{\geq1}$, the set $\Xi'_{s,k,\iota}$ is described by using the following objects:
\begin{table}[H]
  \begin{tabular}{|c|c|} \hline
  Type of $\mathfrak{g}$ & combinatorial objects describing $\Xi'_{s,k,\iota}$ \\ \hline
  ${\rm A}^{(1)}_{n-1}$ & extended Young diagram $T$ with $y_{\infty}=k$ \\ 
   &  (the assignment $\ovl{L}_{s,k,\iota}(T)\in \Xi'_{s,k,\iota}$ is related to ${\rm A}^{(1)}_{n-1}$) \\ \hline
  ${\rm D}^{(2)}_{n}$ & extended Young diagram $T$ with $y_{\infty}=k$ \\
   & (the assignment $L_{s,k,\iota}(T)\in \Xi'_{s,k,\iota}$ is related to ${\rm C}^{(1)}_{n-1}$) \\ \hline
  ${\rm A}^{(2)}_{2n-2}$ & revised extended Young diagrams in ${\rm REYD}_{{\rm A}^{(2)},k}$ if $k\in I\setminus\{1\}$\\ 
   & proper Young walls of ground state $Y_{\Lambda_1}$ of type ${\rm A}^{(2)}_{2n-2}$ if $k=1$\\ \hline
  ${\rm C}^{(1)}_{n-1}$ & revised extended Young diagrams in ${\rm REYD}_{{\rm D}^{(2)},k}$ if $k\in I\setminus\{1,n\}$\\ 
   & proper Young walls of ground state $Y_{\Lambda_k}$ of type ${\rm D}^{(2)}_{n}$ if $k=1$ or $n$\\ \hline
  \end{tabular}
\end{table}
\nd
In this way, the set $\Xi'_{s,k,\iota}$ of inequalities is described by a combinatorial object related to the representation $V(\Lambda_k)$
of $U_q(^L\mathfrak{g})$. Here $^L\mathfrak{g}$ is the affine Lie algebra whose generalized Cartan matrix is the transposed matrix of that of
$\mathfrak{g}$.

The organization of this article is as follows. In Sect.2, after a concise reminder on 
crystals, we review the crystals and polyhedral realizations. We also give an algorithm to compute
the polyhedral realizations for $B(\infty)$, which is a slight modification of Nakashima-Zelevinsky's algorithm in \cite{NZ}.
In Sect.3, we recall combinatorial objects such as extended Young diagrams and Young walls. To state the main results,
one also introduce `revised' extended Young diagrams. Sect.4 is devoted to present our main results. We will
express the polyhedral realizations in terms of the combinatorial objects in Sect.3. 
In Sect.5, we prove a closedness of the combinatorial objects
under the action of operators $S'_j$ defined in the modification of Nakashima-Zelevinsky's algorithm.
We completed the proof of main results in Sect.6.

\vspace{2mm}

\nd
{\bf Acknowledgements}
The author wishes to thank Daisuke Sagaki for useful discussions.
This work was supported by JSPS KAKENHI Grant Number JP20J00186.

\section{Polyhedral realizations of $B(\infty)$ and an algorithm}

\subsection{Notation}

Let $\mathfrak{g}$ be a symmetrizable Kac-Moody algebra over $\mathbb{Q}$ with the index set $I=\{1,2,\cdots,n\}$ and
a generalized Cartan matrix $A=(a_{i,j})_{i,j\in I}$. Let  $\mathfrak{h}$ be a Cartan subalgebra, $P\subset \mathfrak{h}^*$ a weight lattice,
$\{\alpha_i\}_{i\in I}$ a set of simple roots and $\{h_i\}_{i\in I}$ a set of simple coroots. 
Let
$\langle \cdot,\cdot \rangle : \mathfrak{h} \times \mathfrak{h}^*\rightarrow \mathbb{Q}$
be the canonical pairing, $P^*:=\{h\in\mathfrak{h} | \langle h,P \rangle\subset\mathbb{Z}\}$ and
$P^+:=\{\lambda\in P | \langle h_i,\lambda \rangle \in\mathbb{Z}_{\geq0} \text{ for all }i\in I \}$.
In particular, it holds $\langle h_{i},\alpha_j \rangle=a_{i,j}$. For each $i\in I$, the fundamental weight $\Lambda_i\in P^+$ is defined as
$\langle h_j,\Lambda_i \rangle=\delta_{i,j}$.
The quantized universal enveloping algebra $U_q(\mathfrak{g})$ is an associative $\mathbb{Q}(q)$-algebra
with generators $e_i$, $f_i$ ($i\in I$) and $q^h$ ($h\in P^*$) satisfying the usual relations. Let $U_q^-(\mathfrak{g})$ be
the subalgebra of $U_q(\mathfrak{g})$ generated by $f_i$ ($i\in I$).

It is known that the irreducible integrable highest weight module $V(\lambda)$ of $U_q(\mathfrak{g})$
has a crystal base $(L(\lambda),B(\lambda))$ for $\lambda\in P^+$. The algebra $U_q^-(\mathfrak{g})$ has
a crystal base $(L(\infty),B(\infty))$ (see \cite{K0,K1}).

\subsection{Crystals}

Let us review the definition of crystals following \cite{K3}:

\begin{defn}
A {\it crystal} is a set $\mathcal{B}$ together with the maps
${\rm wt}:\mathcal{B}\rightarrow P$,
$\varepsilon_i,\varphi_i:\mathcal{B}\rightarrow \mathbb{Z}\sqcup \{-\infty\}$
and $\tilde{e}_i$,$\tilde{f}_i:\mathcal{B}\rightarrow \mathcal{B}\sqcup\{0\}$
($i\in I$) which satisfy the following relation: For $b,b'\in\mathcal{B}$ and $i\in I$,
\begin{itemize}
\item $\varphi_i(b)=\varepsilon_i(b)+\langle h_i,{\rm wt}(b)\rangle$,
\item ${\rm wt}(\tilde{e}_ib)={\rm wt}(b)+\alpha_i$ if $\tilde{e}_i(b)\in\mathcal{B}$,
\quad ${\rm wt}(\tilde{f}_ib)={\rm wt}(b)-\alpha_i$ if $\tilde{f}_i(b)\in\mathcal{B}$,
\item $\varepsilon_i(\tilde{e}_i(b))=\varepsilon_i(b)-1,\ \ 
\varphi_i(\tilde{e}_i(b))=\varphi_i(b)+1$\ if $\tilde{e}_i(b)\in\mathcal{B}$, 
\item $\varepsilon_i(\tilde{f}_i(b))=\varepsilon_i(b)+1,\ \ 
\varphi_i(\tilde{f}_i(b))=\varphi_i(b)-1$\ if $\tilde{f}_i(b)\in\mathcal{B}$, 
\item $\tilde{f}_i(b)=b'$ if and only if $b=\tilde{e}_i(b')$,
\item if $\varphi_i(b)=-\infty$ then $\tilde{e}_i(b)=\tilde{f}_i(b)=0$.
\end{itemize}
Here, $0$ and $-\infty$ are additional elements which do not belong to $\mathcal{B}$ and $\mathbb{Z}$, respectively.
We call $\tilde{e}_i$,$\tilde{f}_i$ {\it Kashiwara operators}.
\end{defn}
\nd
Crystal bases $B(\infty)$ and $B(\lambda)$ are important crystals.

\begin{defn}
Let $\mathcal{B}_1$, $\mathcal{B}_2$ be crystals.
A map $\psi:\mathcal{B}_1\sqcup\{0\}\rightarrow\mathcal{B}_2\sqcup\{0\}$
satisfying the following conditions is said to be a {\it strict morphism} from $\mathcal{B}_1$
to $\mathcal{B}_2$:
\begin{itemize}
\item $\psi(0)=0$,
\item For $i\in I$, if $b\in\mathcal{B}_1$ and $\psi(b)\in \mathcal{B}_2$ then
\[
\psi(0)=0, \quad
{\rm wt}(\psi(b))={\rm wt}(b),\quad
\varepsilon_i(\psi(b))=\varepsilon_i(b),\quad
\varphi_i(\psi(b))=\varphi_i(b),
\]
\item $\tilde{e}_i(\psi(b))=\psi(\tilde{e}_i(b))$ and $\tilde{f}_i(\psi(b))=\psi(\tilde{f}_i(b))$ for $i\in I$ and $b\in\mathcal{B}_1$, where $\eit(0)=\fit(0)=0$.
\end{itemize}
An injective strict morphism $\psi:\mathcal{B}_1\sqcup\{0\}\rightarrow\mathcal{B}_2\sqcup\{0\}$ is said to be {\it strict embedding} of crystals
and denoted by $\psi:\mathcal{B}_1 \hookrightarrow \mathcal{B}_2$.
\end{defn}

\subsection{Polyhedral realizations of $B(\infty)$}\label{2-a}

We define
\[
\ZZ^{\ify}
:=\{(\cd,a_k,\cd,a_2,a_1)| a_k\in\ZZ
\,\,{\rm and}\,\,a_k=0\,\,{\rm for}\,\,k\gg 0\}
\]
and take an infinite sequence
$\io=(\cd,i_k,\cd,i_2,i_1)$ of indices from $I$ such that
\begin{equation}
{\hbox{
$i_k\ne i_{k+1}$ and $\sharp\{k: i_k=i\}=\ify$ for any $i\in I$.}}
\label{seq-con}
\end{equation}
Following \cite{NZ}, one can define a crystal structure on $\ZZ^{\ify}$ 
corresponding to $\iota$ as follows:
For $\textbf{a}=(\cd,a_k,\cd,a_2,a_1)\in\ZZ^{\ify}$ and $k\in\mathbb{Z}_{\geq1}$, we set
\begin{equation}\label{sigmak}
\sigma_k(\textbf{a}):=a_k+\sum_{j\in\mathbb{Z} : j>k} \langle h_{i_k},\alpha_{i_j}\rangle a_j.
\end{equation}
Because $a_j=0$ if $j$ is sufficiently larger than $0$, the above definition is well-defined.
For $i\in I$, we define ${\rm wt} : \ZZ^{\ify}\rightarrow P$, $\varepsilon_i:\ZZ^{\ify}\rightarrow\mathbb{Z}$
and $\varphi_i:\ZZ^{\ify}\rightarrow\mathbb{Z}$ as
\[
{\rm wt}(\textbf{a}):=-\sum^{\infty}_{j=1} a_j\alpha_{i_j},\quad
\varepsilon_i(\textbf{a}):={\rm max}\{\sigma_k(\textbf{a}) | k\in\mathbb{Z}_{\geq1},\ i_k=i \},\quad
\varphi_i(\textbf{a}):=
\langle h_i, {\rm wt}(\textbf{a}) \rangle
+\varepsilon_i(\textbf{a}).
\]
Putting
\[
M^{(i)}=M^{(i)}(\textbf{a}):=\{k\in \mathbb{Z}_{\geq1} | i_k=i,\ \sigma_k(\textbf{a})=\varepsilon_i(\textbf{a})\},
\]
we see that $M^{(i)}$ is a finite set if and only if $\varepsilon_i(\textbf{a})>0$.
One defines $\tilde{f}_i:\ZZ^{\ify}\rightarrow\ZZ^{\ify}$ and
$\tilde{e}_i:\ZZ^{\ify}\rightarrow\ZZ^{\ify}\sqcup\{0\}$ as
\begin{equation}\label{kashi-def1}
(\tilde{f}_i(\textbf{a}))_k:=a_k+\delta_{k,{\rm min}M^{(i)}},
\end{equation}
\begin{equation}\label{kashi-def2}
(\tilde{e}_i(\textbf{a}))_k:=a_k-\delta_{k,{\rm max}M^{(i)}} \ \text{if } \varepsilon_i(\textbf{a})>0
\end{equation}
and $\tilde{e}_i(\textbf{a})=0$ if $\varepsilon_i(\textbf{a})=0$.

\begin{thm}\cite{NZ}
The set $\ZZ^{\ify}$ together with the above maps
$\{\tilde{e}_i\}_{i\in I},\{\tilde{f}_i\}_{i\in I},\{\varepsilon_i\}_{i\in I},\{\varphi_i\}_{i\in I}$ and ${\rm wt}$
is a crystal. 
\end{thm}
\nd
Let $\ZZ^{\ify}_{\io}$ denote the above crystal.

\begin{thm}\cite{K3,NZ}
\label{emb}
There is a unique strict embedding of crystals
\begin{equation}
\Psi_{\io}:B(\ify)\hookrightarrow \ZZ^{\ify}_{\geq 0}
\subset \ZZ^{\ify}_{\io},
\label{psi}
\end{equation}
such that 
$\Psi_{\io} (u_{\ify}) = \textbf{0}$, where
$u_{\ify}\in B(\ify)$ is
the highest weight vector, $\textbf{0}:=(\cd,0,\cd,0,0)$
and $\ZZ^{\ify}_{\geq 0}:=\{(\cdots,a_k,\cdots,a_2,a_1)\in\mathbb{Z}^{\infty} | a_k\geq0\}$.
\end{thm}

\begin{defn}
The image ${\rm Im}(\Psi_{\iota}) (\cong B(\infty))$ is called a \it{polyhedral realization} of $B(\infty)$.
\end{defn}

\subsection{Modified Nakashima-Zelevinsky's algorithm}\label{poly-uqm}

In \cite{NZ}, an algorithm to describe polyhedral realizations ${\rm Im}(\Psi_{\iota})$ is given.
In this subsection, we introduce a modified algorithm of it. Let us fix an infinite sequence 
\[
\io=(\cd,i_k,\cd,i_2,i_1)
\]
of indices
satisfying (\ref{seq-con}).
We consider a vector space
$\QQ^{\ify}:=\{\textbf{a}=
(\cd,a_k,\cd,a_2,a_1)| a_k \in \QQ\,\,{\rm and }\,\,
a_k = 0\,\,{\rm for}\,\, k \gg 0\}$ and its dual space $(\QQ^{\ify})^*$.
For $k\in \mathbb{Z}_{\geq1}$, the element $x_k\in(\QQ^{\ify})^*$ is defined as
$x_k(\cd,a_k,\cd,a_2,a_1)=a_k$. For $k\in \mathbb{Z}_{<1}$, we set $x_k:=0$.
Using this notation, we will write each element 
$\varphi\in(\QQ^{\ify})^*$ as $\varphi=\sum_{k\in\mathbb{Z}_{\geq1}}c_kx_k$ with some $c_k\in\mathbb{Q}$.
For $k\in \mathbb{Z}_{\geq1}$, let
\[
k^{(+)}:={\rm min}\{l\in\mathbb{Z}_{\geq1}\ |\ l>k\,\,{\rm and }\,\,i_k=i_l\},\ \ 
k^{(-)}:={\rm max} \{l\in\mathbb{Z}_{\geq1}\ |\ l<k\,\,{\rm and }\,\,i_k=i_l\}\cup\{0\},
\]
\begin{equation}
\beta_k:=x_k+\sum_{k<j<\kp}\lan h_{i_k},\al_{i_j}\ran x_j+x_{\kp}\in (\QQ^{\ify})^*,
\label{betak}
\end{equation}
and $\beta_0:=0$. Note that it holds $\beta_k(\textbf{a})=\sigma_k(\textbf{a})-\sigma_{k^{(+)}}(\textbf{a})$
by (\ref{sigmak}) for $\textbf{a}\in\mathbb{Z}^{\infty}$.
We define the operator 
$S_k'=S_{k,\io}':(\QQ^{\ify})^*\rightarrow (\QQ^{\ify})^*$ as follows:
For $\vp=\sum_{k\in\mathbb{Z}_{\geq1}}c_kx_k\in(\QQ^{\ify})^*$,
\begin{equation}
S_k'(\vp):=
\begin{cases}
\vp-\beta_k & {\mbox{ if }}\ c_k>0,\\
 \vp+\beta_{\km} & {\mbox{ if }}\ c_k< 0,\\
 \vp &  {\mbox{ if }}\ c_k= 0.
\end{cases}
\label{Sk}
\end{equation}
We often write $S_k'(\vp)$ as $S_k'\vp$. Let us define
\begin{eqnarray}
\Xi_{\io}' &:=  &\{S_{j_l}'\cd S_{j_2}'S_{j_1}'x_{j_0}\,|\,
l\in\mathbb{Z}_{\geq0},j_0,j_1,\cd,j_l\in\mathbb{Z}_{\geq1}\}, \label{xiiodef}\\
\Sigma_{\io}' & := &
\{\textbf{a}\in \ZZ^{\ify}\subset \QQ^{\ify}\,|\,\vp(\textbf{a})\geq0\,\,{\rm for}\,\,
{\rm any}\,\,\vp\in \Xi_{\io}'\}.
\end{eqnarray}
We say $\io$ satisfies the $\Xi'$-{\it positivity condition}
when it holds
\begin{equation}
{\hbox{if $\km=0$ then $c_k\geq0$ for any 
$\vp=\sum_k c_kx_k\in \Xi_{\io}'$}}.
\label{posi}
\end{equation}
\begin{thm}\label{polyhthm}
Let $\io$ be a sequence of indices satisfying $(\ref{seq-con})$ 
and $(\ref{posi})$. Then it holds 
${\rm Im}(\Psi_{\io})=\Sigma_{\io}'$.
\end{thm}

\nd
{\it Proof.} We prove our claim by a slight modification of the proof of Theorem 3.1 in \cite{NZ}.
It follows from Theorem \ref{emb} that ${\rm Im}(\Psi_{\io})\subset \mathbb{Z}^{\infty}_{\geq0}$
is a subset of $\mathbb{Z}^{\infty}_{\iota}$ obtained from $\Psi_{\iota}(u_{\infty})=\textbf{0}$
by applying Kashiwara operators $\tilde{f}_i$ ($i\in I$). By the definition of $\Sigma_{\io}'$,
it holds $\textbf{0}=(\cdots,0,\cdots,0,0)\in\Sigma_{\io}'$. Thus, to see
${\rm Im}(\Psi_{\io})\subset\Sigma_{\io}'$, we need to check $\Sigma_{\io}'$ is closed
under the action of all $\tilde{f}_i$. Let us take $\textbf{a}=(\cdots,a_k,\cdots,a_2,a_1)\in\Sigma_{\io}'$ and $i\in I$
and show $\tilde{f}_i(\textbf{a})\in\Sigma_{\io}'$, that is, for any $\vp=\sum_{l\in\mathbb{Z}_{\geq1}}c_lx_l \in\Xi'_{\iota}$,
show
\begin{equation}\label{pr-1}
\vp(\tilde{f}_i(\textbf{a}))\geq0.
\end{equation}
It follows by (\ref{betak}), (\ref{Sk}) that
all coefficients $c_l$ of $\vp$ are integers.
The element $\tilde{f}_i(\textbf{a})$ is in the form
\[
\tilde{f}_i(\textbf{a})=(\cdots,a_{k+1},a_k+1,a_{k-1},\cdots,a_2,a_1)
\]
with some $k\in \mathbb{Z}_{\geq1}$ such that $i_k=i$ by (\ref{kashi-def1}). 
Since we know $\varphi(\textbf{a})\geq0$, it holds $\vp(\tilde{f}_i(\textbf{a}))=\varphi(\textbf{a})+c_k\geq c_k$
so that we may suppose that $c_k<0$. Our assumption (\ref{posi}) means $k^{(-)}\geq1$.
Considering (\ref{kashi-def1}), it holds $\sigma_k(\textbf{a})>\sigma_{k^{(-)}}(\textbf{a})$.
Hence, it follows $\beta_{k^{(-)}}(\textbf{a})=\sigma_{k^{(-)}}(\textbf{a})-\sigma_k(\textbf{a})\leq-1$.
Thus, we see that
\[
\vp(\tilde{f}_i(\textbf{a}))=\varphi(\textbf{a})+c_k
\geq \varphi(\textbf{a})-c_k\beta_{k^{(-)}}(\textbf{a})
= S_k'^{|c_k|}(\vp)(\textbf{a})\geq0
\]
by $S_k'^{|c_k|}(\vp)\in \Xi'_{\iota}$ and  $\textbf{a}\in\Sigma_{\io}'$, which implies (\ref{pr-1}) so that ${\rm Im}(\Psi_{\io})\subset\Sigma_{\io}'$ holds.

Next, we prove $\Sigma_{\io}'\subset{\rm Im}(\Psi_{\io})$. First, we show $\tilde{e}_i\Sigma_{\io}'\subset
\Sigma_{\io}'\cup\{0\}$ for any $i\in I$. For any element $\textbf{a}=(\cdots,a_3,a_2,a_1)\in \Sigma_{\io}'$,
if $\tilde{e}_i(\textbf{a})\neq0$ then
$\tilde{e}_i(\textbf{a})$ is in the form
\[
\tilde{e}_i(\textbf{a})=(\cdots,a_{k+1},a_k-1,a_{k-1},\cdots,a_2,a_1)
\]
with some $k\in \mathbb{Z}_{\geq1}$ such that $i_k=i$ by (\ref{kashi-def2}). We need to prove $\vp(\tilde{e}_i(\textbf{a}))\geq0$
for any $\vp=\sum_{l\in\mathbb{Z}_{\geq1}}c_lx_l \in\Xi'_{\iota}$.
By $\vp(\tilde{e}_i(\textbf{a}))=\vp(\textbf{a})-c_k\geq -c_k$, we may assume $c_k>0$.
It follows from (\ref{kashi-def2}) that $\sigma_{k}(\textbf{a})>\sigma_{k^{(+)}}(\textbf{a})$ 
so that $\beta_k(\textbf{a})=\sigma_{k}(\textbf{a})-\sigma_{k^{(+)}}(\textbf{a})\geq1$. Thus,
one can verify
\[
\vp(\tilde{e}_i(\textbf{a}))=\vp(\textbf{a})-c_k
\geq \vp(\textbf{a})-c_k\beta_k(\textbf{a})
=S_k'^{c_k}(\vp)(\textbf{a})\geq0
\]
since $S_k'^{c_k}(\vp)\in\Xi'_{\iota}$ and  $\textbf{a}\in\Sigma_{\io}'$. Thus, we obtain $\tilde{e}_i\Sigma_{\io}'\subset
\Sigma_{\io}'\cup\{0\}$.

For any $\textbf{a}=(\cdots,a_3,a_2,a_1)\in \Sigma_{\io}'\setminus\{\textbf{0}\}\subset\mathbb{Z}^{\infty}_{\geq0}$,
there exists $i\in I$ such that $\tilde{e}_{i}\textbf{a}\neq0$. In fact,
putting $j:={\rm max}\{l\in\mathbb{Z}_{\geq1} | a_l>0 \}$ and $i:=i_j$, it holds $\tilde{e}_{i}\textbf{a}\neq0$.
Taking $\tilde{e}_i\Sigma_{\io}'\subset
\Sigma_{\io}'\cup\{0\}$ and $\Sigma_{\io}'\subset\mathbb{Z}^{\infty}_{\geq0}$
into account, there exists a sequence $i_1,\cdots,i_l\in I$ such that
\[
\tilde{e}_{i_l}\cdots\tilde{e}_{i_2}\tilde{e}_{i_1}\textbf{a}=\textbf{0},
\] 
which yields $\textbf{a}=\tilde{f}_{i_1}\tilde{f}_{i_2}\cdots\tilde{f}_{i_l}\textbf{0}\in {\rm Im}(\Psi_{\io})$.
Therefore, it follows $\Sigma_{\io}'\subset{\rm Im}(\Psi_{\io})$. \qed

\begin{rem}
In \cite{NZ},
piecewise-linear operators $S_k$ on $(\QQ^{\ify})^*$
are defined for $k\in\mathbb{Z}_{\geq1}$.
Using these operators, one can define a set
\[
\Xi_{\io}:= \{S_{j_l}\cd S_{j_2}S_{j_1}x_{j_0} |\,
l\in\mathbb{Z}_{\geq0},j_0,j_1,\cd,j_l\in\mathbb{Z}_{\geq1}\}.
\]
We can easily verify that
if $c_k\in\mathbb{Z}_{\neq0}$ then $S_k'^{|c_k|}(\vp):=
\underbrace{S_k'\cdots S_k'}_{|c_k|\text{times}}(\vp)=S_k(\vp)$
and if $c_k=0$ then $S_k'(\vp)=S_k(\vp)=\vp$.
Hence, it follows $\Xi_{\io}\subset \Xi'_{\io}$. 
Therefore, if $\io$ satisfies the $\Xi'$-positivity condition
then the {\it positivity assumption} in \cite{NZ} is satisfied.
\end{rem}

\section{Extended Young diagrams and Young walls}

In what follows, we consider the case $\mathfrak{g}$ is of affine type ${\rm A}_{n-1}^{(1)}$, ${\rm C}_{n-1}^{(1)}$, ${\rm A}_{2n-2}^{(2)}$ or ${\rm D}_{n}^{(2)}$.
The numbering of vertices in affine Dynkin diagrams are as follows:
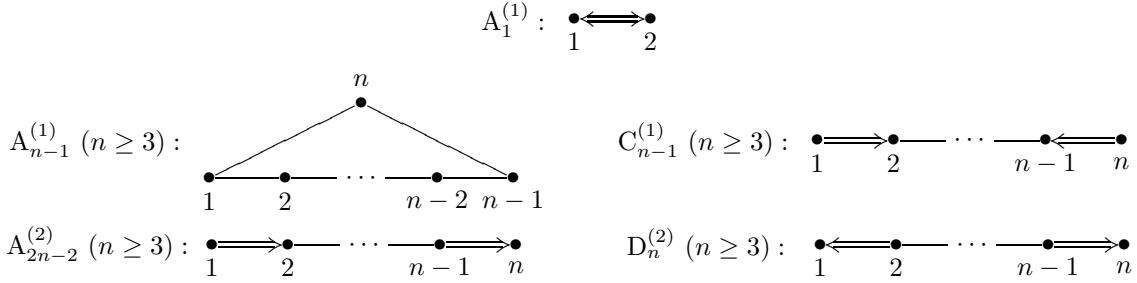
\[
\begin{xy}
(-8,0) *{{\rm A}_{1}^{(1)} : }="A1",
(0,0) *{\bullet}="1",
(0,-3) *{1}="1a",
(10,0)*{\bullet}="2",
(10,-3) *{2}="2a",
\ar@{<=>} "1";"2"^{}
\end{xy}
\]
\[
\begin{xy}
(-15,5) *{{\rm A}_{n-1}^{(1)} \ (n\geq 3) : }="A1",
(20,10) *{\bullet}="n",
(20,13) *{n}="na",
(0,0) *{\bullet}="1",
(0,-3) *{1}="1a",
(10,0)*{\bullet}="2",
(10,-3)*{2}="2a",
(20,0)*{\ \cdots\ }="3",
(30,0)*{\bullet}="4",
(30,-3)*{n-2}="4a",
(40,0)*{\bullet}="5",
(40,-3)*{n-1}="5a",
(65,5) *{{\rm C}_{n-1}^{(1)} \ (n\geq 3) : }="C1",
(80,5) *{\bullet}="11",
(80,2) *{1}="11a",
(90,5)*{\bullet}="22",
(90,2)*{2}="22a",
(100,5)*{\ \cdots\ }="33",
(110,5)*{\bullet}="44",
(110,2)*{n-1}="44a",
(120,5)*{\bullet}="55",
(120,2)*{n}="55a",
\ar@{-} "1";"2"^{}
\ar@{-} "2";"3"^{}
\ar@{-} "3";"4"^{}
\ar@{-} "4";"5"^{}
\ar@{-} "1";"n"^{}
\ar@{-} "n";"5"^{}
\ar@{=>} "11";"22"^{}
\ar@{-} "22";"33"^{}
\ar@{-} "33";"44"^{}
\ar@{<=} "44";"55"^{}
\end{xy}
\]
\[
\begin{xy}
(-15,5) *{{\rm A}_{2n-2}^{(2)}\ (n\geq 3) : }="A2",
(0,5) *{\bullet}="1",
(0,2) *{1}="1a",
(10,5)*{\bullet}="2",
(10,2)*{2}="2a",
(20,5)*{\ \cdots\ }="3",
(30,5)*{\bullet}="4",
(30,2)*{n-1}="4a",
(40,5)*{\bullet}="5",
(40,2)*{n}="5a",
(65,5) *{{\rm D}_{n}^{(2)} \ (n\geq 3) : }="D2",
(80,5) *{\bullet}="11",
(80,2) *{1}="11a",
(90,5)*{\bullet}="22",
(90,2)*{2}="22a",
(100,5)*{\ \cdots\ }="33",
(110,5)*{\bullet}="44",
(110,2)*{n-1}="44a",
(120,5)*{\bullet}="55",
(120,2)*{n}="55a",
\ar@{=>} "1";"2"^{}
\ar@{-} "2";"3"^{}
\ar@{-} "3";"4"^{}
\ar@{=>} "4";"5"^{}
\ar@{<=} "11";"22"^{}
\ar@{-} "22";"33"^{}
\ar@{-} "33";"44"^{}
\ar@{=>} "44";"55"^{}
\end{xy}
\]
Replacing our numberings $1,2,\cdots,n-1,n,n+1$ of ${\rm A}_{n}^{(1)}$, ${\rm C}_{n}^{(1)}$ and ${\rm D}_{n+1}^{(2)}$
with $0,1,2,\cdots,n$, we get the numbering in \cite{JMMO, Kang, KMM}.
Replacing our numberings $1,2,\cdots,n-1,n,n+1$ of ${\rm A}_{2n}^{(2)}$ with $n,n-1,\cdots,1,0$, we obtain the numbering in \cite{Kang, KMM}.

\subsection{Extended Young diagrams}\label{EYD-sub}

\begin{defn}\cite{Ha,JMMO}
For a fixed integer $y_{\infty}$, a sequence $(y_k)_{k\in\mathbb{Z}_{\geq0}}$
is called an {\it extended Young diagram} of charge $y_{\infty}$ if it holds
\begin{itemize}
\item $y_k\in\mathbb{Z}$, $y_k\leq y_{k+1}$ for all $k\in\mathbb{Z}_{\geq0}$,
\item $y_k=y_{\infty}$ for $k\gg0$.
\end{itemize}
\end{defn}

Each extended Young diagram is described as an infinite Young diagram
drawn on $\mathbb{R}_{\geq0}\times \mathbb{R}_{\leq y_{\infty}}$ as follows:
For $(y_k)_{k\in\mathbb{Z}_{\geq0}}$, we draw a line between the points $(k,y_k)$ and $(k+1,y_k)$ and when
$y_k<y_{k+1}$ draw
a line
between $(k+1,y_k)$ and $(k+1,y_{k+1})$
for each $k\in\mathbb{Z}_{\geq0}$.
\begin{ex}\label{ex-1}
Let $T=(y_k)_{k\in\mathbb{Z}_{\geq0}}$ be
the extended Young diagram of charge
$y_{\infty}=1$ defined as
$y_0=-3$, $y_1=-2$, $y_2=y_3=-1$, $y_4=0$, $y_5=1, y_6=1,\cdots$.
Then $Y$ is described as
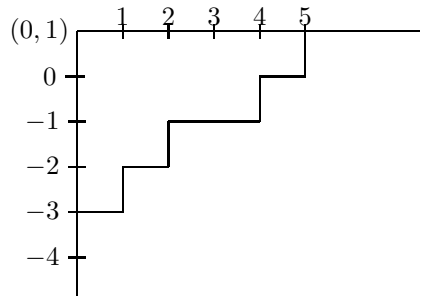
\[
\begin{xy}
(0,0) *{}="1",
(45,0)*{}="2",
(0,-35)*{}="3",
(-5,0)*{(0,1)}="4",
(6,2) *{1}="10",
(6,-1) *{}="1010",
(12,2) *{2}="11",
(12,-1) *{}="1111",
(18,2) *{3}="12",
(18,-1) *{}="1212",
(24,2) *{4}="13",
(24,-1) *{}="1313",
(30,2) *{5}="14",
(30,-1) *{}="1414",
(-3,-6)*{0\ }="5",
(1,-6)*{}="55",
(-4,-12)*{-1\ }="6",
(1,-12)*{}="66",
(-4,-18)*{-2\ }="7",
(1,-18)*{}="77",
(-4,-24)*{-3\ }="8",
(1,-24)*{}="88",
(6,-24)*{}="810",
(6,-18)*{}="8107",
(12,-18)*{}="810711",
(12,-12)*{}="8107116",
(24,-12)*{}="810711613",
(24,-6)*{}="8107116130",
(30,-6)*{}="81071161300",
(30,0)*{}="810711613000",
(-4,-30)*{-4\ }="9",
(1,-30)*{}="99",
\ar@{-} "1";"2"^{}
\ar@{-} "1";"3"^{}
\ar@{-} "5";"55"^{}
\ar@{-} "6";"66"^{}
\ar@{-} "7";"77"^{}
\ar@{-} "8";"810"^{}
\ar@{-} "810";"8107"^{}
\ar@{-} "8107";"810711"^{}
\ar@{-} "810711";"8107116"^{}
\ar@{-} "8107116";"810711613"^{}
\ar@{-} "810711613";"8107116130"^{}
\ar@{-} "8107116130";"81071161300"^{}
\ar@{-} "81071161300";"810711613000"^{}
\ar@{-} "9";"99"^{}
\ar@{-} "10";"1010"^{}
\ar@{-} "11";"1111"^{}
\ar@{-} "12";"1212"^{}
\ar@{-} "13";"1313"^{}
\ar@{-} "14";"1414"^{}
\end{xy}
\]
\end{ex}
\nd
Note that if $y_k< y_{k+1}$ then
the points $(k+1,y_k)$ and $(k+1,y_{k+1})$ are corners. 
\begin{defn}\cite{JMMO}
For an extended Young diagram
$(y_k)_{k\in\mathbb{Z}_{\geq0}}$, if $y_k< y_{k+1}$ then
we say $(k+1,y_k)$ is a {\it convex corner} and $(k+1,y_{k+1})$
is a {\it concave corner}. The point $(0,y_0)$ is also called a concave corner.
A corner $(i,j)$ is called a $d$-{\it diagonal corner}, where $d=i+j$.
\end{defn}

\begin{ex}
Let us consider the same
extended Young diagram as in Example \ref{ex-1}.
The points $(1,-3)$, $(2,-2)$, $(4,-1)$ and $(5,0)$ are convex corners and
$(0,-3)$, $(1,-2)$, $(2,-1)$, $(4,0)$ and $(5,1)$ are concave corners.

\end{ex}

\begin{defn}\cite{JMMO,KMM}\label{pi1-def}
\begin{enumerate}
\item
The map $\overline{\ }:\mathbb{Z}\rightarrow \{1,2,\cdots,n\}=I$ is defined as
\[
\overline{l+rn}=l
\]
for any $r\in\mathbb{Z}$ and $l\in\{1,2,\cdots,n\}$.
\item
We define a map
$\{1,2,\cdots,2n-2\}\rightarrow \{1,2,\cdots,n\}$
as
\[
l\mapsto l,\ 2n-l\mapsto l \qquad (2\leq l\leq n-1),
\]
\[
1\mapsto1,\ n\mapsto n
\]
and extend it to a map $\pi:\mathbb{Z}\rightarrow\{1,2,\cdots,n\}=I$ by periodicity $2n-2$.
\item We define a map $\{1,2,3,\cdots,2n-1\}\rightarrow \{1,2,\cdots,n\}$ as
\[
l\mapsto l,\ \ 2n-l\mapsto l \quad (1\leq l\leq n-1),\ \ n\mapsto n
\]
and extend it to the map $\pi_1:\mathbb{Z}\rightarrow \{1,2,\cdots,n\}=I$ with periodicity $2n-1$.
\item We define a map $\{1,2,3,\cdots,2n\}\rightarrow \{1,2,\cdots,n\}$ as
\[
l\mapsto l,\ \ 2n+1-l\mapsto l \quad (1\leq l\leq n),
\]
and extend it to the map $\pi_2:\mathbb{Z}\rightarrow \{1,2,\cdots,n\}=I$ with periodicity $2n$.
\end{enumerate}
\end{defn}

The maps $\overline{\ }$, $\pi$, $\pi_1$ and $\pi_2$ in Definition \ref{pi1-def} were
introduced to define the action of Chevalley generators (or Kashiwara operators) of type ${\rm A}^{(1)}_n$, ${\rm C}^{(1)}_n$, ${\rm A}^{(2)}_{2n}$
or ${\rm D}^{(2)}_{n+1}$ on extended Young diagrams and each corner $(i,j)$ is colored by
$\ovl{i+j}$, $\pi(i+j)$, $\pi_1(i+j)$ or $\pi_2(i+j)\in I$, respectively \cite{JMMO, KMM, P}.
Roughly speaking, concave corners colored by $m$ are replaced by convex corners by the action of $f_m$ or 
Kashiwara operator $\tilde{f}_m$. 
In Proposition \ref{prop-closednessAD}, \ref{A2closed} and \ref{D2closed},
we will consider a similar replacement for the action of operators $S'$.

\subsection{Revised extended Young diagrams}

To describe inequalities of type ${\rm A}^{(2)}_{2n-2}$ and ${\rm C}^{(1)}_{n-1}$,
we need to introduce `revised extended Young diagrams'.

\begin{defn}\label{AEYD}
For $k\in I\setminus\{1\}$, let ${\rm REYD}_{{\rm A}^{(2)},k}$ be the set of sequences $(y_t)_{t\in\mathbb{Z}}$ such that
\begin{itemize}
\item[$(1)$] $y_t\in\mathbb{Z}$ for $t\in\mathbb{Z}$,  
\item[$(2)$] $y_{t}=k$ for $t\gg0$ and $y_t=k+t$ for $t\ll0$,
\item[$(3)$] for $t\in\mathbb{Z}$ such that $k+t\not\equiv 0$ (mod $2n-1$), it holds either $y_{t+1}=y_t$ or $y_{t+1}=y_t+1$,
\item[$(4)$] for $t\in\mathbb{Z}_{>0}$ such that $k+t\equiv 0$ (mod $2n-1$), it holds $y_{t+1}\geq y_t$,
\item[$(5)$] for $t\in\mathbb{Z}_{<0}$ such that $k+t\equiv 0$ (mod $2n-1$), it holds $y_{t+1}\leq y_t+1$.
\end{itemize}
\end{defn}

\begin{defn}\label{DEYD}
For $k\in I\setminus\{1,n\}$, let ${\rm REYD}_{{\rm D}^{(2)},k}$ be the set of sequences $(y_t)_{t\in\mathbb{Z}}$ such that
\begin{itemize}
\item[$(1)$] $y_t\in\mathbb{Z}$ for $t\in\mathbb{Z}$,  
\item[$(2)$] $y_{t}=k$ for $t\gg0$ and $y_t=k+t$ for $t\ll0$,
\item[$(3)$] for $t\in\mathbb{Z}$ such that $k+t\not\equiv 0, n$ (mod $2n$), it holds either $y_{t+1}=y_t$ or $y_{t+1}=y_t+1$,
\item[$(4)$] for $t\in\mathbb{Z}_{>0}$ such that $k+t\equiv 0$ or $n$ (mod $2n$), it holds $y_{t+1}\geq y_t$,
\item[$(5)$] for $t\in\mathbb{Z}_{<0}$ such that $k+t\equiv 0$ or $n$ (mod $2n$), it holds $y_{t+1}\leq y_t+1$.
\end{itemize}
\end{defn}
Each element in ${\rm REYD}_{{\rm A}^{(2)},k}$ and ${\rm REYD}_{{\rm D}^{(2)},k}$  is described as a diagram drawn on
$\mathbb{R}\times \mathbb{R}_{\leq k}$
by a similar rule to ordinary extended Young diagrams. For example, 
let $n=3$, $k=2$ and $T=(y_t)_{t\in\mathbb{Z}}$ be
the element in ${\rm REYD}_{{\rm A}^{(2)},2}$
defined as
\[
y_{l}=l+2\ (l\leq -3),\ y_{-2}=0,\ y_{-1}=y_0=y_1=-1,\ y_2=y_3=0,\ y_t=2 (t\geq 4).
\]
Then $T$ is described as
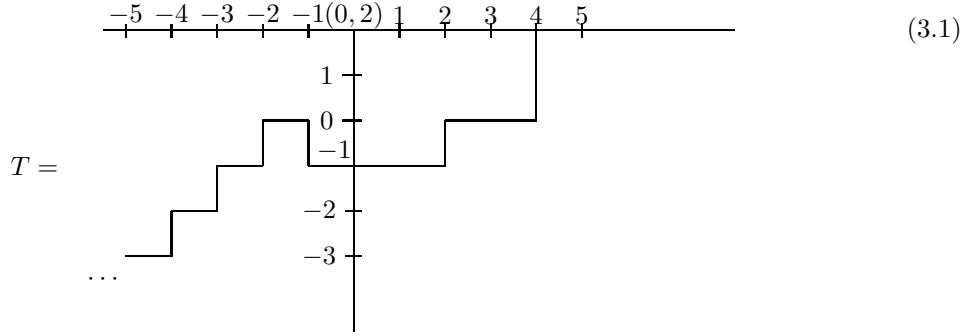
\begin{equation}\label{reydA-ex1}
\begin{xy}
(-42,-18) *{T=}="YY",
(-33,0) *{}="-6",
(-6,2) *{-1}="-1",
(-12,2) *{-2}="-2",
(-18,2) *{-3}="-3",
(-24,2) *{-4}="-4",
(-30,2) *{-5}="-5",
(-6,-1) *{}="-1a",
(-12,-1) *{}="-2a",
(-18,-1) *{}="-3a",
(-24,-1) *{}="-4a",
(-30,-1) *{}="-5a",
(0,0) *{}="1",
(50,0)*{}="2",
(0,-40)*{}="3",
(0,2)*{(0,2)}="4",
(6,2) *{1}="10",
(6,-1) *{}="1010",
(12,2) *{2}="11",
(12,-1) *{}="1111",
(18,2) *{3}="12",
(18,-1) *{}="1212",
(24,2) *{4}="13",
(24,-1) *{}="1313",
(30,2) *{5}="14",
(30,-1) *{}="1414",
(-3,-6)*{1\ }="5",
(1,-6)*{}="55",
(-3,-12)*{0\ }="6",
(-6,-12)*{}="6a",
(-12,-12)*{}="6aa",
(-12,-18)*{}="6aaa",
(-18,-18)*{}="6aaaa",
(-18,-24)*{}="st1",
(-24,-24)*{}="st2",
(-24,-30)*{}="st3",
(-30,-30)*{}="st4",
(-33,-33)*{\cdots}="stdot",
(1,-12)*{}="66",
(-6,-18)*{}="7",
(-2,-16)*{-1\ }="7a",
(1,-18)*{}="77",
(-4,-24)*{-2\ }="8",
(1,-24)*{}="88",
(6,-18)*{}="8107",
(12,-18)*{}="810711",
(12,-12)*{}="8107116",
(24,-12)*{}="810711613",
(24,-6)*{}="8107116130",
(30,-6)*{}="81071161300",
(30,0)*{}="810711613000",
(-4,-30)*{-3\ }="9",
(1,-30)*{}="99",
\ar@{-} "1";"-6"^{}
\ar@{-} "1";"2"^{}
\ar@{-} "1";"3"^{}
\ar@{-} "-1";"-1a"^{}
\ar@{-} "-2";"-2a"^{}
\ar@{-} "-3";"-3a"^{}
\ar@{-} "-4";"-4a"^{}
\ar@{-} "-5";"-5a"^{}
\ar@{-} "5";"55"^{}
\ar@{-} "6";"66"^{}
\ar@{-} "8107";"7"^{}
\ar@{-} "8107";"810711"^{}
\ar@{-} "810711";"8107116"^{}
\ar@{-} "8107116";"810711613"^{}
\ar@{-} "810711613";"8107116130"^{}
\ar@{-} "8107116130";"13"^{}
\ar@{-} "7";"6a"^{}
\ar@{-} "6aa";"6a"^{}
\ar@{-} "6aaa";"6aa"^{}
\ar@{-} "6aaaa";"6aaa"^{}
\ar@{-} "st1";"6aaaa"^{}
\ar@{-} "st1";"st2"^{}
\ar@{-} "st2";"st3"^{}
\ar@{-} "st3";"st4"^{}
\ar@{-} "8";"88"^{}
\ar@{-} "9";"99"^{}
\ar@{-} "10";"1010"^{}
\ar@{-} "11";"1111"^{}
\ar@{-} "12";"1212"^{}
\ar@{-} "13";"1313"^{}
\ar@{-} "14";"1414"^{}
\end{xy}
\end{equation}

\begin{rem}\label{reydrem}
Each revised extended Young diagram can be drawn in
the quarter plane $\mathbb{R}_{\geq0}\times \mathbb{R}_{\leq k}$
just as in ordinary extended Young diagrams:
Shifting $j$-th row of the revised extended Young diagram to the right by $j-1$ for all $j\in\mathbb{Z}_{\geq1}$,
it becomes a diagram drawn in $\mathbb{R}_{\geq0}\times \mathbb{R}_{\leq k}$. For example,
the diagram $T$ in (\ref{reydA-ex1}) is changed as follows:
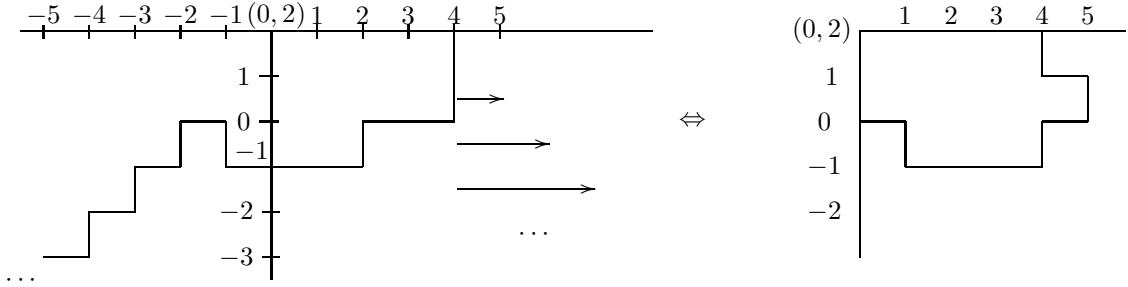
\begin{equation*}
\begin{xy}
(-33,0) *{}="-6",
(-6,2) *{-1}="-1",
(-12,2) *{-2}="-2",
(-18,2) *{-3}="-3",
(-24,2) *{-4}="-4",
(-30,2) *{-5}="-5",
(-6,-1) *{}="-1a",
(-12,-1) *{}="-2a",
(-18,-1) *{}="-3a",
(-24,-1) *{}="-4a",
(-30,-1) *{}="-5a",
(0,0) *{}="1",
(50,0)*{}="2",
(0,-33)*{}="3",
(0,2)*{\ (0,2)}="4",
(6,2) *{1}="10",
(6,-1) *{}="1010",
(12,2) *{2}="11",
(12,-1) *{}="1111",
(18,2) *{3}="12",
(18,-1) *{}="1212",
(24,2) *{4}="13",
(24.5,-9) *{}="arrow1",
(30.5,-9) *{}="arrow1-1",
(24.5,-15) *{}="arrow2",
(36.5,-15) *{}="arrow2-1",
(24.5,-21) *{}="arrow3",
(42.5,-21) *{}="arrow3-1",
(34.5,-27) *{\cdots}="arrow4",
(24,-1) *{}="1313",
(30,2) *{5}="14",
(30,-1) *{}="1414",
(-3,-6)*{1\ }="5",
(1,-6)*{}="55",
(-3,-12)*{0\ }="6",
(-6,-12)*{}="6a",
(-12,-12)*{}="6aa",
(-12,-18)*{}="6aaa",
(-18,-18)*{}="6aaaa",
(-18,-24)*{}="st1",
(-24,-24)*{}="st2",
(-24,-30)*{}="st3",
(-30,-30)*{}="st4",
(-33,-33)*{\cdots}="stdot",
(1,-12)*{}="66",
(-6,-18)*{}="7",
(-2,-16)*{-1\ }="7a",
(1,-18)*{}="77",
(-4,-24)*{-2\ }="8",
(1,-24)*{}="88",
(6,-18)*{}="8107",
(12,-18)*{}="810711",
(12,-12)*{}="8107116",
(24,-12)*{}="810711613",
(24,-6)*{}="8107116130",
(30,-6)*{}="81071161300",
(30,0)*{}="810711613000",
(-4,-30)*{-3\ }="9",
(1,-30)*{}="99",
\ar@{->} "arrow1";"arrow1-1"^{}
\ar@{->} "arrow2";"arrow2-1"^{}
\ar@{->} "arrow3";"arrow3-1"^{}
\ar@{-} "1";"-6"^{}
\ar@{-} "1";"2"^{}
\ar@{-} "1";"3"^{}
\ar@{-} "-1";"-1a"^{}
\ar@{-} "-2";"-2a"^{}
\ar@{-} "-3";"-3a"^{}
\ar@{-} "-4";"-4a"^{}
\ar@{-} "-5";"-5a"^{}
\ar@{-} "5";"55"^{}
\ar@{-} "6";"66"^{}
\ar@{-} "8107";"7"^{}
\ar@{-} "8107";"810711"^{}
\ar@{-} "810711";"8107116"^{}
\ar@{-} "8107116";"810711613"^{}
\ar@{-} "810711613";"8107116130"^{}
\ar@{-} "8107116130";"13"^{}
\ar@{-} "7";"6a"^{}
\ar@{-} "6aa";"6a"^{}
\ar@{-} "6aaa";"6aa"^{}
\ar@{-} "6aaaa";"6aaa"^{}
\ar@{-} "st1";"6aaaa"^{}
\ar@{-} "st1";"st2"^{}
\ar@{-} "st2";"st3"^{}
\ar@{-} "st3";"st4"^{}
\ar@{-} "8";"88"^{}
\ar@{-} "9";"99"^{}
\ar@{-} "10";"1010"^{}
\ar@{-} "11";"1111"^{}
\ar@{-} "12";"1212"^{}
\ar@{-} "13";"1313"^{}
\ar@{-} "14";"1414"^{}
\end{xy}\quad
\begin{xy}
(0,0) *{}="1",
(35,0)*{}="2",
(0,-30)*{}="3",
(-5,0)*{(0,2)}="4",
(6,2) *{1}="10",
(6,-1) *{}="1010",
(12,2) *{2}="11",
(12,-1) *{}="1111",
(18,2) *{3}="12",
(18,-1) *{}="1212",
(24,2) *{4}="13",
(24,-1) *{}="1313",
(30,2) *{5}="1515",
(30,-1) *{}="1414",
(-3,-6)*{1\ }="5",
(1,-6)*{}="55",
(0,-12)*{}="l1",
(6,-12)*{}="l2",
(6,-18)*{}="l3",
(24,-18)*{}="l4",
(24,-12)*{}="l5",
(30,-12)*{}="l6",
(30,-6)*{}="l7",
(24,-6)*{}="l8",
(24,0)*{}="l9",
(-4,-12)*{0\ }="6",
(-22,-12)*{\Leftrightarrow}="arrow",
(1,-12)*{}="66",
(-4,-18)*{-1\ }="7",
(1,-18)*{}="77",
(-4,-24)*{-2\ }="8",
(1,-24)*{}="88",
(1,-30)*{}="99",
\ar@{-} "1";"2"^{}
\ar@{-} "1";"3"^{}
\ar@{-} "l1";"l2"^{}
\ar@{-} "l2";"l3"^{}
\ar@{-} "l3";"l4"^{}
\ar@{-} "l4";"l5"^{}
\ar@{-} "l5";"l6"^{}
\ar@{-} "l6";"l7"^{}
\ar@{-} "l7";"l8"^{}
\ar@{-} "l8";"l9"^{}
\end{xy}
\end{equation*}
Here, in the right diagram, we regard
the left plane $\mathbb{R}_{\leq0}\times \mathbb{R}_{\leq 2}$
is filled by boxes. 
As the right diagram drawn in the quarter plane $\mathbb{R}_{\geq0}\times \mathbb{R}_{\leq 2}$,
there might be a space above a box.
It is difficult to describe such diagram as a sequence of integers.
To avoid treating such diagrams, we consider diagrams written in the
half plane $\mathbb{R}\times \mathbb{R}_{\leq k}$.
Although we will describe them in the half plane, 
revised extended Young diagrams are essentially considered as diagrams drawn
in the quarter plane just as in the subsection \ref{EYD-sub}.

\end{rem}

\begin{defn}\label{ad-rem-pt}
We put ${\rm REYD}_{k}={\rm REYD}_{{\rm A}^{(2)},k}$ or ${\rm REYD}_{k}={\rm REYD}_{{\rm D}^{(2)},k}$ in Definition \ref{AEYD}, \ref{DEYD} .
Let $T=(y_t)_{t\in\mathbb{Z}}$ be a sequence in ${\rm REYD}_{k}$ and $i\in\mathbb{Z}$.
\begin{enumerate}
\item Let $T'=(y_t')_{t\in\mathbb{Z}}$ be the sequence such that
$y_i'=y_i-1$ and $y_t'=y_t$ $(t\neq i)$. If $T'\in{\rm REYD}_{k}$ then we say the point $(i,y_i)$
is an admissible point of $T$.
\item Let $T''=(y_t'')_{t\in\mathbb{Z}}$ be the sequence such that
$y_{i-1}''=y_{i-1}+1$ and $y_t''=y_t$ $(t\neq i-1)$. If $T''\in{\rm REYD}_{k}$ then we say the point $(i,y_{i-1})$
is a removable point of $T$.
\end{enumerate}
\end{defn}

For example, if $T$ is the element of ${\rm REYD}_{{\rm A}^{(2)},2}$ in (\ref{reydA-ex1}) then $(i,y_i)=(-2,0)$ is an admissible point
since $T'$ is an element of ${\rm REYD}_{{\rm A}^{(2)},2}$:
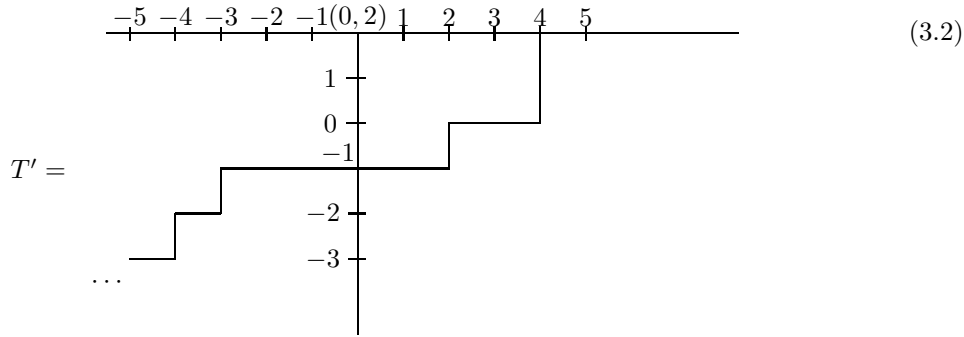
\begin{equation}\label{reydA-ex2}
\begin{xy}
(-42,-18) *{T'=}="YY",
(-33,0) *{}="-6",
(-6,2) *{-1}="-1",
(-12,2) *{-2}="-2",
(-18,2) *{-3}="-3",
(-24,2) *{-4}="-4",
(-30,2) *{-5}="-5",
(-6,-1) *{}="-1a",
(-12,-1) *{}="-2a",
(-18,-1) *{}="-3a",
(-24,-1) *{}="-4a",
(-30,-1) *{}="-5a",
(0,0) *{}="1",
(50,0)*{}="2",
(0,-40)*{}="3",
(0,2)*{(0,2)}="4",
(6,2) *{1}="10",
(6,-1) *{}="1010",
(12,2) *{2}="11",
(12,-1) *{}="1111",
(18,2) *{3}="12",
(18,-1) *{}="1212",
(24,2) *{4}="13",
(24,-1) *{}="1313",
(30,2) *{5}="14",
(30,-1) *{}="1414",
(-3,-6)*{1\ }="5",
(1,-6)*{}="55",
(-3,-12)*{0\ }="6",
(-12,-18)*{}="6aaa",
(-18,-18)*{}="6aaaa",
(-18,-24)*{}="st1",
(-24,-24)*{}="st2",
(-24,-30)*{}="st3",
(-30,-30)*{}="st4",
(-33,-33)*{\cdots}="stdot",
(1,-12)*{}="66",
(-6,-18)*{}="7",
(-2,-16)*{-1\ }="7a",
(1,-18)*{}="77",
(-4,-24)*{-2\ }="8",
(1,-24)*{}="88",
(6,-18)*{}="8107",
(12,-18)*{}="810711",
(12,-12)*{}="8107116",
(24,-12)*{}="810711613",
(24,-6)*{}="8107116130",
(30,-6)*{}="81071161300",
(30,0)*{}="810711613000",
(-4,-30)*{-3\ }="9",
(1,-30)*{}="99",
\ar@{-} "1";"-6"^{}
\ar@{-} "1";"2"^{}
\ar@{-} "1";"3"^{}
\ar@{-} "-1";"-1a"^{}
\ar@{-} "-2";"-2a"^{}
\ar@{-} "-3";"-3a"^{}
\ar@{-} "-4";"-4a"^{}
\ar@{-} "-5";"-5a"^{}
\ar@{-} "5";"55"^{}
\ar@{-} "6";"66"^{}
\ar@{-} "8107";"7"^{}
\ar@{-} "8107";"810711"^{}
\ar@{-} "810711";"8107116"^{}
\ar@{-} "8107116";"810711613"^{}
\ar@{-} "810711613";"8107116130"^{}
\ar@{-} "8107116130";"13"^{}
\ar@{-} "7";"6aaa"^{}
\ar@{-} "6aaaa";"6aaa"^{}
\ar@{-} "st1";"6aaaa"^{}
\ar@{-} "st1";"st2"^{}
\ar@{-} "st2";"st3"^{}
\ar@{-} "st3";"st4"^{}
\ar@{-} "8";"88"^{}
\ar@{-} "9";"99"^{}
\ar@{-} "10";"1010"^{}
\ar@{-} "11";"1111"^{}
\ar@{-} "12";"1212"^{}
\ar@{-} "13";"1313"^{}
\ar@{-} "14";"1414"^{}
\end{xy}
\end{equation}
Similarly, the points $(-1,-1)$, $(2,0)$ and $(4,2)$ are also admissible
and $(4,0)$ and $(2,-1)$ are removable in $T$.

It is easy to see that if $(i,j)$ is an admissible point in an element of ${\rm REYD}_k$ then
\begin{equation}\label{pos-lem1}
i\geq 0 \text{ means } j\leq k, \quad i<0 \text{ means }j\leq k+i.
\end{equation}
If $(i,j)$ is a removable point then
\begin{equation}\label{pos-lem2}
i\geq 1 \text{ means } j< k, \quad i\leq0 \text{ means }j< k+i-1.
\end{equation}

\subsection{Young walls of type ${\rm A}^{(2)}_{2n-2}$ and ${\rm D}^{(2)}_{n}$}

Following \cite{Kang}, let us recall the notion of Young walls of type ${\rm A}^{(2)}_{2n-2}$ and ${\rm D}^{(2)}_{n}$.
In the original paper, the Young wall consists of $I$-colored blocks of three different
shapes:
\begin{enumerate}
\item[(1)] block with unit width, unit height and unit thickness:
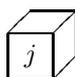
\[
\begin{xy}
(3,3) *{j}="0",
(0,0) *{}="1",
(6,0)*{}="2",
(6,6)*{}="3",
(0,6)*{}="4",
(3,9)*{}="5",
(9,9)*{}="6",
(9,3)*{}="7",
\ar@{-} "1";"2"^{}
\ar@{-} "1";"4"^{}
\ar@{-} "2";"3"^{}
\ar@{-} "3";"4"^{}
\ar@{-} "5";"4"^{}
\ar@{-} "5";"6"^{}
\ar@{-} "3";"6"^{}
\ar@{-} "2";"7"^{}
\ar@{-} "6";"7"^{}
\end{xy}
\]
\item[(2)] block with unit width, unit height and half-unit thickness:
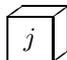
\[
\begin{xy}
(3,3) *{j}="0",
(0,0) *{}="1",
(6,0)*{}="2",
(6,6)*{}="3",
(0,6)*{}="4",
(2,7.5)*{}="5",
(8,7.5)*{}="6",
(8,1.5)*{}="7",
\ar@{-} "1";"2"^{}
\ar@{-} "1";"4"^{}
\ar@{-} "2";"3"^{}
\ar@{-} "3";"4"^{}
\ar@{-} "5";"4"^{}
\ar@{-} "5";"6"^{}
\ar@{-} "3";"6"^{}
\ar@{-} "2";"7"^{}
\ar@{-} "6";"7"^{}
\end{xy}
\]
\item[(3)] block with unit width, half-unit height and unit thickness:
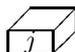
\[
\begin{xy}
(3,1.5) *{j}="0",
(0,0) *{}="1",
(6,0)*{}="2",
(6,3)*{}="3",
(0,3)*{}="4",
(3,6)*{}="5",
(9,6)*{}="6",
(9,3)*{}="7",
\ar@{-} "1";"2"^{}
\ar@{-} "1";"4"^{}
\ar@{-} "2";"3"^{}
\ar@{-} "3";"4"^{}
\ar@{-} "5";"4"^{}
\ar@{-} "5";"6"^{}
\ar@{-} "3";"6"^{}
\ar@{-} "2";"7"^{}
\ar@{-} "6";"7"^{}
\end{xy}
\]
\end{enumerate}
In this article, blocks of second shape will not be used. Just as in \cite{Kang},
we simply describe blocks (1) with color $j\in I$
as
\begin{equation}\label{smpl1}
\begin{xy}
(3,3) *{j}="0",
(0,0) *{}="1",
(6,0)*{}="2",
(6,6)*{}="3",
(0,6)*{}="4",
\ar@{-} "1";"2"^{}
\ar@{-} "1";"4"^{}
\ar@{-} "2";"3"^{}
\ar@{-} "3";"4"^{}
\end{xy}
\end{equation}
and (3) with color $j\in I$ as
\begin{equation}\label{smpl2}
\begin{xy}
(1.5,1.5) *{\ \ j}="0",
(0,0) *{}="1",
(6,0)*{}="2",
(6,3.5)*{}="3",
(0,3.5)*{}="4",
\ar@{-} "1";"2"^{}
\ar@{-} "1";"4"^{}
\ar@{-} "2";"3"^{}
\ar@{-} "3";"4"^{}
\end{xy}
\end{equation}
The blocks (\ref{smpl1}) and (\ref{smpl2}) are called $j$-blocks.
If colored blocks are stacked as follows
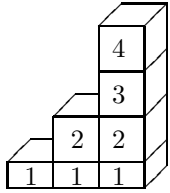
\[
\begin{xy}
(-9.5,1.5) *{\ 1}="000",
(-3.5,1.5) *{\ 1}="00",
(1.5,1.5) *{\ \ 1}="0",
(1.5,6.5) *{\ \ 2}="02",
(-4,6.5) *{\ \ 2}="002",
(1.5,12.5) *{\ \ 3}="03",
(1.5,18.5) *{\ \ 4}="04",
(9,3)*{}="2-a",
(9,6.5)*{}="3-a",
(9,12.5)*{}="3-1-a",
(9,18.5)*{}="3-2-a",
(9,24.5)*{}="3-3-a",
(3,24.5)*{}="4-3-a",
(-3,12.5)*{}="5-1-a",
(0,12.5)*{}="5-1-ab",
(-9,6.5)*{}="7-a",
(-6,6.5)*{}="7-ab",
(0,0) *{}="1",
(6,0)*{}="2",
(6,3.5)*{}="3",
(0,3.5)*{}="4",
(-6,3.5)*{}="5",
(-6,0)*{}="6",
(-12,3.5)*{}="7",
(-12,0)*{}="8",
(6,9.5)*{}="3-1",
(6,15.5)*{}="3-2",
(6,21.5)*{}="3-3",
(0,9.5)*{}="4-1",
(0,15.5)*{}="4-2",
(0,21.5)*{}="4-3",
(-6,9.5)*{}="5-1",
\ar@{-} "1";"2"^{}
\ar@{-} "1";"4"^{}
\ar@{-} "2";"3"^{}
\ar@{-} "3";"4"^{}
\ar@{-} "5";"6"^{}
\ar@{-} "5";"4"^{}
\ar@{-} "1";"6"^{}
\ar@{-} "7";"8"^{}
\ar@{-} "7";"5"^{}
\ar@{-} "6";"8"^{}
\ar@{-} "3";"3-1"^{}
\ar@{-} "4-1";"3-1"^{}
\ar@{-} "4-1";"4"^{}
\ar@{-} "5-1";"5"^{}
\ar@{-} "4-1";"5-1"^{}
\ar@{-} "4-1";"4-2"^{}
\ar@{-} "3-2";"3-1"^{}
\ar@{-} "3-2";"4-2"^{}
\ar@{-} "3-2";"3-3"^{}
\ar@{-} "4-3";"4-2"^{}
\ar@{-} "4-3";"3-3"^{}
\ar@{-} "2";"2-a"^{}
\ar@{-} "3";"3-a"^{}
\ar@{-} "3-1";"3-1-a"^{}
\ar@{-} "3-2";"3-2-a"^{}
\ar@{-} "3-3";"3-3-a"^{}
\ar@{-} "2-a";"3-3-a"^{}
\ar@{-} "4-3";"4-3-a"^{}
\ar@{-} "3-3-a";"4-3-a"^{}
\ar@{-} "5-1";"5-1-a"^{}
\ar@{-} "5-1-ab";"5-1-a"^{}
\ar@{-} "7";"7-a"^{}
\ar@{-} "7-ab";"7-a"^{}
\end{xy}
\]
it is simply described as
\begin{equation}\label{smpl3}
\begin{xy}
(-9.5,1.5) *{\ 1}="000",
(-3.5,1.5) *{\ 1}="00",
(1.5,1.5) *{\ \ 1}="0",
(1.5,6.5) *{\ \ 2}="02",
(-4,6.5) *{\ \ 2}="002",
(1.5,12.5) *{\ \ 3}="03",
(1.5,18.5) *{\ \ 4}="04",
(0,0) *{}="1",
(6,0)*{}="2",
(6,3.5)*{}="3",
(0,3.5)*{}="4",
(-6,3.5)*{}="5",
(-6,0)*{}="6",
(-12,3.5)*{}="7",
(-12,0)*{}="8",
(6,9.5)*{}="3-1",
(6,15.5)*{}="3-2",
(6,21.5)*{}="3-3",
(0,9.5)*{}="4-1",
(0,15.5)*{}="4-2",
(0,21.5)*{}="4-3",
(-6,9.5)*{}="5-1",
\ar@{-} "1";"2"^{}
\ar@{-} "1";"4"^{}
\ar@{-} "2";"3"^{}
\ar@{-} "3";"4"^{}
\ar@{-} "5";"6"^{}
\ar@{-} "5";"4"^{}
\ar@{-} "1";"6"^{}
\ar@{-} "7";"8"^{}
\ar@{-} "7";"5"^{}
\ar@{-} "6";"8"^{}
\ar@{-} "3";"3-1"^{}
\ar@{-} "4-1";"3-1"^{}
\ar@{-} "4-1";"4"^{}
\ar@{-} "5-1";"5"^{}
\ar@{-} "4-1";"5-1"^{}
\ar@{-} "4-1";"4-2"^{}
\ar@{-} "3-2";"3-1"^{}
\ar@{-} "3-2";"4-2"^{}
\ar@{-} "3-2";"3-3"^{}
\ar@{-} "4-3";"4-2"^{}
\ar@{-} "4-3";"3-3"^{}
\end{xy}
\end{equation}
We will call these diagrams `walls'.
One use the following colored blocks for each case: 
\begin{itemize}
\item In type ${\rm A}^{(2)}_{2n-2}$-case, one use
\[
\begin{xy}
(1.5,1.5) *{\ \ 1}="0",
(0,0) *{}="1",
(6,0)*{}="2",
(6,3.5)*{}="3",
(0,3.5)*{}="4",
\ar@{-} "1";"2"^{}
\ar@{-} "1";"4"^{}
\ar@{-} "2";"3"^{}
\ar@{-} "3";"4"^{}
\end{xy}
\]
and
\[
\begin{xy}
(3,3) *{j}="0",
(0,0) *{}="1",
(6,0)*{}="2",
(6,6)*{}="3",
(0,6)*{}="4",
\ar@{-} "1";"2"^{}
\ar@{-} "1";"4"^{}
\ar@{-} "2";"3"^{}
\ar@{-} "3";"4"^{}
\end{xy}
\]
for $j=2,3,\cdots,n$.
\item
In type ${\rm D}^{(2)}_{n}$-case, one use
\[
\begin{xy}
(1.5,1.5) *{\ \ 1}="0",
(0,0) *{}="1",
(6,0)*{}="2",
(6,3.5)*{}="3",
(0,3.5)*{}="4",
(16,1.5) *{\ \ \ \ \ n}="00",
(16,0) *{}="11",
(22,0)*{}="22",
(22,3.5)*{}="33",
(16,3.5)*{}="44",
\ar@{-} "1";"2"^{}
\ar@{-} "1";"4"^{}
\ar@{-} "2";"3"^{}
\ar@{-} "3";"4"^{}
\ar@{-} "11";"22"^{}
\ar@{-} "11";"44"^{}
\ar@{-} "22";"33"^{}
\ar@{-} "33";"44"^{}
\end{xy}
\]
and
\[
\begin{xy}
(3,3) *{j}="0",
(0,0) *{}="1",
(6,0)*{}="2",
(6,6)*{}="3",
(0,6)*{}="4",
\ar@{-} "1";"2"^{}
\ar@{-} "1";"4"^{}
\ar@{-} "2";"3"^{}
\ar@{-} "3";"4"^{}
\end{xy}
\]
for $j=2,3,\cdots,n-1$.
\end{itemize}
To define the Young walls, we consider `ground state walls'.
\begin{defn}\cite{Kang}\label{def-gsw}
In type ${\rm A}^{(2)}_{2n-2}$-case, the ground state wall $Y_{\Lambda_1}$ is defined as
\[
Y_{\Lambda_1}=
\begin{xy}
(-15.5,1.5) *{\ \cdots}="0000",
(-9.5,1.5) *{\ 1}="000",
(-3.5,1.5) *{\ 1}="00",
(1.5,1.5) *{\ \ 1}="0",
(0,0) *{}="1",
(6,0)*{}="2",
(6,3.5)*{}="3",
(0,3.5)*{}="4",
(-6,3.5)*{}="5",
(-6,0)*{}="6",
(-12,3.5)*{}="7",
(-12,0)*{}="8",
\ar@{-} "1";"2"^{}
\ar@{-} "1";"4"^{}
\ar@{-} "2";"3"^{}
\ar@{-} "3";"4"^{}
\ar@{-} "5";"6"^{}
\ar@{-} "5";"4"^{}
\ar@{-} "1";"6"^{}
\ar@{-} "7";"8"^{}
\ar@{-} "7";"5"^{}
\ar@{-} "6";"8"^{}
\end{xy}
\]
Here, in $Y_{\Lambda_1}$, the block which has half-unit height with color $1$ extends infinitely to the left.

In type ${\rm D}^{(2)}_{n}$-case, the ground state walls $Y_{\Lambda_1}$ and  $Y_{\Lambda_n}$ are defined as
\[
Y_{\Lambda_1}=
\begin{xy}
(-15.5,1.5) *{\ \cdots}="0000",
(-9.5,1.5) *{\ 1}="000",
(-3.5,1.5) *{\ 1}="00",
(1.5,1.5) *{\ \ 1}="0",
(0,0) *{}="1",
(6,0)*{}="2",
(6,3.5)*{}="3",
(0,3.5)*{}="4",
(-6,3.5)*{}="5",
(-6,0)*{}="6",
(-12,3.5)*{}="7",
(-12,0)*{}="8",
\ar@{-} "1";"2"^{}
\ar@{-} "1";"4"^{}
\ar@{-} "2";"3"^{}
\ar@{-} "3";"4"^{}
\ar@{-} "5";"6"^{}
\ar@{-} "5";"4"^{}
\ar@{-} "1";"6"^{}
\ar@{-} "7";"8"^{}
\ar@{-} "7";"5"^{}
\ar@{-} "6";"8"^{}
\end{xy}
\]
and
\[
Y_{\Lambda_n}=
\begin{xy}
(-15.5,1.5) *{\ \cdots}="0000",
(-9.5,1.5) *{\ n}="000",
(-3.5,1.5) *{\ n}="00",
(1.5,1.5) *{\ \ n}="0",
(0,0) *{}="1",
(6,0)*{}="2",
(6,3.5)*{}="3",
(0,3.5)*{}="4",
(-6,3.5)*{}="5",
(-6,0)*{}="6",
(-12,3.5)*{}="7",
(-12,0)*{}="8",
\ar@{-} "1";"2"^{}
\ar@{-} "1";"4"^{}
\ar@{-} "2";"3"^{}
\ar@{-} "3";"4"^{}
\ar@{-} "5";"6"^{}
\ar@{-} "5";"4"^{}
\ar@{-} "1";"6"^{}
\ar@{-} "7";"8"^{}
\ar@{-} "7";"5"^{}
\ar@{-} "6";"8"^{}
\end{xy}
\]
\end{defn}

In \cite{Kang}, the generalized Cartan matrix of type ${\rm A}_{2n-2}^{(2)}$
is the transposed matrix of ours so that $V(\Lambda_1)$ is the level $1$ representation.
Thus, strictly speaking, the above $\Lambda_1$ of type ${\rm A}_{2n-2}^{(2)}$ is the fundamental weight of $^L\mathfrak{g}$.

\begin{defn}\cite{Kang}\label{def-YW}
A wall is called a {\it Young wall} of ground state $\lambda=\Lambda_1$ (resp. $\lambda\in\{\Lambda_1,\Lambda_n\}$)
of type ${\rm A}^{(2)}_{2n-2}$ (resp. ${\rm D}^{(2)}_{n}$) if it satisfies the following:
\begin{enumerate}
\item The wall is built on top of the ground state wall $Y_{\lambda}$. Finitely many colored blocks
are stacked on $Y_{\lambda}$.
\item The colored blocks are stacked in the patterns we give below for each type and $\lambda$.
\item Let $h_j$ be the height of $j$-th column of the wall from the right. Then it holds $h_j\geq h_{j+1}$. 
\end{enumerate}
Here, the patterns mentioned in (ii) are as follows: 

In type ${\rm A}^{(2)}_{2n-2}$-case, $\lambda=\Lambda_1$:
\[
\begin{xy}
(-15.5,-2) *{\ 1}="000-1",
(-9.5,-2) *{\ 1}="00-1",
(-3.5,-2) *{\ 1}="0-1",
(1.5,-2) *{\ \ 1}="0-1",
(-15.5,1.5) *{\ 1}="0000",
(-9.5,1.5) *{\ 1}="000",
(-3.5,1.5) *{\ 1}="00",
(1.5,1.5) *{\ \ 1}="0",
(1.5,6.5) *{\ \ 2}="02",
(-4,6.5) *{\ \ 2}="002",
(-10,6.5) *{\ \ 2}="0002",
(-16,6.5) *{\ \ 2}="00002",
(1.5,12.5) *{\ \ 3}="03",
(-4,12.5) *{\ \ 3}="003",
(-10,12.5) *{\ \ 3}="0003",
(-16,12.5) *{\ \ 3}="00003",
(1.5,20.5) *{\ \ \vdots}="04",
(-4,20.5) *{\ \ \vdots}="004",
(-10,20.5) *{\ \ \vdots}="0004",
(-16,20.5) *{\ \ \vdots}="00004",
(1.5,27.5) *{\ \ _{n-1}}="05",
(-4,27.5) *{\ \ _{n-1}}="005",
(-10,27.5) *{\ \ _{n-1}}="0005",
(-16,27.5) *{\ \ _{n-1}}="00005",
(1.5,33.5) *{\ \ n}="06",
(-4,33.5) *{\ \ n}="006",
(-25,34) *{\cdots}="dots",
(-10,33.5) *{\ \ n}="0006",
(-16,33.5) *{\ \ n}="00006",
(1.5,39.5) *{\ \ _{n-1}}="07",
(-4,39.5) *{\ \ _{n-1}}="007",
(-10,39.5) *{\ \ _{n-1}}="0007",
(-16,39.5) *{\ \ _{n-1}}="00007",
(1.5,48) *{\ \ \vdots}="08",
(-4,48) *{\ \ \vdots}="008",
(-10,48) *{\ \ \vdots}="0008",
(-16,48) *{\ \ \vdots}="00008",
(1.5,55) *{\ \ 2}="09",
(-4,55) *{\ \ 2}="009",
(-10,55) *{\ \ 2}="0009",
(-16,55) *{\ \ 2}="00009",
(1.5,60) *{\ \ 1}="010",
(-4,60) *{\ \ 1}="0010",
(-10,60) *{\ \ 1}="00010",
(-16,60) *{\ \ 1}="000010",
(1.5,64) *{\ \ 1}="011",
(-4,64) *{\ \ 1}="0011",
(-10,64) *{\ \ 1}="00011",
(-16,64) *{\ \ 1}="000011",
(1.5,69) *{\ \ 2}="012",
(-4,69) *{\ \ 2}="0012",
(-10,69) *{\ \ 2}="00012",
(-16,69) *{\ \ 2}="000012",
(-18,-3.5) *{}="0.5-l",
(-20,-3.5) *{}="0-l",
(-20,0) *{}="1-l",
(-20,3.5) *{}="2-l",
(-20,9.5) *{}="3-l",
(-20,15.5) *{}="4-l",
(-20,24.5) *{}="5-l",
(-20,30.5) *{}="6-l",
(-20,36.5) *{}="7-l",
(-20,42.5) *{}="8-l",
(-20,52) *{}="9-l",
(-20,58) *{}="10-l",
(-20,62) *{}="11-l",
(-20,66) *{}="12-l",
(-20,72) *{}="13-l",
(-18,74) *{}="13.5-l",
(6,-3.5)*{}="r--1",
(6,0)*{}="r",
(6,3.5)*{}="r-0",
(0,-3.5) *{}="b1",
(-6,-3.5)*{}="b2",
(-12,-3.5)*{}="b3",
(0,74) *{}="t1",
(-6,74)*{}="t2",
(-12,74)*{}="t3",
(6,9.5)*{}="r-1",
(6,15.5)*{}="r-2",
(6,24.5)*{}="r-3",
(6,30.5)*{}="r-4",
(6,36.5)*{}="r-5",
(6,42.5)*{}="r-6",
(6,52)*{}="r-7",
(6,58)*{}="r-8",
(6,62)*{}="r-9",
(6,66)*{}="r-10",
(6,72)*{}="r-11",
(6,74)*{}="r-11.5",
\ar@{-} "b1";"t1"^{}
\ar@{-} "b2";"t2"^{}
\ar@{-} "b3";"t3"^{}
\ar@{-} "0.5-l";"13.5-l"^{}
\ar@{-} "r-11.5";"r--1"^{}
\ar@{-} "0-l";"r--1"^{}
\ar@{-} "1-l";"r"^{}
\ar@{-} "2-l";"r-0"^{}
\ar@{-} "3-l";"r-1"^{}
\ar@{-} "4-l";"r-2"^{}
\ar@{-} "5-l";"r-3"^{}
\ar@{-} "6-l";"r-4"^{}
\ar@{-} "7-l";"r-5"^{}
\ar@{-} "8-l";"r-6"^{}
\ar@{-} "9-l";"r-7"^{}
\ar@{-} "10-l";"r-8"^{}
\ar@{-} "11-l";"r-9"^{}
\ar@{-} "12-l";"r-10"^{}
\ar@{-} "13-l";"r-11"^{}
\end{xy}
\]
In type ${\rm D}^{(2)}_{n}$-case,
\[
\begin{xy}
(-40,35) *{\lambda=\Lambda_1:}="weight",
(-15.5,-2) *{\ 1}="000-1",
(-9.5,-2) *{\ 1}="00-1",
(-3.5,-2) *{\ 1}="0-1",
(1.5,-2) *{\ \ 1}="0-1",
(-15.5,1.5) *{\ 1}="0000",
(-9.5,1.5) *{\ 1}="000",
(-3.5,1.5) *{\ 1}="00",
(1.5,1.5) *{\ \ 1}="0",
(1.5,6.5) *{\ \ 2}="02",
(-4,6.5) *{\ \ 2}="002",
(-10,6.5) *{\ \ 2}="0002",
(-16,6.5) *{\ \ 2}="00002",
(1.5,12.5) *{\ \ 3}="03",
(-4,12.5) *{\ \ 3}="003",
(-10,12.5) *{\ \ 3}="0003",
(-16,12.5) *{\ \ 3}="00003",
(1.5,20.5) *{\ \ \vdots}="04",
(-4,20.5) *{\ \ \vdots}="004",
(-10,20.5) *{\ \ \vdots}="0004",
(-16,20.5) *{\ \ \vdots}="00004",
(1.5,27.5) *{\ \ _{n-1}}="05",
(-4,27.5) *{\ \ _{n-1}}="005",
(-10,27.5) *{\ \ _{n-1}}="0005",
(-16,27.5) *{\ \ _{n-1}}="00005",
(1.5,32) *{\ \ _n}="06",
(-4,32) *{\ \ _n}="006",
(-25,34) *{\cdots}="dots",
(-10,32) *{\ \ _n}="0006",
(-16,32) *{\ \ _n}="00006",
(1.5,35) *{\ \ _n}="06.5",
(-4,35) *{\ \ _n}="006.5",
(-10,35) *{\ \ _n}="0006.5",
(-16,35) *{\ \ _n}="00006.5",
(1.5,39.5) *{\ \ _{n-1}}="07",
(-4,39.5) *{\ \ _{n-1}}="007",
(-10,39.5) *{\ \ _{n-1}}="0007",
(-16,39.5) *{\ \ _{n-1}}="00007",
(1.5,48) *{\ \ \vdots}="08",
(-4,48) *{\ \ \vdots}="008",
(-10,48) *{\ \ \vdots}="0008",
(-16,48) *{\ \ \vdots}="00008",
(1.5,55) *{\ \ 2}="09",
(-4,55) *{\ \ 2}="009",
(-10,55) *{\ \ 2}="0009",
(-16,55) *{\ \ 2}="00009",
(1.5,60) *{\ \ 1}="010",
(-4,60) *{\ \ 1}="0010",
(-10,60) *{\ \ 1}="00010",
(-16,60) *{\ \ 1}="000010",
(1.5,64) *{\ \ 1}="011",
(-4,64) *{\ \ 1}="0011",
(-10,64) *{\ \ 1}="00011",
(-16,64) *{\ \ 1}="000011",
(1.5,69) *{\ \ 2}="012",
(-4,69) *{\ \ 2}="0012",
(-10,69) *{\ \ 2}="00012",
(-16,69) *{\ \ 2}="000012",
(-18,-3.5) *{}="0.5-l",
(-20,-3.5) *{}="0-l",
(-20,0) *{}="1-l",
(-20,3.5) *{}="2-l",
(-20,9.5) *{}="3-l",
(-20,15.5) *{}="4-l",
(-20,24.5) *{}="5-l",
(-20,30.5) *{}="6-l",
(-20,33.5) *{}="65-l",
(-20,36.5) *{}="7-l",
(-20,42.5) *{}="8-l",
(-20,52) *{}="9-l",
(-20,58) *{}="10-l",
(-20,62) *{}="11-l",
(-20,66) *{}="12-l",
(-20,72) *{}="13-l",
(-18,74) *{}="13.5-l",
(6,-3.5)*{}="r--1",
(6,0)*{}="r",
(6,3.5)*{}="r-0",
(0,-3.5) *{}="b1",
(-6,-3.5)*{}="b2",
(-12,-3.5)*{}="b3",
(0,74) *{}="t1",
(-6,74)*{}="t2",
(-12,74)*{}="t3",
(6,9.5)*{}="r-1",
(6,15.5)*{}="r-2",
(6,24.5)*{}="r-3",
(6,30.5)*{}="r-4",
(6,33.5)*{}="r-45",
(6,36.5)*{}="r-5",
(6,42.5)*{}="r-6",
(6,52)*{}="r-7",
(6,58)*{}="r-8",
(6,62)*{}="r-9",
(6,66)*{}="r-10",
(6,72)*{}="r-11",
(6,74)*{}="r-11.5",
\ar@{-} "b1";"t1"^{}
\ar@{-} "b2";"t2"^{}
\ar@{-} "b3";"t3"^{}
\ar@{-} "0.5-l";"13.5-l"^{}
\ar@{-} "r-11.5";"r--1"^{}
\ar@{-} "0-l";"r--1"^{}
\ar@{-} "1-l";"r"^{}
\ar@{-} "2-l";"r-0"^{}
\ar@{-} "3-l";"r-1"^{}
\ar@{-} "4-l";"r-2"^{}
\ar@{-} "5-l";"r-3"^{}
\ar@{-} "6-l";"r-4"^{}
\ar@{-} "65-l";"r-45"^{}
\ar@{-} "7-l";"r-5"^{}
\ar@{-} "8-l";"r-6"^{}
\ar@{-} "9-l";"r-7"^{}
\ar@{-} "10-l";"r-8"^{}
\ar@{-} "11-l";"r-9"^{}
\ar@{-} "12-l";"r-10"^{}
\ar@{-} "13-l";"r-11"^{}
\end{xy}\qquad \qquad
\begin{xy}
(-40,35) *{\lambda=\Lambda_n:}="weight",
(-15.5,-2) *{\ _n}="000-1",
(-9.5,-2) *{\ _n}="00-1",
(-3.5,-2) *{\ _n}="0-1",
(1.5,-2) *{\ \ _n}="0-1",
(-15.5,1.5) *{\ _n}="0000",
(-9.5,1.5) *{\ _n}="000",
(-3.5,1.5) *{\ _n}="00",
(1.5,1.5) *{\ \ _n}="0",
(1.5,6.5) *{\ \ _{n-1}}="02",
(-4,6.5) *{\ \ _{n-1}}="002",
(-10,6.5) *{\ \ _{n-1}}="0002",
(-16,6.5) *{\ \ _{n-1}}="00002",
(1.5,12.5) *{\ \ _{n-2}}="03",
(-4,12.5) *{\ \ _{n-2}}="003",
(-10,12.5) *{\ \ _{n-2}}="0003",
(-16,12.5) *{\ \ _{n-2}}="00003",
(1.5,20.5) *{\ \ \vdots}="04",
(-4,20.5) *{\ \ \vdots}="004",
(-10,20.5) *{\ \ \vdots}="0004",
(-16,20.5) *{\ \ \vdots}="00004",
(1.5,27.5) *{\ \ 2}="05",
(-4,27.5) *{\ \ 2}="005",
(-10,27.5) *{\ \ 2}="0005",
(-16,27.5) *{\ \ 2}="00005",
(1.5,32) *{\ \ 1}="06",
(-4,32) *{\ \ 1}="006",
(-25,34) *{\cdots}="dots",
(-10,32) *{\ \ 1}="0006",
(-16,32) *{\ \ 1}="00006",
(1.5,35) *{\ \ 1}="06.5",
(-4,35) *{\ \ 1}="006.5",
(-10,35) *{\ \ 1}="0006.5",
(-16,35) *{\ \ 1}="00006.5",
(1.5,39.5) *{\ \ 2}="07",
(-4,39.5) *{\ \ 2}="007",
(-10,39.5) *{\ \ 2}="0007",
(-16,39.5) *{\ \ 2}="00007",
(1.5,48) *{\ \ \vdots}="08",
(-4,48) *{\ \ \vdots}="008",
(-10,48) *{\ \ \vdots}="0008",
(-16,48) *{\ \ \vdots}="00008",
(1.5,55) *{\ \ _{n-1}}="09",
(-4,55) *{\ \ _{n-1}}="009",
(-10,55) *{\ \ _{n-1}}="0009",
(-16,55) *{\ \ _{n-1}}="00009",
(1.5,60) *{\ \ _{n}}="010",
(-4,60) *{\ \ _{n}}="0010",
(-10,60) *{\ \ _{n}}="00010",
(-16,60) *{\ \ _{n}}="000010",
(1.5,64) *{\ \ _{n}}="011",
(-4,64) *{\ \ _{n}}="0011",
(-10,64) *{\ \ _{n}}="00011",
(-16,64) *{\ \ _{n}}="000011",
(1.5,69) *{\ \ _{n-1}}="012",
(-4,69) *{\ \ _{n-1}}="0012",
(-10,69) *{\ \ _{n-1}}="00012",
(-16,69) *{\ \ _{n-1}}="000012",
(-18,-3.5) *{}="0.5-l",
(-20,-3.5) *{}="0-l",
(-20,0) *{}="1-l",
(-20,3.5) *{}="2-l",
(-20,9.5) *{}="3-l",
(-20,15.5) *{}="4-l",
(-20,24.5) *{}="5-l",
(-20,30.5) *{}="6-l",
(-20,33.5) *{}="65-l",
(-20,36.5) *{}="7-l",
(-20,42.5) *{}="8-l",
(-20,52) *{}="9-l",
(-20,58) *{}="10-l",
(-20,62) *{}="11-l",
(-20,66) *{}="12-l",
(-20,72) *{}="13-l",
(-18,74) *{}="13.5-l",
(6,-3.5)*{}="r--1",
(6,0)*{}="r",
(6,3.5)*{}="r-0",
(0,-3.5) *{}="b1",
(-6,-3.5)*{}="b2",
(-12,-3.5)*{}="b3",
(0,74) *{}="t1",
(-6,74)*{}="t2",
(-12,74)*{}="t3",
(6,9.5)*{}="r-1",
(6,15.5)*{}="r-2",
(6,24.5)*{}="r-3",
(6,30.5)*{}="r-4",
(6,33.5)*{}="r-45",
(6,36.5)*{}="r-5",
(6,42.5)*{}="r-6",
(6,52)*{}="r-7",
(6,58)*{}="r-8",
(6,62)*{}="r-9",
(6,66)*{}="r-10",
(6,72)*{}="r-11",
(6,74)*{}="r-11.5",
\ar@{-} "b1";"t1"^{}
\ar@{-} "b2";"t2"^{}
\ar@{-} "b3";"t3"^{}
\ar@{-} "0.5-l";"13.5-l"^{}
\ar@{-} "r-11.5";"r--1"^{}
\ar@{-} "0-l";"r--1"^{}
\ar@{-} "1-l";"r"^{}
\ar@{-} "2-l";"r-0"^{}
\ar@{-} "3-l";"r-1"^{}
\ar@{-} "4-l";"r-2"^{}
\ar@{-} "5-l";"r-3"^{}
\ar@{-} "6-l";"r-4"^{}
\ar@{-} "65-l";"r-45"^{}
\ar@{-} "7-l";"r-5"^{}
\ar@{-} "8-l";"r-6"^{}
\ar@{-} "9-l";"r-7"^{}
\ar@{-} "10-l";"r-8"^{}
\ar@{-} "11-l";"r-9"^{}
\ar@{-} "12-l";"r-10"^{}
\ar@{-} "13-l";"r-11"^{}
\end{xy}
\]
Here, the first row of each pattern from the bottom is
the ground state wall.
\end{defn}

\begin{ex}\label{ex-A2-0}
The following is a Young wall of type ${\rm A}^{(2)}_{4}$ of ground state $\Lambda_1$. 
\begin{equation*}
\begin{xy}
(-15.5,-2) *{\ 1}="000-1",
(-9.5,-2) *{\ 1}="00-1",
(-3.5,-2) *{\ 1}="0-1",
(1.5,-2) *{\ \ 1}="0-1",
(-21.5,-2) *{\dots}="00000",
(-9.5,1.5) *{\ 1}="000",
(-3.5,1.5) *{\ 1}="00",
(1.5,1.5) *{\ \ 1}="0",
(1.5,6.5) *{\ \ 2}="02",
(-4,6.5) *{\ \ 2}="002",
(1.5,12.5) *{\ \ 3}="03",
(1.5,18.5) *{\ \ 2}="04",
(0,0) *{}="1",
(0,-3.5) *{}="1-u",
(6,0)*{}="2",
(6,-3.5)*{}="2-u",
(6,3.5)*{}="3",
(0,3.5)*{}="4",
(-6,3.5)*{}="5",
(-6,0)*{}="6",
(-6,-3.5)*{}="6-u",
(-12,3.5)*{}="7",
(-12,0)*{}="8",
(-12,-3.6)*{}="8-u",
(-18,3.5)*{}="9",
(-18,0)*{}="10-a",
(-18,0)*{}="10",
(-18,-3.5)*{}="10-u",
(6,9.5)*{}="3-1",
(6,15.5)*{}="3-2",
(6,21.5)*{}="3-3",
(0,9.5)*{}="4-1",
(0,15.5)*{}="4-2",
(0,21.5)*{}="4-3",
(-6,9.5)*{}="5-1",
\ar@{-} "8";"10-a"^{}
\ar@{-} "10-u";"10-a"^{}
\ar@{-} "10-u";"2-u"^{}
\ar@{-} "8";"8-u"^{}
\ar@{-} "6";"6-u"^{}
\ar@{-} "2";"2-u"^{}
\ar@{-} "1";"1-u"^{}
\ar@{-} "1";"2"^{}
\ar@{-} "1";"4"^{}
\ar@{-} "2";"3"^{}
\ar@{-} "3";"4"^{}
\ar@{-} "5";"6"^{}
\ar@{-} "5";"4"^{}
\ar@{-} "1";"6"^{}
\ar@{-} "7";"8"^{}
\ar@{-} "7";"5"^{}
\ar@{-} "6";"8"^{}
\ar@{-} "3";"3-1"^{}
\ar@{-} "4-1";"3-1"^{}
\ar@{-} "4-1";"4"^{}
\ar@{-} "5-1";"5"^{}
\ar@{-} "4-1";"5-1"^{}
\ar@{-} "4-1";"4-2"^{}
\ar@{-} "3-2";"3-1"^{}
\ar@{-} "3-2";"4-2"^{}
\ar@{-} "3-2";"3-3"^{}
\ar@{-} "4-3";"4-2"^{}
\ar@{-} "4-3";"3-3"^{}
\end{xy}
\end{equation*}
\end{ex}

\begin{defn}\cite{Kang}\label{def-YW2}
Let $Y$ be a Young wall of type ${\rm A}^{(2)}_{2n-2}$ or ${\rm D}^{(2)}_{n}$ of ground state $\lambda$.
\begin{enumerate}
\item A column of $Y$ is called a {\it full column} if its height is a multiple of the unit length.
\item $Y$ is said to be {\it proper} if none of two full columns of $Y$ have the same height.
\end{enumerate}
\end{defn}

\begin{defn}\cite{Kang}\label{def-YW2a}
Let $Y$ be a proper Young wall of type ${\rm A}^{(2)}_{2n-2}$ or ${\rm D}^{(2)}_{n}$.
\begin{enumerate}
\item A block colored by $i\in I$ in $Y$ is said to be {\it removable} $i$-{\it block} 
if the wall obtained from $Y$ by removing this block 
remains a proper Young wall. 
\item If we obtain a proper Young wall by adding
an block colored by $i\in I$ to a place of $Y$ then the place is said to be an
$i$-{\it admissible} {\it slot}.
\end{enumerate}
\end{defn}
\begin{defn}
Let $Y$ be a proper Young wall of type ${\rm A}^{(2)}_{2n-2}$ or ${\rm D}^{(2)}_{n}$
and $t=1$ or $t=n$. 
\begin{enumerate}
\item[(i)]
Let $Y'$ be a wall
obtained by adding two $t$-blocks of shape (\ref{smpl2}) to a column of $Y$:
\[
Y=
\begin{xy}
(14,10) *{\leftarrow A}="A",
(-3,3) *{\cdots}="dot1",
(12,3) *{\cdots}="dot2",
(0,-4) *{}="-1",
(8,-4)*{}="-2",
(4,-6)*{\vdots}="dot3",
(0,0) *{}="1",
(8,0)*{}="2",
(8,8)*{}="3",
(0,8)*{}="4",
(8,12)*{}="5",
(0,12)*{}="6",
(8,16)*{}="7",
(0,16)*{}="8",
\ar@{-} "1";"-1"^{}
\ar@{-} "-2";"2"^{}
\ar@{-} "1";"2"^{}
\ar@{-} "1";"4"^{}
\ar@{-} "2";"3"^{}
\ar@{-} "3";"4"^{}
\ar@{--} "3";"5"^{}
\ar@{--} "4";"6"^{}
\ar@{--} "6";"5"^{}
\ar@{--} "7";"5"^{}
\ar@{--} "6";"8"^{}
\ar@{--} "7";"8"^{}
\end{xy}\qquad
Y'=
\begin{xy}
(14,14) *{}="B",
(4,10)*{t}="t1",
(4,14)*{t}="t2",
(-3,3) *{\cdots}="dot1",
(12,3) *{\cdots}="dot2",
(0,-4) *{}="-1",
(8,-4)*{}="-2",
(4,-6)*{\vdots}="dot3",
(0,0) *{}="1",
(8,0)*{}="2",
(8,8)*{}="3",
(0,8)*{}="4",
(8,12)*{}="5",
(0,12)*{}="6",
(8,16)*{}="7",
(0,16)*{}="8",
\ar@{-} "1";"-1"^{}
\ar@{-} "-2";"2"^{}
\ar@{-} "1";"2"^{}
\ar@{-} "1";"4"^{}
\ar@{-} "2";"3"^{}
\ar@{-} "3";"4"^{}
\ar@{-} "3";"5"^{}
\ar@{-} "4";"6"^{}
\ar@{-} "6";"5"^{}
\ar@{-} "7";"5"^{}
\ar@{-} "6";"8"^{}
\ar@{-} "7";"8"^{}
\end{xy}
\]
In $Y$, a slot is named $A$ as above. If $Y'$ is a proper Young wall then
the slot $A$ in $Y$ is said to be {\it double} $t$-{\it admissible}.
\item[(ii)]
Let $Y''$ be a wall
obtained by removing two $t$-blocks of shape (\ref{smpl2}) from a column of $Y$:
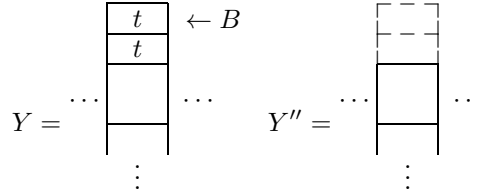
\[
Y=
\begin{xy}
(14,14) *{\leftarrow B}="B",
(4,10)*{t}="t1",
(4,14)*{t}="t2",
(-3,3) *{\cdots}="dot1",
(12,3) *{\cdots}="dot2",
(0,-4) *{}="-1",
(8,-4)*{}="-2",
(4,-6)*{\vdots}="dot3",
(0,0) *{}="1",
(8,0)*{}="2",
(8,8)*{}="3",
(0,8)*{}="4",
(8,12)*{}="5",
(0,12)*{}="6",
(8,16)*{}="7",
(0,16)*{}="8",
\ar@{-} "1";"-1"^{}
\ar@{-} "-2";"2"^{}
\ar@{-} "1";"2"^{}
\ar@{-} "1";"4"^{}
\ar@{-} "2";"3"^{}
\ar@{-} "3";"4"^{}
\ar@{-} "3";"5"^{}
\ar@{-} "4";"6"^{}
\ar@{-} "6";"5"^{}
\ar@{-} "7";"5"^{}
\ar@{-} "6";"8"^{}
\ar@{-} "7";"8"^{}
\end{xy}\quad
Y''=
\begin{xy}
(14,10) *{}="A",
(-3,3) *{\cdots}="dot1",
(12,3) *{\cdots}="dot2",
(0,-4) *{}="-1",
(8,-4)*{}="-2",
(4,-6)*{\vdots}="dot3",
(0,0) *{}="1",
(8,0)*{}="2",
(8,8)*{}="3",
(0,8)*{}="4",
(8,12)*{}="5",
(0,12)*{}="6",
(8,16)*{}="7",
(0,16)*{}="8",
\ar@{-} "1";"-1"^{}
\ar@{-} "-2";"2"^{}
\ar@{-} "1";"2"^{}
\ar@{-} "1";"4"^{}
\ar@{-} "2";"3"^{}
\ar@{-} "3";"4"^{}
\ar@{--} "3";"5"^{}
\ar@{--} "4";"6"^{}
\ar@{--} "6";"5"^{}
\ar@{--} "7";"5"^{}
\ar@{--} "6";"8"^{}
\ar@{--} "7";"8"^{}
\end{xy}
\]
A block in $Y$ is named $B$ as above.
If $Y''$ is a proper Young wall then the block $B$ in $Y$ is said to be {\it double} $t$-{\it removable}.
\item[(iii)] Other admissible slot (resp. removable block)
is said to be single admissible (resp. single removable).
\end{enumerate}
\end{defn}

\begin{ex}\label{ex-A2-1}
Let us consider the Young wall of type ${\rm A}^{(2)}_{4}$ in Example \ref{ex-A2-0}.
By $h_1=4$, $h_2=2$, $h_3=1$ and $h_l=\frac{1}{2}$ ($l\geq4$), 
the first, second and third columns are full columns, but the fourth column is not since its height is half-unit.
We see that the Young wall is proper.
Admissible slots and removable blocks are as follows: 
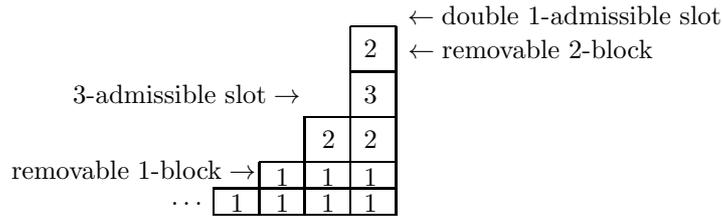
\begin{equation*}
\begin{xy}
(-15.5,-2) *{\ 1}="000-1",
(-9.5,-2) *{\ 1}="00-1",
(-3.5,-2) *{\ 1}="0-1",
(1.5,-2) *{\ \ 1}="0-1",
(-21.5,-2) *{\dots}="00000",
(-9.5,1.5) *{\ 1}="000",
(-3.5,1.5) *{\ 1}="00",
(1.5,1.5) *{\ \ 1}="0",
(1.5,6.5) *{\ \ 2}="02",
(-4,6.5) *{\ \ 2}="002",
(1.5,12.5) *{\ \ 3}="03",
(1.5,18.5) *{\ \ 2}="04",
(22.5,18.5) *{\ \ \leftarrow\text{removable }2\text{-block}}="arrow1",
(27,23) *{\ \ \leftarrow\text{double }1\text{-admissible slot}}="arrow2",
(-21.5,12.5) *{3\text{-admissible slot}\rightarrow}="arrow3",
(-28.5,2.5) *{\text{removable }1\text{-block}\rightarrow}="arrow4",
(0,0) *{}="1",
(0,-3.5) *{}="1-u",
(6,0)*{}="2",
(6,-3.5)*{}="2-u",
(6,3.5)*{}="3",
(0,3.5)*{}="4",
(-6,3.5)*{}="5",
(-6,0)*{}="6",
(-6,-3.5)*{}="6-u",
(-12,3.5)*{}="7",
(-12,0)*{}="8",
(-12,-3.6)*{}="8-u",
(-18,3.5)*{}="9",
(-18,0)*{}="10-a",
(-18,0)*{}="10",
(-18,-3.5)*{}="10-u",
(6,9.5)*{}="3-1",
(6,15.5)*{}="3-2",
(6,21.5)*{}="3-3",
(0,9.5)*{}="4-1",
(0,15.5)*{}="4-2",
(0,21.5)*{}="4-3",
(-6,9.5)*{}="5-1",
\ar@{-} "8";"10-a"^{}
\ar@{-} "10-u";"10-a"^{}
\ar@{-} "10-u";"2-u"^{}
\ar@{-} "8";"8-u"^{}
\ar@{-} "6";"6-u"^{}
\ar@{-} "2";"2-u"^{}
\ar@{-} "1";"1-u"^{}
\ar@{-} "1";"2"^{}
\ar@{-} "1";"4"^{}
\ar@{-} "2";"3"^{}
\ar@{-} "3";"4"^{}
\ar@{-} "5";"6"^{}
\ar@{-} "5";"4"^{}
\ar@{-} "1";"6"^{}
\ar@{-} "7";"8"^{}
\ar@{-} "7";"5"^{}
\ar@{-} "6";"8"^{}
\ar@{-} "3";"3-1"^{}
\ar@{-} "4-1";"3-1"^{}
\ar@{-} "4-1";"4"^{}
\ar@{-} "5-1";"5"^{}
\ar@{-} "4-1";"5-1"^{}
\ar@{-} "4-1";"4-2"^{}
\ar@{-} "3-2";"3-1"^{}
\ar@{-} "3-2";"4-2"^{}
\ar@{-} "3-2";"3-3"^{}
\ar@{-} "4-3";"4-2"^{}
\ar@{-} "4-3";"3-3"^{}
\end{xy}
\end{equation*}
It has only one double $1$-admissible slot and other admissible slots and removable blocks are single.
\end{ex}

In \cite{Kang}, the notion of {\it reducedness} on proper Young walls is defined and
it is proved that the set of all reduced proper Young walls of ground state $\lambda$
has an affine crystal structure and
is isomorphic to crystal graph $B(\lambda)$ of irreducible highest weight representation
with dominant integral highest weight $\lambda$ of level $1$ of $U_q(\mathfrak{g})$. 
In this article, however, we will not use the reducedness. According to [Corollary 2.5, \cite{KK}],
the set of all proper Young walls of ground state $\lambda$ is isomorphic to the direct sum of
$B(\lambda-m\delta)$ $(m\in\mathbb{Z}_{\geq0})$ with some multiplicities
(we will use this set in Theorem \ref{thmA2}, \ref{thmC1}).

\section{Combinatorial descriptions of inequalities}

\subsection{Setting and notation}\label{seno}

\begin{defn}\label{adapt}\cite{KaN}
Let $A=(a_{i,j})$ be the symmetrizable generalized Cartan matrix of $\mathfrak{g}$.
We suppose that a sequence $\io=(\cdots,i_3,i_2,i_1)$
satisfies $(\ref{seq-con})$.
We say $\iota$ is {\it adapted} to $A$ if the following condition holds : 
For $i,j\in I$ with $a_{i,j}<0$, the subsequence of $\iota$ consisting of all $i$, $j$ is
\[
(\cdots,i,j,i,j,i,j,i,j)\quad {\rm or}\quad (\cdots,j,i,j,i,j,i,j,i).
\]
In the case the Cartan matrix is fixed, the sequence $\iota$ is shortly said to be {\it adapted}.
\end{defn}

\begin{ex}
We consider the case
 $\mathfrak{g}$ is of type ${\rm A}^{(1)}_2$, $\iota=(\cdots,2,1,3,2,1,3,2,1,3)$.
It holds $a_{1,2}=a_{2,3}=a_{1,3}=-1$. 
\begin{itemize}
\item
The subsequence consisting of $1$, $2$ is $(\cdots,2,1,2,1,2,1)$.
\item
The subsequence consisting of $2$, $3$ is $(\cdots,2,3,2,3,2,3)$.
\item
The subsequence consisting of $1$, $3$ is $(\cdots,1,3,1,3,1,3)$.
\end{itemize}

Hence $\iota$ is an adapted sequence.
\end{ex}

\hspace{-7mm}
In the rest of article, we fix a sequence $\iota=(\cdots,i_3,i_2,i_1)$ which is adapted to $A$.
Let $(p_{i,j})_{i,j\in I;a_{i,j}<0}$ be the set of integers such that
\begin{equation}\label{pij}
p_{i,j}=\begin{cases}
1 & {\rm if}\ {\rm the\ subsequence\ of\ }\iota{\rm\ consisting\ of}\ i,j\ {\rm is}\ (\cdots,j,i,j,i,j,i),\\
0 & {\rm if}\ {\rm the\ subsequence\ of\ }\iota{\rm\ consisting\ of}\ i,j\ {\rm is}\ (\cdots,i,j,i,j,i,j).
\end{cases}
\end{equation}
Note that if $a_{i,j}<0$ then
\begin{equation}\label{pij2}
p_{i,j}+p_{j,i}=1.
\end{equation}
We identify each single index $j\in\mathbb{Z}_{\geq1}$ with a double index $(s,l)\in \mathbb{Z}_{\geq1}\times I$
if $i_j=l$ and $l$ is appearing $s$ times in $i_j$, $i_{j-1}$, $\cdots,i_1$.
For example, in the case $\iota=(\cdots,2,1,3,2,1,3,2,1,3)$, single indices
$\cdots,6,5,4,3,2,1$ are identified with 
\[
\cdots,(2,2),(2,1),(2,3),(1,2),(1,1),(1,3).
\]
The notation $x_{j}$, $\beta_{j}$ and $S'_{j}$ in the subsection \ref{poly-uqm} are also rewritten as
\[
x_j=x_{s,l},\quad \beta_j=\beta_{s,l},\quad S'_j=S'_{s,l}.
\]
By this identification and the ordinary order on $\mathbb{Z}_{\geq1}$, that is, $1<2<3<4<5<6<\cdots$,
we define an order on $\mathbb{Z}_{\geq1}\times I$. In the case $\iota=(\cdots,2,1,3,2,1,3,2,1,3)$, the order is
\[
\cdots>(2,2)>(2,1)>(2,3)>(1,2)>(1,1)>(1,3).
\]
Using the notation in (\ref{pij}), one can verify
\begin{equation}\label{pij3}
\beta_{s,l}=x_{s,l}+x_{s+1,l}+\sum_{j\in I; a_{l,j}<0} a_{l,j}x_{s+p_{j,l},j}.
\end{equation}
The above
$\beta_{s,l}$ is regarded as an analog of
simple root $\alpha_l$ of $^L\mathfrak{g}$
since $\alpha_l$ is expressed by 
\[
\alpha_l=\Lambda_l+\Lambda_l+\sum_{j\in I; a_{l,j}<0} a_{l,j}\Lambda_{j}
\]
on $\bigoplus_{i\in I}\mathbb{Z}h_i$.
It is a reason inequalities defining ${\rm Im}(\Psi_{\iota})$
is expressed by combinatorial objects related to $^L\mathfrak{g}$ rather than $\mathfrak{g}$
in Theorem \ref{thmA1}, \ref{thmD2}, \ref{thmA2}, \ref{thmC1}.

\subsection{Type ${\rm A}_{n-1}^{(1)}$-case}

For $k\in I$ and $t\in\mathbb{Z}$, let $P^k(t)\in\mathbb{Z}_{\geq0}$ be the non-negative integer defined as follows:
We set $P^k(k):=0$ and inductively define as
\begin{equation}\label{A1pk1}
P^k(t):=P^k(t-1)+p_{\ovl{t},\ovl{t-1}}\ (\text{for } t>k),
\end{equation}
\begin{equation}\label{A1pk2}
P^k(t):=P^k(t+1)+p_{\ovl{t},\ovl{t+1}}\ (\text{for } t<k),
\end{equation}
where the notation $\ovl{t}$ is defined in Definition \ref{pi1-def} (i). 
For an integer point $(i,j)$, $s\in\mathbb{Z}_{\geq1}$ and $k\in I$, we put
\begin{equation}\label{ovlL1}
\overline{L}_{s,k,\iota}(i,j):=x_{s+P^k(i+j)+{\rm min}\{k-j,i\},\overline{i+j}}\in (\mathbb{Q}^{\infty})^*.
\end{equation}
Note that if $(i,j)$ is a corner of an extended Young diagram with $y_{\infty}=k$ then by $i\geq0$ and $j\leq k$, it holds
\begin{equation}\label{pos-A1}
s+P^k(i+j)+{\rm min}\{k-j,i\}\geq s \geq1.
\end{equation}
For an extended Young diagram $T$ with $y_{\infty}=k$, one define
\begin{equation}\label{ovlL2}
\overline{L}_{s,k,\iota}(T)
:=\sum_{P:{\rm concave\ corner\ of\ }T} \overline{L}_{s,k,\iota}(P) - \sum_{P:{\rm convex\ corner\ of\ }T} \overline{L}_{s,k,\iota}(P)
\in (\mathbb{Q}^{\infty})^*.
\end{equation}
Let ${\rm EYD}_{k}$ be the set of extended Young diagrams with $y_{\infty}=k$ for $k\in I$.

\begin{thm}\label{thmA1}
If $\mathfrak{g}$ is of type ${\rm A}_{n-1}^{(1)}$ $(n\geq2)$ and $\iota$ is adapted, 
then $\iota$ satisfies the $\Xi'$-positivity condition and
\[
{\rm Im}(\Psi_{\iota})
=\left\{ 
\textbf{a}\in\mathbb{Z}^{\infty} \left|
\begin{array}{l}
\text{for any }k\in I,\ s\in\mathbb{Z}_{\geq1}\\
\text{and}\ T\in{\rm EYD}_{k},\ \text{it holds}\ 
\overline{L}_{s,k,\iota}(T)(\textbf{a})\geq0
\end{array} \right.
\right\}.
\]
\end{thm}

\begin{ex}
Considering the case $\mathfrak{g}$ is of type ${\rm A}_2^{(1)}$ and $\iota=(\cdots,3,1,2,3,1,2)$, let us
compute a part of inequalities defining ${\rm Im}(\Psi_{\iota})$. Note that $\iota$ is adapted.
We take an arbitrary $s\in\mathbb{Z}_{\geq1}$.
For $k=1,2,3$, the following diagrams are elements in ${\rm EYD}_{k}$:
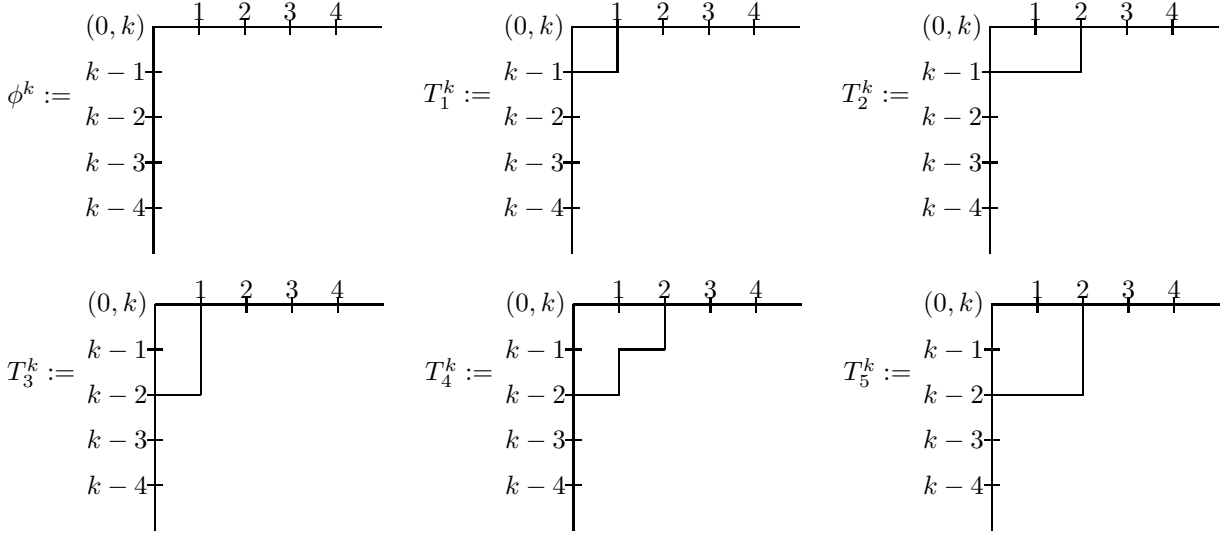
\[
\begin{xy}
(-15,-9) *{\phi^k:=}="Y0",
(0,0) *{}="1",
(30,0)*{}="2",
(0,-30)*{}="3",
(-5,0)*{(0,k)}="4",
(6,2) *{1}="10",
(6,-1) *{}="1010",
(12,2) *{2}="11",
(12,-1) *{}="1111",
(18,2) *{3}="12",
(18,-1) *{}="1212",
(24,2) *{4}="13",
(24,-1) *{}="1313",
(-5,-6)*{k-1}="5",
(1,-6)*{}="55",
(-5,-12)*{k-2}="6",
(1,-12)*{}="66",
(-5,-18)*{k-3}="7",
(1,-18)*{}="77",
(-5,-24)*{k-4}="8",
(1,-24)*{}="88",
(40,-9) *{T_1^k:=}="Y1",
(55,0) *{}="1a",
(85,0)*{}="2a",
(55,-30)*{}="3a",
(50,0)*{(0,k)}="4a",
(61,2) *{1}="10a",
(61,-1) *{}="1010a",
(61,-6) *{}="1-0a",
(67,2) *{2}="11a",
(67,-1) *{}="1111a",
(73,2) *{3}="12a",
(73,-1) *{}="1212a",
(79,2) *{4}="13a",
(79,-1) *{}="1313a",
(50,-6)*{k-1}="5a",
(56,-6)*{}="55a",
(50,-12)*{k-2}="6a",
(56,-12)*{}="66a",
(50,-18)*{k-3}="7a",
(56,-18)*{}="77a",
(50,-24)*{k-4}="8a",
(56,-24)*{}="88a",
(95,-9) *{T_2^k:=}="Y2",
(110,0) *{}="1b",
(140,0)*{}="2b",
(110,-30)*{}="3b",
(105,0)*{(0,k)}="4b",
(116,2) *{1}="10b",
(116,-1) *{}="1010b",
(122,-6) *{}="1-0b",
(122,2) *{2}="11b",
(122,-1) *{}="1111b",
(128,2) *{3}="12b",
(128,-1) *{}="1212b",
(134,2) *{4}="13b",
(134,-1) *{}="1313b",
(105,-6)*{k-1}="5b",
(111,-6)*{}="55b",
(105,-12)*{k-2}="6b",
(111,-12)*{}="66b",
(105,-18)*{k-3}="7b",
(111,-18)*{}="77b",
(105,-24)*{k-4}="8b",
(111,-24)*{}="88b",
\ar@{-} "1";"2"^{}
\ar@{-} "1";"3"^{}
\ar@{-} "5";"55"^{}
\ar@{-} "6";"66"^{}
\ar@{-} "7";"77"^{}
\ar@{-} "8";"88"^{}
\ar@{-} "10";"1010"^{}
\ar@{-} "11";"1111"^{}
\ar@{-} "12";"1212"^{}
\ar@{-} "13";"1313"^{}
\ar@{-} "1a";"2a"^{}
\ar@{-} "1a";"3a"^{}
\ar@{-} "5a";"55a"^{}
\ar@{-} "6a";"66a"^{}
\ar@{-} "7a";"77a"^{}
\ar@{-} "8a";"88a"^{}
\ar@{-} "10a";"1010a"^{}
\ar@{-} "11a";"1111a"^{}
\ar@{-} "12a";"1212a"^{}
\ar@{-} "13a";"1313a"^{}
\ar@{-} "55a";"1-0a"^{}
\ar@{-} "1010a";"1-0a"^{}
\ar@{-} "1b";"2b"^{}
\ar@{-} "1b";"3b"^{}
\ar@{-} "5b";"55b"^{}
\ar@{-} "6b";"66b"^{}
\ar@{-} "7b";"77b"^{}
\ar@{-} "8b";"88b"^{}
\ar@{-} "10b";"1010b"^{}
\ar@{-} "11b";"1111b"^{}
\ar@{-} "12b";"1212b"^{}
\ar@{-} "13b";"1313b"^{}
\ar@{-} "55b";"1-0b"^{}
\ar@{-} "11b";"1-0b"^{}
\end{xy}
\]
\[
\begin{xy}
(-15,-9) *{T_3^k:=}="Y3",
(0,0) *{}="1",
(30,0)*{}="2",
(0,-30)*{}="3",
(-5,0)*{(0,k)}="4",
(6,2) *{1}="10",
(6,-1) *{}="1010",
(12,2) *{2}="11",
(12,-1) *{}="1111",
(18,2) *{3}="12",
(18,-1) *{}="1212",
(24,2) *{4}="13",
(24,-1) *{}="1313",
(-5,-6)*{k-1}="5",
(1,-6)*{}="55",
(6,-12)*{}="6-12",
(-5,-12)*{k-2}="6",
(1,-12)*{}="66",
(-5,-18)*{k-3}="7",
(1,-18)*{}="77",
(-5,-24)*{k-4}="8",
(1,-24)*{}="88",
(40,-9) *{T_4^k:=}="Y4",
(55,0) *{}="1a",
(85,0)*{}="2a",
(55,-30)*{}="3a",
(50,0)*{(0,k)}="4a",
(61,2) *{1}="10a",
(61,-1) *{}="1010a",
(61,-6) *{}="6-6a",
(61,-12)*{}="6-12a",
(67,-6)*{}="12-6a",
(67,2) *{2}="11a",
(67,-1) *{}="1111a",
(73,2) *{3}="12a",
(73,-1) *{}="1212a",
(79,2) *{4}="13a",
(79,-1) *{}="1313a",
(50,-6)*{k-1}="5a",
(56,-6)*{}="55a",
(50,-12)*{k-2}="6a",
(56,-12)*{}="66a",
(50,-18)*{k-3}="7a",
(56,-18)*{}="77a",
(50,-24)*{k-4}="8a",
(56,-24)*{}="88a",
(95,-9) *{T_5^k:=}="Y5",
(110,0) *{}="1b",
(140,0)*{}="2b",
(110,-30)*{}="3b",
(105,0)*{(0,k)}="4b",
(116,2) *{1}="10b",
(116,-1) *{}="1010b",
(122,-12) *{}="1-0b",
(122,2) *{2}="11b",
(122,-1) *{}="1111b",
(128,2) *{3}="12b",
(128,-1) *{}="1212b",
(134,2) *{4}="13b",
(134,-1) *{}="1313b",
(105,-6)*{k-1}="5b",
(111,-6)*{}="55b",
(105,-12)*{k-2}="6b",
(111,-12)*{}="66b",
(105,-18)*{k-3}="7b",
(111,-18)*{}="77b",
(105,-24)*{k-4}="8b",
(111,-24)*{}="88b",
\ar@{-} "1";"2"^{}
\ar@{-} "1";"3"^{}
\ar@{-} "5";"55"^{}
\ar@{-} "6";"66"^{}
\ar@{-} "6";"6-12"^{}
\ar@{-} "1010";"6-12"^{}
\ar@{-} "7";"77"^{}
\ar@{-} "8";"88"^{}
\ar@{-} "10";"1010"^{}
\ar@{-} "11";"1111"^{}
\ar@{-} "12";"1212"^{}
\ar@{-} "13";"1313"^{}
\ar@{-} "1a";"2a"^{}
\ar@{-} "1a";"3a"^{}
\ar@{-} "5a";"55a"^{}
\ar@{-} "6a";"66a"^{}
\ar@{-} "7a";"77a"^{}
\ar@{-} "8a";"88a"^{}
\ar@{-} "10a";"1010a"^{}
\ar@{-} "11a";"1111a"^{}
\ar@{-} "12a";"1212a"^{}
\ar@{-} "13a";"1313a"^{}
\ar@{-} "6-12a";"66a"^{}
\ar@{-} "6-12a";"6-6a"^{}
\ar@{-} "12-6a";"6-6a"^{}
\ar@{-} "12-6a";"1111a"^{}
\ar@{-} "1b";"2b"^{}
\ar@{-} "1b";"3b"^{}
\ar@{-} "5b";"55b"^{}
\ar@{-} "6b";"66b"^{}
\ar@{-} "7b";"77b"^{}
\ar@{-} "8b";"88b"^{}
\ar@{-} "10b";"1010b"^{}
\ar@{-} "11b";"1111b"^{}
\ar@{-} "12b";"1212b"^{}
\ar@{-} "13b";"1313b"^{}
\ar@{-} "6b";"1-0b"^{}
\ar@{-} "11b";"1-0b"^{}
\end{xy}
\]
The element $\phi^k$ has only one concave corner $(0,k)$ and no convex corner.
By (\ref{ovlL2}), it holds $\ovl{L}_{s,k,\iota}(\phi^k)=x_{s,k}$.
In $T^k_1$, the points $(1,k)$ and $(0,k-1)$ are concave corners and $(1,k-1)$ is a convex corner. Thus,
\[
\ovl{L}_{s,k,\iota}(T^k_1)=x_{s+P^k(k+1),\ovl{k+1}}+x_{s+P^k(k-1),\ovl{k-1}}-x_{s+1,k}.
\]
Similarly, it follows
\[
\ovl{L}_{s,k,\iota}(T^k_2)=x_{s+P^k(k+2),\ovl{k+2}}+x_{s+P^k(k-1),\ovl{k-1}}-x_{s+P^k(k+1)+1,\ovl{k+1}},
\]
\[
\ovl{L}_{s,k,\iota}(T^k_3)=x_{s+P^k(k+1),\ovl{k+1}}+x_{s+P^k(k-2),\ovl{k-2}}-x_{s+P^k(k-1)+1,\ovl{k-1}},
\]
\[
\ovl{L}_{s,k,\iota}(T^k_4)=x_{s+P^k(k-2),\ovl{k-2}}+x_{s+1,k}+x_{s+P^k(k+2),\ovl{k+2}}-x_{s+P^k(k-1)+1,\ovl{k-1}}-x_{s+P^k(k+1)+1,\ovl{k+1}},
\]
\[
\ovl{L}_{s,k,\iota}(T^k_5)=x_{s+P^k(k-2),\ovl{k-2}}+x_{s+P^k(k+2),\ovl{k+2}}-x_{s+2,k}.
\]
We obtain
\[
\cdots,\ P^1(-1)=1,\ P^1(0)=0,\ P^1(1)=0,\ P^1(2)=1,\ P^1(3)=1,\cdots,
\]
\[
\cdots,\ P^2(0)=0,\ P^2(1)=0,\ P^2(2)=0,\ P^2(3)=0,\ P^2(4)=1,\cdots,
\]
\[
\cdots,\ P^3(1)=1,\ P^3(2)=1,\ P^3(3)=0,\ P^3(4)=1,\ P^3(5)=2,\cdots .
\]
Thus,
\[
\ovl{L}_{s,1,\iota}(\phi^1)=x_{s,1},\quad
\ovl{L}_{s,1,\iota}(T^1_1)=x_{s+1,2}+x_{s,3}-x_{s+1,1},\quad
\ovl{L}_{s,1,\iota}(T^1_2)=x_{s+1,3}+x_{s,3}-x_{s+2,2},
\]
\[
\ovl{L}_{s,1,\iota}(T^1_3)=2x_{s+1,2}-x_{s+1,3},\quad
\ovl{L}_{s,1,\iota}(T^1_4)=x_{s+1,2}+x_{s+1,1}-x_{s+2,2},\quad
\ovl{L}_{s,1,\iota}(T^1_5)=x_{s+1,2}+x_{s+1,3}-x_{s+2,1},
\]
\[
\ovl{L}_{s,2,\iota}(\phi^2)=x_{s,2},\quad
\ovl{L}_{s,2,\iota}(T^2_1)=x_{s,1}+x_{s,3}-x_{s+1,2},\quad
\ovl{L}_{s,2,\iota}(T^2_2)=x_{s,1}+x_{s+1,1}-x_{s+1,3},
\]
\[
\ovl{L}_{s,2,\iota}(T^2_3)=2x_{s,3}-x_{s+1,1},\quad
\ovl{L}_{s,2,\iota}(T^2_4)=x_{s,3}+x_{s+1,2}-x_{s+1,3},\quad
\ovl{L}_{s,2,\iota}(T^2_5)=x_{s,3}+x_{s+1,1}-x_{s+2,2},
\]
\[
\ovl{L}_{s,3,\iota}(\phi^3)=x_{s,3},\quad
\ovl{L}_{s,3,\iota}(T^3_1)=x_{s+1,1}+x_{s+1,2}-x_{s+1,3},\quad
\ovl{L}_{s,3,\iota}(T^3_2)=x_{s+2,2}+x_{s+1,2}-x_{s+2,1},
\]
\[
\ovl{L}_{s,3,\iota}(T^3_3)=2x_{s+1,1}-x_{s+2,2},\quad
\ovl{L}_{s,3,\iota}(T^3_4)=x_{s+1,1}+x_{s+1,3}-x_{s+2,1},\quad
\ovl{L}_{s,3,\iota}(T^3_5)=x_{s+1,1}+x_{s+2,2}-x_{s+2,3}.
\]
By Theorem \ref{thmA1}, we get a part of inequalities defining ${\rm Im}(\Psi_{\iota})$:
\[
{\rm Im}(\Psi_{\iota})
=\left\{ 
\textbf{a}\in\mathbb{Z}^{\infty} \left| 
\begin{array}{l}
s\in\mathbb{Z}_{\geq1},\ a_{s,1}\geq0,
a_{s+1,2}+a_{s,3}-a_{s+1,1}\geq0,
a_{s+1,3}+a_{s,3}-a_{s+2,2}\geq0,\\
2a_{s+1,2}-a_{s+1,3}\geq0,
a_{s+1,2}+a_{s+1,1}-a_{s+2,2}\geq0,
a_{s+1,2}+a_{s+1,3}-a_{s+2,1}\geq0,\cdots\\
a_{s,2}\geq0,\ 
a_{s,1}+a_{s,3}-a_{s+1,2}\geq0,
a_{s,1}+a_{s+1,1}-a_{s+1,3}\geq0,\\
2a_{s,3}-a_{s+1,1}\geq0,
a_{s,3}+a_{s+1,2}-a_{s+1,3}\geq0,
a_{s,3}+a_{s+1,1}-a_{s+2,2}\geq0,\cdots\\
a_{s,3}\geq0,
a_{s+1,1}+a_{s+1,2}-a_{s+1,3}\geq0,
a_{s+2,2}+a_{s+1,2}-a_{s+2,1}\geq0,\\
2a_{s+1,1}-a_{s+2,2}\geq0,
a_{s+1,1}+a_{s+1,3}-a_{s+2,1}\geq0,
a_{s+1,1}+a_{s+2,2}-a_{s+2,3}\geq0,\cdots
\end{array}
\right.
\right\}.
\]
Since we considered only finitely many diagrams, the above inequalities are not all.
The omitted other inequalities are represented by `$\cdots$'.

\end{ex}

\subsection{Type ${\rm D}_{n}^{(2)}$-case}

For $k\in I$ and $t\in\mathbb{Z}$, let $P^k(t)\in\mathbb{Z}_{\geq0}$ be the non-negative integer defined as follows:
We set $P^k(k):=0$ and inductively define as
\[
P^k(t):=P^k(t-1)+p_{\pi(t),\pi(t-1)}\ (\text{for } t>k),
\]
\[
P^k(t):=P^k(t+1)+p_{\pi(t),\pi(t+1)}\ (\text{for } t<k),
\]
where $\pi$ is defined in Definition \ref{pi1-def} (ii).
For an integer point $(i,j)$, $s\in\mathbb{Z}_{\geq1}$ and $k\in I$, we put
\begin{equation}\label{LL1}
L_{s,k,\iota}(i,j):=x_{s+P^k(i+j)+{\rm min}\{k-j,i\},\pi(i+j)}
\end{equation}
just as in (\ref{ovlL1}).
Note that since the map $\pi$ is introduced in \cite{KMM} for the representation theory of type ${\rm C}^{(1)}_{n-1}$,
the assignment (\ref{LL1}) is related to ${\rm C}^{(1)}_{n-1}$ rather
than ${\rm D}_{n}^{(2)}$.
As checked in (\ref{pos-A1}),
if $(i,j)$ is a corner of an extended Young diagram with $y_{\infty}=k$ then it holds
\begin{equation*}
s+P^k(i+j)+{\rm min}\{k-j,i\}\geq s \geq1.
\end{equation*}
For an extended Young diagram $T$ with $y_{\infty}=k$, we set
\begin{equation}\label{LL2}
L_{s,k,\iota}(T):=
\sum_{P:{\rm concave\ corner\ of\ }T} L_{s,k,\iota}(P) - \sum_{P:{\rm convex\ corner\ of\ }T} L_{s,k,\iota}(P)\\
\in (\mathbb{Q}^{\infty})^*.
\end{equation}
Recall that we defined ${\rm EYD}_{k}$ as the set of extended Young diagrams with $y_{\infty}=k$ in the previous subsection.

\begin{thm}\label{thmD2}
If $\mathfrak{g}$ is of type ${\rm D}_n^{(2)}$ $(n\geq3)$ and $\iota$ is adapted, 
then $\iota$ satisfies the $\Xi'$-positivity condition and
\[
{\rm Im}(\Psi_{\iota})
=\left\{ 
\textbf{a}\in\mathbb{Z}^{\infty} \left|
\begin{array}{l}
\text{for any }k\in I,\ s\in\mathbb{Z}_{\geq1}\\
\text{and}\ T\in{\rm EYD}_{k},\ \text{it holds}\ 
L_{s,k,\iota}(T)(\textbf{a})\geq0
\end{array} \right.
\right\}.
\]
\end{thm}

\subsection{Type ${\rm A}_{2n-2}^{(2)}$-case}

\subsubsection{Assignment of inequalities to ${\rm REYD}_{{\rm A}^{(2)},k}$ $(k>1)$}

Let us fix an index $k\in I$ such that $k>1$.
For $t\in\mathbb{Z}$, let $P^k(t)\in\mathbb{Z}_{\geq0}$ be the non-negative integer defined as follows:
We set $P^k(k):=0$ and inductively define as
\[
P^k(t):=P^k(t-1)+p_{\pi_1(t),\pi_1(t-1)}\ (\text{for } t>k),
\]
\[
P^k(t):=P^k(t+1)+p_{\pi_1(t),\pi_1(t+1)}\ (\text{for } t<k),
\]
where we set $p_{1,1}=0$ and $\pi_1$ is defined in Definition \ref{pi1-def} (iii).
For $(i,j)\in \mathbb{Z}\times\mathbb{Z}$ and $s\in\mathbb{Z}_{\geq1}$, one defines
\[
L^1_{s,k,{\rm ad}}(i,j)=x_{s+P^k(i+k)+[i]_-+k-j,\pi_1(i+k)},\quad
L^1_{s,k,{\rm re}}(i,j)=x_{s+P^k(i+k-1)+[i-1]_-+k-j,\pi_1(i+k-1)},
\]
where $[i]_-={\rm min}\{i,0\}$. Note that if $(i,j)$ is admissible in ${\rm REYD}_{{\rm A}^{(2)},k}$ then
\begin{equation}\label{pos-A2-1}
s+P^k(i+k)+[i]_-+k-j\geq s\geq 1
\end{equation}
by (\ref{pos-lem1}). If $(i,j)$ is removable then by (\ref{pos-lem2}),
\begin{equation}\label{pos-A2-2}
s+P^k(i+k-1)+[i-1]_-+k-j\geq s+1\geq 2.
\end{equation}

\begin{defn}\label{ad-rem-pt2}
Let $T=(y_t)_{t\in\mathbb{Z}}$ be a sequence in ${\rm REYD}_{{\rm A}^{(2)},k}$ of Definition \ref{AEYD} and $i\in\mathbb{Z}$. 
\begin{enumerate}
\item We suppose that the point $(i,y_i)$ is admissible.
If $y_{i-1}<y_i=y_{i+1}$ and it holds either $i+k\equiv 1$ and $i<0$ or $i+k\equiv 0$ and $i>0$
then we say the point $(i,y_i)$ is a double $1$-admissible point. 
\item We suppose that the point $(i,y_{i-1})$ is removable. If $y_{i-2}=y_{i-1}<y_{i}$ and it holds
either
$i+k-1\equiv 1$ and $i>1$ or $i+k-1\equiv 0$ and $i<1$
then we say the point $(i,y_{i-1})$ is a double $1$-removable point. 
\item Other admissible (resp. removable) points $(i,y_i)$ (resp. $(i,y_{i-1})$) other than (i) (resp. (ii)) are said to be single $\pi_1(i+k)$-admissible (resp. $\pi_1(i+k-1)$-removable) points.
\end{enumerate}
Here, in (i) and (ii), the notation $a\equiv b$ means $a\equiv b$ (mod $2n-1$). 
\end{defn}

\begin{ex}
For example, let $n=3$, $k=2$ and $T=(y_i)_{i\in\mathbb{Z}}\in {\rm REYD}_{{\rm A}^{(2)},2}$ is as follows:
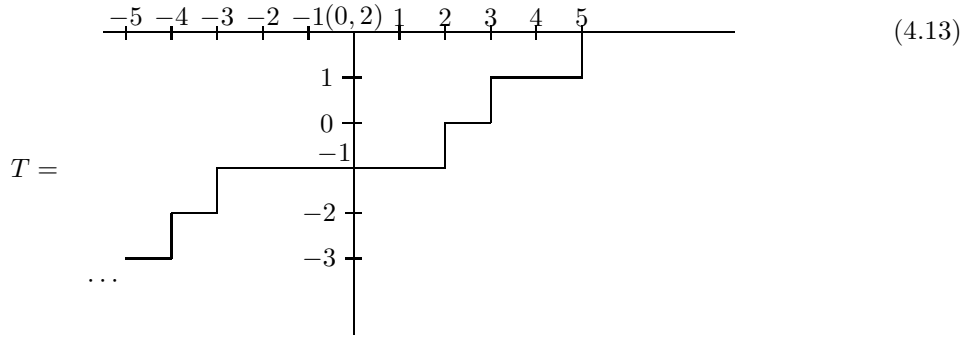
\begin{equation}\label{reydA-ex3}
\begin{xy}
(-42,-18) *{T=}="YY",
(-33,0) *{}="-6",
(-6,2) *{-1}="-1",
(-12,2) *{-2}="-2",
(-18,2) *{-3}="-3",
(-24,2) *{-4}="-4",
(-30,2) *{-5}="-5",
(-6,-1) *{}="-1a",
(-12,-1) *{}="-2a",
(-18,-1) *{}="-3a",
(-24,-1) *{}="-4a",
(-30,-1) *{}="-5a",
(0,0) *{}="1",
(50,0)*{}="2",
(0,-40)*{}="3",
(0,2)*{(0,2)}="4",
(6,2) *{1}="10",
(6,-1) *{}="1010",
(12,2) *{2}="11",
(12,-1) *{}="1111",
(18,2) *{3}="12",
(18,-1) *{}="1212",
(24,2) *{4}="13",
(24,-1) *{}="1313",
(30,2) *{5}="14",
(30,-1) *{}="1414",
(-3,-6)*{1\ }="5",
(1,-6)*{}="55",
(-3,-12)*{0\ }="6",
(-12,-18)*{}="6aaa",
(-18,-18)*{}="6aaaa",
(-18,-24)*{}="st1",
(-24,-24)*{}="st2",
(-24,-30)*{}="st3",
(-30,-30)*{}="st4",
(-33,-33)*{\cdots}="stdot",
(1,-12)*{}="66",
(-6,-18)*{}="7",
(-2,-16)*{-1\ }="7a",
(1,-18)*{}="77",
(-4,-24)*{-2\ }="8",
(1,-24)*{}="88",
(6,-18)*{}="8107",
(12,-18)*{}="810711",
(12,-12)*{}="8107116",
(18,-12)*{}="8107116a",
(18,-6)*{}="810711613",
(30,-6)*{}="8107116130",
(30,0)*{}="810711613000",
(-4,-30)*{-3\ }="9",
(1,-30)*{}="99",
\ar@{-} "1";"-6"^{}
\ar@{-} "1";"2"^{}
\ar@{-} "1";"3"^{}
\ar@{-} "-1";"-1a"^{}
\ar@{-} "-2";"-2a"^{}
\ar@{-} "-3";"-3a"^{}
\ar@{-} "-4";"-4a"^{}
\ar@{-} "-5";"-5a"^{}
\ar@{-} "5";"55"^{}
\ar@{-} "6";"66"^{}
\ar@{-} "8107";"7"^{}
\ar@{-} "8107";"810711"^{}
\ar@{-} "810711";"8107116"^{}
\ar@{-} "8107116a";"8107116"^{}
\ar@{-} "8107116a";"810711613"^{}
\ar@{-} "810711613";"8107116130"^{}
\ar@{-} "8107116130";"810711613000"^{}
\ar@{-} "7";"6aaa"^{}
\ar@{-} "6aaaa";"6aaa"^{}
\ar@{-} "st1";"6aaaa"^{}
\ar@{-} "st1";"st2"^{}
\ar@{-} "st2";"st3"^{}
\ar@{-} "st3";"st4"^{}
\ar@{-} "8";"88"^{}
\ar@{-} "9";"99"^{}
\ar@{-} "10";"1010"^{}
\ar@{-} "11";"1111"^{}
\ar@{-} "12";"1212"^{}
\ar@{-} "13";"1313"^{}
\ar@{-} "14";"1414"^{}
\end{xy}
\end{equation}
Thus, $y_l=l+2$ for $l\in\mathbb{Z}_{\leq-3}$, $y_{-2}=y_{-1}=y_0=y_1=-1$, $y_2=0$, $y_3=y_4=1$
and $y_l=2$ for $l\in\mathbb{Z}_{\geq 5}$.
Then the point $(5,1)$ is a double $1$-removable point.
The point $(3,1)$ is a double $1$-admissible point,
the points $(5,2)$ and $(-3,-1)$ are single $2$-admissible points.
The point $(2,-1)$ is a single $3$-removable point.
Note that the point $(-1,-1)$ is a single $1$-admissible point and also a single $1$-removable point.
In this way, it may happen that a point has both admissibility and removability.

\end{ex}

\nd
For each $T\in{\rm REYD}_{{\rm A}^{(2)},k}$, we set
\begin{eqnarray}
L^1_{s,k,\iota}(T)&:=&
\sum_{t\in I} \left(
\sum_{P:\text{single }t\text{-admissible point of }T} L^1_{s,k,{\rm ad}}(P)
-\sum_{P:\text{single }t\text{-removable point of }T} L^1_{s,k,{\rm re}}(P)\right)\label{L1kdef}\\
& & + 
\sum_{P:\text{double }1\text{-admissible point of }T} 2L^1_{s,k,{\rm ad}}(P)
-\sum_{P:\text{double }1\text{-removable point of }T} 2L^1_{s,k,{\rm re}}(P)\in (\mathbb{Q}^{\infty})^*.\nonumber
\end{eqnarray}
\nd
For instance, if $T$ is the element in (\ref{reydA-ex3}) then
\begin{eqnarray*}
L^1_{s,2,\iota}(T)
&=& L^1_{s,2,{\rm ad}}(5,2)+L^1_{s,2,{\rm ad}}(-3,-1) + L^1_{s,2,{\rm ad}}(-1,-1)
-L^1_{s,2,{\rm re}}(2,-1)-L^1_{s,2,{\rm re}}(-1,-1)\\
& &+2L^1_{s,2,{\rm ad}}(3,1)-2L^1_{s,2,{\rm re}}(5,1)\\
&=&x_{s+P^2(7),2}+x_{s+P^2(-1),2}+x_{s+P^2(1)+2,1}
-x_{s+P^2(3)+3,3}-x_{s+P^2(0)+1,1}\\
& & +2x_{s+P^2(5)+1,1}-2x_{s+P^2(6)+1,1}\\
&=&x_{s+P^2(7),2}+x_{s+P^2(-1),2}+x_{s+P^2(1)+2,1}
-x_{s+P^2(3)+3,3}-x_{s+P^2(0)+1,1}.
\end{eqnarray*}

\subsubsection{Assignment of inequalities to Young walls}

We draw Young walls on $\mathbb{R}_{\leq0}\times \mathbb{R}_{\geq1}$.
For example, the Young wall in Example \ref{ex-A2-1} is drawn as follow:
\begin{equation*}
\begin{xy}
(-15.5,-2) *{\ 1}="000-1",
(-9.5,-2) *{\ 1}="00-1",
(-3.5,-2) *{\ 1}="0-1",
(1.5,-2) *{\ \ 1}="0-1",
(-21.5,-2) *{\dots}="00000",
(-9.5,1.5) *{\ 1}="000",
(-3.5,1.5) *{\ 1}="00",
(1.5,1.5) *{\ \ 1}="0",
(1.5,6.5) *{\ \ 2}="02",
(-4,6.5) *{\ \ 2}="002",
(1.5,12.5) *{\ \ 3}="03",
(1.5,18.5) *{\ \ 2}="04",
(12,-6) *{(0,1)}="origin",
(0,0) *{}="1",
(0,-3.5) *{}="1-u",
(6,0)*{}="2",
(6,-3.5)*{}="2-u",
(6,-7.5)*{}="-1-u",
(6,3.5)*{}="3",
(0,3.5)*{}="4",
(-6,3.5)*{}="5",
(-6,0)*{}="6",
(-6,-3.5)*{}="6-u",
(-12,3.5)*{}="7",
(-12,0)*{}="8",
(-12,-3.6)*{}="8-u",
(-18,3.5)*{}="9",
(-18,0)*{}="10-a",
(-18,0)*{}="10",
(-18,-3.5)*{}="10-u",
(-24,-3.5)*{}="11-u",
(20,-3.5)*{}="0-u",
(6,9.5)*{}="3-1",
(6,15.5)*{}="3-2",
(6,21.5)*{}="3-3",
(6,25)*{}="3-4",
(0,9.5)*{}="4-1",
(0,15.5)*{}="4-2",
(0,21.5)*{}="4-3",
(-6,9.5)*{}="5-1",
\ar@{-} "-1-u";"2-u"^{}
\ar@{->} "2-u";"0-u"^{}
\ar@{-} "10-u";"11-u"^{}
\ar@{-} "8";"10-a"^{}
\ar@{-} "10-u";"10-a"^{}
\ar@{-} "10-u";"2-u"^{}
\ar@{-} "8";"8-u"^{}
\ar@{-} "6";"6-u"^{}
\ar@{-} "2";"2-u"^{}
\ar@{-} "1";"1-u"^{}
\ar@{-} "1";"2"^{}
\ar@{-} "1";"4"^{}
\ar@{-} "2";"3"^{}
\ar@{-} "3";"4"^{}
\ar@{-} "5";"6"^{}
\ar@{-} "5";"4"^{}
\ar@{-} "1";"6"^{}
\ar@{-} "7";"8"^{}
\ar@{-} "7";"5"^{}
\ar@{-} "6";"8"^{}
\ar@{-} "3";"3-1"^{}
\ar@{-} "4-1";"3-1"^{}
\ar@{-} "4-1";"4"^{}
\ar@{-} "5-1";"5"^{}
\ar@{-} "4-1";"5-1"^{}
\ar@{-} "4-1";"4-2"^{}
\ar@{-} "3-2";"3-1"^{}
\ar@{-} "3-2";"4-2"^{}
\ar@{-} "3-2";"3-3"^{}
\ar@{-} "4-3";"4-2"^{}
\ar@{-} "4-3";"3-3"^{}
\ar@{->} "3-3";"3-4"^{}
\end{xy}
\end{equation*}

\nd
Here, the unit length is $1$. 
Considering the map
$\{1,2,\cdots,2n-2\}\rightarrow \{1,2,\cdots,n\}$
defined as
\[
l\mapsto l,\ 2n-l\mapsto l \qquad (2\leq l\leq n-1),
\]
\[
1\mapsto1,\ n\mapsto n
\]
and extend it to a map 
\begin{equation}\label{piprime}
\pi':\mathbb{Z}_{\geq1}\rightarrow\{1,2,\cdots,n\}
\end{equation}
by periodicity $2n-2$.
We inductively define integers $P^1(l)$ ($l\in\mathbb{Z}_{\geq1}$) as
\[
P^1(1):=0,\ \ P^1(l)=P^1(l-1)+p_{\pi'(l),\pi'(l-1)}
\]
and fix an integer $s\in\mathbb{Z}_{\geq1}$.
Let $i\in\mathbb{Z}_{\geq0}$, $l\in\mathbb{Z}_{\geq1}$ and $S$ be a slot or block
\[
S=
\begin{xy}
(-8,15) *{(-i-1,l+1)}="0000",
(17,15) *{(-i,l+1)}="000",
(15,-3) *{(-i,l)}="00",
(-7,-3) *{(-i-1,l)}="0",
(0,0) *{}="1",
(12,0)*{}="2",
(12,12)*{}="3",
(0,12)*{}="4",
\ar@{-} "1";"2"^{}
\ar@{-} "1";"4"^{}
\ar@{-} "2";"3"^{}
\ar@{-} "3";"4"^{}
\end{xy}
\]
in $\mathbb{R}_{\leq0}\times \mathbb{R}_{\geq1}$.
If $S$ is colored by $t\in I\setminus \{1\}$ in the pattern of Definition \ref{def-YW} (ii)
then we set
\begin{equation}\label{l1def1}
L^1_{s,1,{\rm ad}}(S):=x_{s+P^1(l)+i,t},\quad L^1_{s,1,{\rm re}}(S):=x_{s+P^1(l)+i+1,t}.
\end{equation}
Let $i\in\mathbb{Z}_{\geq0}$, $l\in\mathbb{Z}_{\geq1}$ and $S'$ be
a slot or block colored by $t=1$ in $\mathbb{R}_{\leq0}\times \mathbb{R}_{\geq1}$
such that the place is one of the following two: 
\[
S'=
\begin{xy}
(-8,10) *{(-i-1,l+\frac{1}{2})}="0000",
(17,10) *{(-i,l+\frac{1}{2})}="000",
(15,-3) *{(-i,l)}="00",
(-7,-3) *{(-i-1,l)}="0",
(0,0) *{}="1",
(12,0)*{}="2",
(12,6)*{}="3",
(0,6)*{}="4",
(32,3) *{{\rm or}}="or",
(52,10) *{(-i-1,l+1)}="a0000",
(77,10) *{(-i,l+1)}="a000",
(75,-3) *{(-i,l+\frac{1}{2})}="a00",
(54,-3) *{(-i-1,l+\frac{1}{2})}="a0",
(60,0) *{}="a1",
(72,0)*{}="a2",
(72,6)*{}="a3",
(60,6)*{}="a4",
\ar@{-} "1";"2"^{}
\ar@{-} "1";"4"^{}
\ar@{-} "2";"3"^{}
\ar@{-} "3";"4"^{}
\ar@{-} "a1";"a2"^{}
\ar@{-} "a1";"a4"^{}
\ar@{-} "a2";"a3"^{}
\ar@{-} "a3";"a4"^{}
\end{xy}
\]
Then we set
\begin{equation}\label{l1def2}
L^1_{s,1,{\rm ad}}(S'):=x_{s+P^1(l)+i,1},\quad L^1_{s,1,{\rm re}}(S'):=x_{s+P^1(l)+i+1,1}.
\end{equation}
In these cases,
it is easy to see 
\begin{equation}\label{A2YW-pr}
s+P^1(l)+i\geq s\geq1,\quad s+P^1(l)+i+1\geq s+1\geq2.
\end{equation}
Considering the pattern in Definition \ref{def-YW} (ii), it holds
\begin{equation}\label{t-pi}
t=\pi'(l).
\end{equation}

For a proper Young wall $Y$ of type ${\rm A}^{(2)}_{2n-2}$ of ground state $\Lambda_1$, we define
\begin{eqnarray}
L^1_{s,1,\iota}(Y)&:=&
\sum_{t\in I} \left(
\sum_{P:\text{single }t\text{-admissible slot}} L^1_{s,1,{\rm ad}}(P)
-\sum_{P:\text{single removable }t\text{-block}} L^1_{s,1,{\rm re}}(P)\right)\nonumber\\
& & + 
\sum_{P:\text{double }1\text{-admissible slot}} 2L^1_{s,1,{\rm ad}}(P)
-\sum_{P:\text{double removable }1\text{-block}} 2L^1_{s,1,{\rm re}}(P). \label{L11-def}
\end{eqnarray}

Let ${\rm YW}_{{\rm A}^{(2)},1}$ be the set of all proper Young walls of type ${\rm A}^{(2)}_{2n-2}$ of ground state $\Lambda_1$.

\subsubsection{Combinatorial description of ${\rm Im}(\Psi_{\iota})$ of type ${\rm A}_{2n-2}^{(2)}$}

\begin{thm}\label{thmA2}
If $\mathfrak{g}$ is of type ${\rm A}_{2n-2}^{(2)}$ $(n\geq3)$ and $\iota$ is adapted 
then $\iota$ satisfies the $\Xi'$-positivity condition and
\[
{\rm Im}(\Psi_{\iota})
=\left\{ 
\textbf{a}\in\mathbb{Z}^{\infty} \left|
\begin{array}{l}
\text{for any }s\in\mathbb{Z}_{\geq1},\ k\in I\setminus\{1\}\\
\text{and}\ T\in{\rm REYD}_{{\rm A}^{(2)},k},\ \text{it holds}\ 
L^1_{s,k,\iota}(T)(\textbf{a})\geq0 \\
\text{and for any }Y\in{\rm YW}_{{\rm A}^{(2)},1},
\text{it holds}\ L^1_{s,1,\iota}(Y)(\textbf{a})\geq0
\end{array} \right.
\right\}.
\]
\end{thm}

\begin{ex}
Considering the case $\mathfrak{g}$ is of type ${\rm A}_4^{(2)}$ and $\iota=(\cdots,3,1,2,3,1,2)$, let us
compute a part of inequalities defining ${\rm Im}(\Psi_{\iota})$. We see that $\iota$ is adapted.
We get
\[
P^1(1)=0,\ P^1(2)=1,\ P^1(3)=1,\ P^1(4)=2,\cdots,
\]
\[
\cdots,\ P^2(0)=0,\ P^2(1)=0,\ P^2(2)=0,\ P^2(3)=0,\ P^2(4)=1,\cdots,
\]
\[
\cdots,\ P^3(1)=1,\ P^3(2)=1,\ P^3(3)=0,\ P^3(4)=1,\ P^3(5)=1,\cdots.
\]
We take an arbitrary $s\in\mathbb{Z}_{\geq1}$.
The following diagrams are elements in ${\rm REYD}_{{\rm A}^{(2)},2}$:
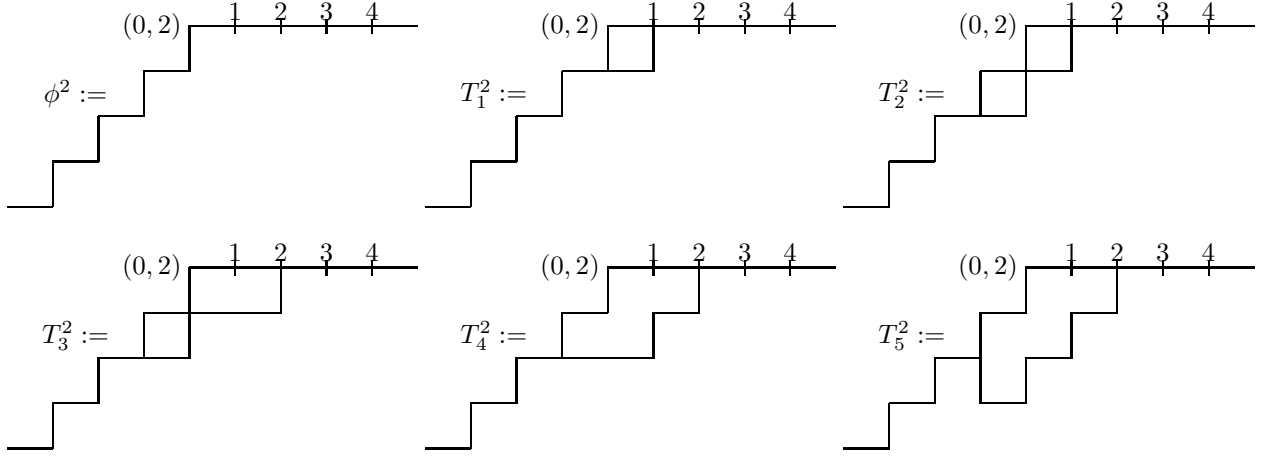
\[
\begin{xy}
(-15,-9) *{\phi^2:=}="Y0",
(0,0) *{}="1",
(30,0)*{}="2",
(-5,0)*{(0,2)}="4",
(6,2) *{1}="10",
(6,-1) *{}="1010",
(12,2) *{2}="11",
(12,-1) *{}="1111",
(18,2) *{3}="12",
(18,-1) *{}="1212",
(24,2) *{4}="13",
(24,-1) *{}="1313",
(0,-6)*{}="5",
(-6,-6)*{}="55",
(-6,-12)*{}="6",
(-12,-12)*{}="66",
(-12,-18)*{}="7",
(-18,-18)*{}="77",
(-18,-24)*{}="8",
(-24,-24)*{}="88",
(40,-9) *{T_1^2:=}="Y1",
(55,0) *{}="1a",
(85,0)*{}="2a",
(50,0)*{(0,2)}="4a",
(61,2) *{1}="10a",
(61,-1) *{}="1010a",
(61,-6) *{}="1010-1a",
(67,2) *{2}="11a",
(67,-1) *{}="1111a",
(73,2) *{3}="12a",
(73,-1) *{}="1212a",
(79,2) *{4}="13a",
(79,-1) *{}="1313a",
(55,-6)*{}="5a",
(49,-6)*{}="55a",
(49,-12)*{}="6a",
(43,-12)*{}="66a",
(43,-18)*{}="7a",
(37,-18)*{}="77a",
(37,-24)*{}="8a",
(31,-24)*{}="88a",
(95,-9) *{T_2^2:=}="Y2",
(110,0) *{}="1b",
(110,-12) *{}="1b-12",
(140,0)*{}="2b",
(105,0)*{(0,2)}="4b",
(116,2) *{1}="10b",
(116,-1) *{}="1010b",
(116,-6) *{}="1010-1b",
(122,2) *{2}="11b",
(122,-1) *{}="1111b",
(128,2) *{3}="12b",
(128,-1) *{}="1212b",
(134,2) *{4}="13b",
(134,-1) *{}="1313b",
(110,-6)*{}="5b",
(104,-6)*{}="55b",
(104,-12)*{}="6b",
(98,-12)*{}="66b",
(98,-18)*{}="7b",
(92,-18)*{}="77b",
(92,-24)*{}="8b",
(86,-24)*{}="88b",
\ar@{-} "1";"2"^{}
\ar@{-} "1";"5"^{}
\ar@{-} "5";"55"^{}
\ar@{-} "55";"6"^{}
\ar@{-} "6";"66"^{}
\ar@{-} "66";"7"^{}
\ar@{-} "7";"77"^{}
\ar@{-} "77";"8"^{}
\ar@{-} "8";"88"^{}
\ar@{-} "10";"1010"^{}
\ar@{-} "11";"1111"^{}
\ar@{-} "12";"1212"^{}
\ar@{-} "13";"1313"^{}
\ar@{-} "1a";"2a"^{}
\ar@{-} "1a";"5a"^{}
\ar@{-} "5a";"55a"^{}
\ar@{-} "55a";"6a"^{}
\ar@{-} "6a";"66a"^{}
\ar@{-} "66a";"7a"^{}
\ar@{-} "7a";"77a"^{}
\ar@{-} "77a";"8a"^{}
\ar@{-} "8a";"88a"^{}
\ar@{-} "10a";"1010a"^{}
\ar@{-} "1010a";"1010-1a"^{}
\ar@{-} "5a";"1010-1a"^{}
\ar@{-} "11a";"1111a"^{}
\ar@{-} "12a";"1212a"^{}
\ar@{-} "13a";"1313a"^{}
\ar@{-} "1b";"2b"^{}
\ar@{-} "1b";"1b-12"^{}
\ar@{-} "66b";"1b-12"^{}
\ar@{-} "1b";"5b"^{}
\ar@{-} "5b";"55b"^{}
\ar@{-} "55b";"6b"^{}
\ar@{-} "6b";"66b"^{}
\ar@{-} "66b";"7b"^{}
\ar@{-} "7b";"77b"^{}
\ar@{-} "77b";"8b"^{}
\ar@{-} "8b";"88b"^{}
\ar@{-} "10b";"1010b"^{}
\ar@{-} "1010b";"1010-1b"^{}
\ar@{-} "5b";"1010-1b"^{}
\ar@{-} "11b";"1111b"^{}
\ar@{-} "12b";"1212b"^{}
\ar@{-} "13b";"1313b"^{}
\end{xy}
\]
\[
\begin{xy}
(-15,-9) *{T_3^2:=}="Y3",
(0,0) *{}="1",
(30,0)*{}="2",
(-5,0)*{(0,2)}="4",
(6,2) *{1}="10",
(6,-1) *{}="1010",
(6,-6) *{}="1010-1",
(12,2) *{2}="11",
(12,-6) *{}="12-6",
(12,-1) *{}="1111",
(18,2) *{3}="12",
(18,-1) *{}="1212",
(24,2) *{4}="13",
(24,-1) *{}="1313",
(0,-6)*{}="5",
(0,-12) *{}="0-12",
(-6,-6)*{}="55",
(-6,-12)*{}="6",
(-6,-18)*{}="6-18",
(0,-18)*{}="0-18",
(-12,-12)*{}="66",
(-12,-18)*{}="7",
(-18,-18)*{}="77",
(-18,-24)*{}="8",
(-24,-24)*{}="88",
(40,-9) *{T_4^2:=}="Y4",
(55,0) *{}="1a",
(85,0)*{}="2a",
(50,0)*{(0,2)}="4a",
(61,2) *{1}="10a",
(61,-1) *{}="1010a",
(61,-6) *{}="6-6a",
(61,-12) *{}="6-12a",
(67,-6) *{}="1010-1a",
(67,2) *{2}="11a",
(67,-1) *{}="1111a",
(73,2) *{3}="12a",
(73,-1) *{}="1212a",
(79,2) *{4}="13a",
(79,-1) *{}="1313a",
(55,-6)*{}="5a",
(49,-6)*{}="55a",
(49,-12)*{}="6a",
(43,-12)*{}="66a",
(43,-18)*{}="7a",
(37,-18)*{}="77a",
(37,-24)*{}="8a",
(31,-24)*{}="88a",
(95,-9) *{T_5^2:=}="Y5",
(110,0) *{}="1b",
(140,0)*{}="2b",
(105,0)*{(0,2)}="4b",
(116,2) *{1}="10b",
(116,-1) *{}="1010b",
(116,-6) *{}="6-6b",
(116,-12) *{}="6-12b",
(116,-18) *{}="6-18b",
(110,-12) *{}="0-12b",
(110,-18) *{}="0-18b",
(104,-18) *{}="-6-18b",
(104,-12) *{}="-6-12b",
(122,-6) *{}="1010-1b",
(122,2) *{2}="11b",
(122,-1) *{}="1111b",
(128,2) *{3}="12b",
(128,-1) *{}="1212b",
(134,2) *{4}="13b",
(134,-1) *{}="1313b",
(110,-6)*{}="5b",
(104,-6)*{}="55b",
(104,-12)*{}="6b",
(98,-12)*{}="66b",
(98,-18)*{}="7b",
(92,-18)*{}="77b",
(92,-24)*{}="8b",
(86,-24)*{}="88b",
\ar@{-} "1";"2"^{}
\ar@{-} "1";"5"^{}
\ar@{-} "5";"55"^{}
\ar@{-} "55";"6"^{}
\ar@{-} "6";"66"^{}
\ar@{-} "5";"0-12"^{}
\ar@{-} "6";"0-12"^{}
\ar@{-} "66";"7"^{}
\ar@{-} "7";"77"^{}
\ar@{-} "77";"8"^{}
\ar@{-} "8";"88"^{}
\ar@{-} "10";"1010"^{}
\ar@{-} "11";"1111"^{}
\ar@{-} "11";"12-6"^{}
\ar@{-} "12";"1212"^{}
\ar@{-} "13";"1313"^{}
\ar@{-} "12-6";"1010-1"^{}
\ar@{-} "5";"1010-1"^{}
\ar@{-} "1a";"2a"^{}
\ar@{-} "1a";"5a"^{}
\ar@{-} "5a";"55a"^{}
\ar@{-} "55a";"6a"^{}
\ar@{-} "6a";"66a"^{}
\ar@{-} "66a";"7a"^{}
\ar@{-} "7a";"77a"^{}
\ar@{-} "77a";"8a"^{}
\ar@{-} "8a";"88a"^{}
\ar@{-} "10a";"1010a"^{}
\ar@{-} "1111a";"1010-1a"^{}
\ar@{-} "6-12a";"6-6a"^{}
\ar@{-} "6-12a";"6a"^{}
\ar@{-} "1010-1a";"6-6a"^{}
\ar@{-} "11a";"1111a"^{}
\ar@{-} "12a";"1212a"^{}
\ar@{-} "13a";"1313a"^{}
\ar@{-} "1b";"2b"^{}
\ar@{-} "1b";"5b"^{}
\ar@{-} "5b";"55b"^{}
\ar@{-} "55b";"6b"^{}
\ar@{-} "6b";"66b"^{}
\ar@{-} "66b";"7b"^{}
\ar@{-} "7b";"77b"^{}
\ar@{-} "77b";"8b"^{}
\ar@{-} "8b";"88b"^{}
\ar@{-} "10b";"1010b"^{}
\ar@{-} "1111b";"1010-1b"^{}
\ar@{-} "6-12b";"6-6b"^{}
\ar@{-} "6-12b";"0-12b"^{}
\ar@{-} "0-12b";"0-18b"^{}
\ar@{-} "-6-18b";"0-18b"^{}
\ar@{-} "-6-18b";"-6-12b"^{}
\ar@{-} "1010-1b";"6-6b"^{}
\ar@{-} "11b";"1111b"^{}
\ar@{-} "12b";"1212b"^{}
\ar@{-} "13b";"1313b"^{}
\end{xy}
\]
The point $(0,2)$ is a single $2$-admissible point in $\phi^2$
 and other points are neither admissible nor removable.
Thus, $L^1_{s,2,\iota}(\phi^2)=x_{s,2}$. In $T_1^2$, the point $(-1,1)$ is a double $1$-admissible point, $(1,2)$ is a single $3$-admissible point,
$(1,1)$ is a single $2$-removable point, which implies
\[
L^1_{s,2,\iota}(T_1^2)=2x_{s+P^2(1),1}+x_{s+P^2(3),3}-x_{s+1,2}=2x_{s,1}+x_{s,3}-x_{s+1,2}.
\]
Similarly, considering $P^2(4)=p_{3,2}+p_{2,3}=1$,
it holds
\[
L^1_{s,2,\iota}(T_2^2)=x_{s+P^2(0),1}+x_{s+P^2(3),3}-x_{s+P^2(1)+1,1}
=x_{s,1}+x_{s,3}-x_{s+1,1},\]
\[
L^1_{s,2,\iota}(T_3^2)=x_{s+P^2(0),1}+2x_{s+1,2}
-x_{s+P^2(3)+1,3}-x_{s+P^2(1)+1,1}
=x_{s,1}+2x_{s+1,2}
-x_{s+1,3}-x_{s+1,1},
\]
\[
L^1_{s,2,\iota}(T_4^2)=x_{s+P^2(1)+1,1}+x_{s+1,2}
+x_{s+P^2(0),1}-x_{s+2,2}
=x_{s+1,1}+x_{s+1,2}
+x_{s,1}-x_{s+2,2}.
\]
We remark that $(-1,0)$ is also an admissible point in $T_4^2$.
In $T_5^2$, the admissible points are $(2,2)$ and $(-2,0)$, removable point is $(0,-1)$ so that
\[
L^1_{s,2,\iota}(T_5^2)=x_{s+P^2(0),1}+x_{s+1,2}
-x_{s+P^2(1)+2,1}=x_{s,1}+x_{s+1,2}
-x_{s+2,1}.
\]
The following diagrams are elements in ${\rm REYD}_{{\rm A}^{(2)},3}$:
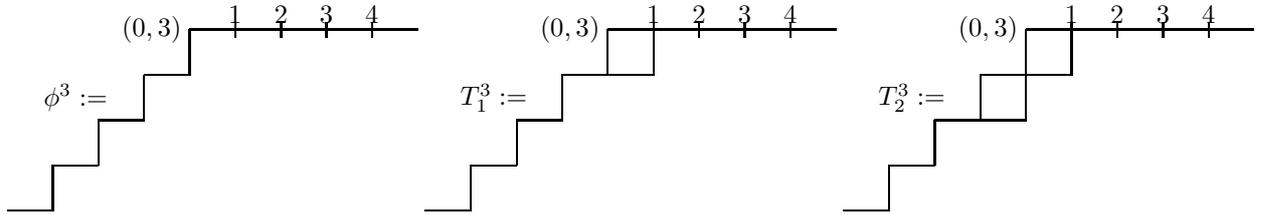
\[
\begin{xy}
(-15,-9) *{\phi^3:=}="Y0",
(0,0) *{}="1",
(30,0)*{}="2",
(-5,0)*{(0,3)}="4",
(6,2) *{1}="10",
(6,-1) *{}="1010",
(12,2) *{2}="11",
(12,-1) *{}="1111",
(18,2) *{3}="12",
(18,-1) *{}="1212",
(24,2) *{4}="13",
(24,-1) *{}="1313",
(0,-6)*{}="5",
(-6,-6)*{}="55",
(-6,-12)*{}="6",
(-12,-12)*{}="66",
(-12,-18)*{}="7",
(-18,-18)*{}="77",
(-18,-24)*{}="8",
(-24,-24)*{}="88",
(40,-9) *{T_1^3:=}="Y1",
(55,0) *{}="1a",
(85,0)*{}="2a",
(50,0)*{(0,3)}="4a",
(61,2) *{1}="10a",
(61,-1) *{}="1010a",
(61,-6) *{}="1010-1a",
(67,2) *{2}="11a",
(67,-1) *{}="1111a",
(73,2) *{3}="12a",
(73,-1) *{}="1212a",
(79,2) *{4}="13a",
(79,-1) *{}="1313a",
(55,-6)*{}="5a",
(49,-6)*{}="55a",
(49,-12)*{}="6a",
(43,-12)*{}="66a",
(43,-18)*{}="7a",
(37,-18)*{}="77a",
(37,-24)*{}="8a",
(31,-24)*{}="88a",
(95,-9) *{T_2^3:=}="Y2",
(110,0) *{}="1b",
(110,-12) *{}="1b-12",
(140,0)*{}="2b",
(105,0)*{(0,3)}="4b",
(116,2) *{1}="10b",
(116,-1) *{}="1010b",
(116,-6) *{}="1010-1b",
(122,2) *{2}="11b",
(122,-1) *{}="1111b",
(128,2) *{3}="12b",
(128,-1) *{}="1212b",
(134,2) *{4}="13b",
(134,-1) *{}="1313b",
(110,-6)*{}="5b",
(104,-6)*{}="55b",
(104,-12)*{}="6b",
(98,-12)*{}="66b",
(98,-18)*{}="7b",
(92,-18)*{}="77b",
(92,-24)*{}="8b",
(86,-24)*{}="88b",
\ar@{-} "1";"2"^{}
\ar@{-} "1";"5"^{}
\ar@{-} "5";"55"^{}
\ar@{-} "55";"6"^{}
\ar@{-} "6";"66"^{}
\ar@{-} "66";"7"^{}
\ar@{-} "7";"77"^{}
\ar@{-} "77";"8"^{}
\ar@{-} "8";"88"^{}
\ar@{-} "10";"1010"^{}
\ar@{-} "11";"1111"^{}
\ar@{-} "12";"1212"^{}
\ar@{-} "13";"1313"^{}
\ar@{-} "1a";"2a"^{}
\ar@{-} "1a";"5a"^{}
\ar@{-} "5a";"55a"^{}
\ar@{-} "55a";"6a"^{}
\ar@{-} "6a";"66a"^{}
\ar@{-} "66a";"7a"^{}
\ar@{-} "7a";"77a"^{}
\ar@{-} "77a";"8a"^{}
\ar@{-} "8a";"88a"^{}
\ar@{-} "10a";"1010a"^{}
\ar@{-} "1010a";"1010-1a"^{}
\ar@{-} "5a";"1010-1a"^{}
\ar@{-} "11a";"1111a"^{}
\ar@{-} "12a";"1212a"^{}
\ar@{-} "13a";"1313a"^{}
\ar@{-} "1b";"2b"^{}
\ar@{-} "1b";"1b-12"^{}
\ar@{-} "66b";"1b-12"^{}
\ar@{-} "1b";"5b"^{}
\ar@{-} "5b";"55b"^{}
\ar@{-} "55b";"6b"^{}
\ar@{-} "6b";"66b"^{}
\ar@{-} "66b";"7b"^{}
\ar@{-} "7b";"77b"^{}
\ar@{-} "77b";"8b"^{}
\ar@{-} "8b";"88b"^{}
\ar@{-} "10b";"1010b"^{}
\ar@{-} "1010b";"1010-1b"^{}
\ar@{-} "5b";"1010-1b"^{}
\ar@{-} "11b";"1111b"^{}
\ar@{-} "12b";"1212b"^{}
\ar@{-} "13b";"1313b"^{}
\end{xy}
\]
and considering $P^3(4)=p_{2,3}=P^3(2)$, it holds
\[
L^1_{s,3,\iota}(\phi^3)=x_{s,3},\quad
L^1_{s,3,\iota}(T_1^3)=2x_{s+P^3(2),2}-x_{s+1,3}=2x_{s+1,2}-x_{s+1,3},
\]
\[
L^1_{s,3,\iota}(T_2^3)=2x_{s+P^3(1),1}+x_{s+P^3(2),2}-x_{s+P^3(2)+1,2}
=2x_{s+1,1}+x_{s+1,2}-x_{s+2,2}.
\]
The following Young walls are elements in ${\rm YW}_{{\rm A}^{(2)},1}$:
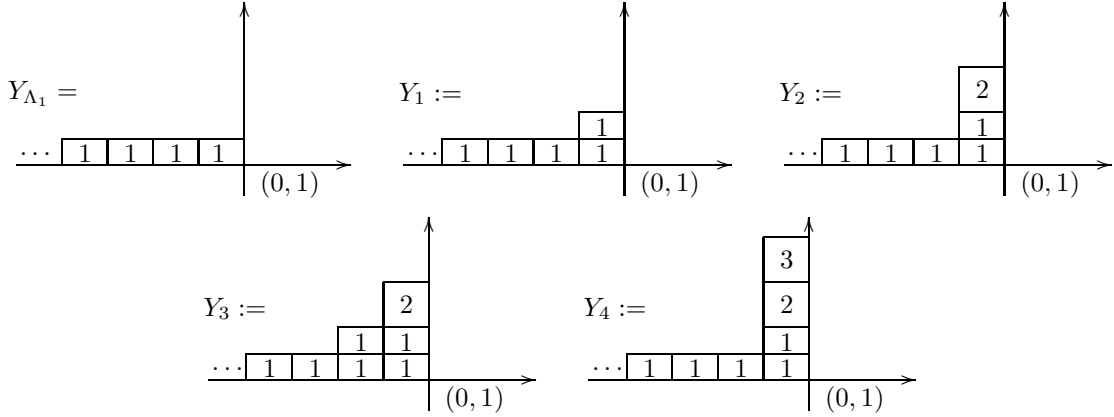
\[
\begin{xy}
(-20.5,6) *{Y_{\Lambda_1}=}="0000-1",
(-15.5,-2) *{\ 1}="000-1",
(-9.5,-2) *{\ 1}="00-1",
(-3.5,-2) *{\ 1}="0-1",
(1.5,-2) *{\ \ 1}="0-1s",
(-21.5,-2) *{\dots}="00000",
(12,-6) *{(0,1)}="origin",
(0,0) *{}="1",
(6,18) *{}="y",
(0,-3.5) *{}="1-u",
(6,0)*{}="2",
(6,-3.5)*{}="2-u",
(6,-7.5)*{}="-1-u",
(-6,0)*{}="6",
(-6,-3.5)*{}="6-u",
(-12,0)*{}="8",
(-12,-3.6)*{}="8-u",
(-18,3.5)*{}="9",
(-18,0)*{}="10-t",
(-18,0)*{}="10",
(-18,-3.5)*{}="10-u",
(-24,-3.5)*{}="11-u",
(20,-3.5)*{}="0-u",
(30.5,6) *{Y_1:=}="0000-1a",
(52.5,1.5) *{\ 1}="000-1-1a",
(34.5,-2) *{\ 1}="000-1a",
(40.5,-2) *{\ 1}="00-1a",
(46.5,-2) *{\ 1}="0-1a",
(52.5,-2) *{\ 1}="0-1sa",
(29.5,-2) *{\dots}="00000a",
(62,-6) *{(0,1)}="origina",
(50,0) *{}="1a",
(56,18) *{}="ya",
(56,3.5)*{}="6+35a",
(50,3.5)*{}="0+35a",
(50,0)*{}="0+0a",
(50,-3.5) *{}="1-ua",
(56,0)*{}="2a",
(56,-3.5)*{}="2-ua",
(56,-7.5)*{}="-1-ua",
(44,0)*{}="6a",
(44,-3.5)*{}="6-ua",
(38,0)*{}="8a",
(38,-3.6)*{}="8-ua",
(32,3.5)*{}="9a",
(32,0)*{}="10-ta",
(32,0)*{}="10ta",
(32,-3.5)*{}="10-uta",
(27,-3.5)*{}="11-uta",
(70,-3.5)*{}="0-ua",
(80.5,6) *{Y_2:=}="0000-1b",
(102.5,1.5) *{\ 1}="000-1-1b",
(102.5,6.5) *{\ 2}="000-1-1-2b",
(84.5,-2) *{\ 1}="000-1b",
(90.5,-2) *{\ 1}="00-1b",
(96.5,-2) *{\ 1}="0-1b",
(102.5,-2) *{\ 1}="0-1sb",
(79.5,-2) *{\dots}="00000b",
(112,-6) *{(0,1)}="originb",
(100,0) *{}="1b",
(106,18) *{}="yb",
(106,9.5)*{}="6+95b",
(100,9.5)*{}="0+95b",
(106,3.5)*{}="6+35b",
(100,3.5)*{}="0+35b",
(100,0)*{}="0+0b",
(100,-3.5) *{}="1-ub",
(106,0)*{}="2b",
(106,-3.5)*{}="2-ub",
(106,-7.5)*{}="-1-ub",
(94,0)*{}="6b",
(94,-3.5)*{}="6-ub",
(88,0)*{}="8b",
(88,-3.6)*{}="8-ub",
(82,3.5)*{}="9b",
(82,0)*{}="10-tb",
(82,0)*{}="10tb",
(82,-3.5)*{}="10-utb",
(77,-3.5)*{}="11-utb",
(120,-3.5)*{}="0-ub",
\ar@{-} "-1-u";"2-u"^{}
\ar@{->} "2-u";"0-u"^{}
\ar@{->} "2";"y"^{}
\ar@{-} "10-u";"11-u"^{}
\ar@{-} "8";"10-t"^{}
\ar@{-} "10-u";"10-t"^{}
\ar@{-} "10-u";"2-u"^{}
\ar@{-} "8";"8-u"^{}
\ar@{-} "6";"6-u"^{}
\ar@{-} "2";"2-u"^{}
\ar@{-} "1";"1-u"^{}
\ar@{-} "1";"2"^{}
\ar@{-} "1";"6"^{}
\ar@{-} "6";"8"^{}
\ar@{-} "-1-ua";"2-ua"^{}
\ar@{->} "2-ua";"0-ua"^{}
\ar@{->} "2a";"ya"^{}
\ar@{-} "10-uta";"11-uta"^{}
\ar@{-} "8a";"10-ta"^{}
\ar@{-} "10-uta";"10-ta"^{}
\ar@{-} "10-uta";"2-ua"^{}
\ar@{-} "8a";"8-ua"^{}
\ar@{-} "6a";"6-ua"^{}
\ar@{-} "2a";"2-ua"^{}
\ar@{-} "1a";"1-ua"^{}
\ar@{-} "1a";"2a"^{}
\ar@{-} "1a";"6a"^{}
\ar@{-} "6a";"8a"^{}
\ar@{-} "6+35a";"0+35a"^{}
\ar@{-} "0+0a";"0+35a"^{}
\ar@{-} "-1-ub";"2-ub"^{}
\ar@{->} "2-ub";"0-ub"^{}
\ar@{->} "2b";"yb"^{}
\ar@{-} "10-utb";"11-utb"^{}
\ar@{-} "8b";"10-tb"^{}
\ar@{-} "10-utb";"10-tb"^{}
\ar@{-} "10-utb";"2-ub"^{}
\ar@{-} "8b";"8-ub"^{}
\ar@{-} "6b";"6-ub"^{}
\ar@{-} "2b";"2-ub"^{}
\ar@{-} "1b";"1-ub"^{}
\ar@{-} "1b";"2b"^{}
\ar@{-} "1b";"6b"^{}
\ar@{-} "6b";"8b"^{}
\ar@{-} "6+35b";"0+35b"^{}
\ar@{-} "0+0b";"0+35b"^{}
\ar@{-} "0+95b";"6+95b"^{}
\ar@{-} "0+35b";"0+95b"^{}
\end{xy}
\]
\[
\begin{xy}
(30.5,6) *{Y_3:=}="0000-1a",
(52.5,1.5) *{\ 1}="000-1-1a",
(52.5,6.5) *{\ 2}="000-1-1-2a",
(34.5,-2) *{\ 1}="000-1a",
(40.5,-2) *{\ 1}="00-1a",
(46.5,-2) *{\ 1}="0-1a",
(46.5,1.5) *{\ 1}="0-1-1a",
(52.5,-2) *{\ 1}="0-1sa",
(29.5,-2) *{\dots}="00000a",
(62,-6) *{(0,1)}="origina",
(50,0) *{}="1a",
(56,18) *{}="ya",
(56,9.5)*{}="6+95a",
(50,9.5)*{}="0+95a",
(56,3.5)*{}="6+35a",
(50,3.5)*{}="0+35a",
(50,0)*{}="0+0a",
(50,-3.5) *{}="1-ua",
(56,0)*{}="2a",
(56,-3.5)*{}="2-ua",
(56,-7.5)*{}="-1-ua",
(44,3.5)*{}="-6+35a",
(44,0)*{}="6a",
(44,-3.5)*{}="6-ua",
(38,0)*{}="8a",
(38,-3.6)*{}="8-ua",
(32,3.5)*{}="9a",
(32,0)*{}="10-ta",
(32,0)*{}="10ta",
(32,-3.5)*{}="10-uta",
(27,-3.5)*{}="11-uta",
(70,-3.5)*{}="0-ua",
(80.5,6) *{Y_4:=}="0000-1b",
(102.5,1.5) *{\ 1}="000-1-1b",
(102.5,6.5) *{\ 2}="000-1-1-2b",
(102.5,12.5) *{\ 3}="000-1-1-2-3b",
(84.5,-2) *{\ 1}="000-1b",
(90.5,-2) *{\ 1}="00-1b",
(96.5,-2) *{\ 1}="0-1b",
(102.5,-2) *{\ 1}="0-1sb",
(79.5,-2) *{\dots}="00000b",
(112,-6) *{(0,1)}="originb",
(100,0) *{}="1b",
(106,18) *{}="yb",
(106,15.5)*{}="6+155b",
(100,15.5)*{}="0+155b",
(106,9.5)*{}="6+95b",
(100,9.5)*{}="0+95b",
(106,3.5)*{}="6+35b",
(100,3.5)*{}="0+35b",
(100,0)*{}="0+0b",
(100,-3.5) *{}="1-ub",
(106,0)*{}="2b",
(106,-3.5)*{}="2-ub",
(106,-7.5)*{}="-1-ub",
(94,0)*{}="6b",
(94,-3.5)*{}="6-ub",
(88,0)*{}="8b",
(88,-3.6)*{}="8-ub",
(82,3.5)*{}="9b",
(82,0)*{}="10-tb",
(82,0)*{}="10tb",
(82,-3.5)*{}="10-utb",
(77,-3.5)*{}="11-utb",
(120,-3.5)*{}="0-ub",
\ar@{-} "-1-ua";"2-ua"^{}
\ar@{->} "2-ua";"0-ua"^{}
\ar@{->} "2a";"ya"^{}
\ar@{-} "10-uta";"11-uta"^{}
\ar@{-} "8a";"10-ta"^{}
\ar@{-} "10-uta";"10-ta"^{}
\ar@{-} "10-uta";"2-ua"^{}
\ar@{-} "8a";"8-ua"^{}
\ar@{-} "6a";"6-ua"^{}
\ar@{-} "2a";"2-ua"^{}
\ar@{-} "1a";"1-ua"^{}
\ar@{-} "1a";"2a"^{}
\ar@{-} "1a";"6a"^{}
\ar@{-} "6a";"8a"^{}
\ar@{-} "-6+35a";"0+35a"^{}
\ar@{-} "-6+35a";"6a"^{}
\ar@{-} "6+35a";"0+35a"^{}
\ar@{-} "0+0a";"0+35a"^{}
\ar@{-} "0+95a";"6+95a"^{}
\ar@{-} "0+35a";"0+95a"^{}
\ar@{-} "-1-ub";"2-ub"^{}
\ar@{->} "2-ub";"0-ub"^{}
\ar@{->} "2b";"yb"^{}
\ar@{-} "10-utb";"11-utb"^{}
\ar@{-} "8b";"10-tb"^{}
\ar@{-} "10-utb";"10-tb"^{}
\ar@{-} "10-utb";"2-ub"^{}
\ar@{-} "8b";"8-ub"^{}
\ar@{-} "6b";"6-ub"^{}
\ar@{-} "2b";"2-ub"^{}
\ar@{-} "1b";"1-ub"^{}
\ar@{-} "1b";"2b"^{}
\ar@{-} "1b";"6b"^{}
\ar@{-} "6b";"8b"^{}
\ar@{-} "6+35b";"0+35b"^{}
\ar@{-} "0+0b";"0+35b"^{}
\ar@{-} "0+95b";"6+95b"^{}
\ar@{-} "0+35b";"0+95b"^{}
\ar@{-} "0+155b";"6+155b"^{}
\ar@{-} "0+95b";"0+155b"^{}
\end{xy}
\]
Since $Y_{\Lambda_1}$ has a single $1$-admissible slot, we have $L^1_{s,1,\iota}(Y_{\Lambda_1})=x_{s,1}$. $Y_1$ has a single $2$-admissible slot
and a single removable $1$-block so that $L^1_{s,1,\iota}(Y_1)=x_{s+P^1(2),2}-x_{s+1,1}=x_{s+1,2}-x_{s+1,1}$. Similarly,
\[
L^1_{s,1,\iota}(Y_2)=x_{s+P^1(3),3}+x_{s+1,1}-x_{s+P^1(2)+1,2}=x_{s+1,3}+x_{s+1,1}-x_{s+2,2}, \]
\[
L^1_{s,1,\iota}(Y_3)=x_{s+P^1(3),3}-x_{s+2,1}=x_{s+1,3}-x_{s+2,1},
\]
\[
L^1_{s,1,\iota}(Y_4)=x_{s+P^1(4),2}+x_{s+1,1}-x_{s+P^1(3)+1,3}=x_{s+2,2}+x_{s+1,1}-x_{s+2,3}.
\]
Thus, by Theorem \ref{thmA2}, we get a part of inequalities defining ${\rm Im}(\Psi_{\iota})$:
\[
{\rm Im}(\Psi_{\iota})
=\left\{ 
\textbf{a}\in\mathbb{Z}^{\infty} \left|
\begin{array}{l}
s\in\mathbb{Z}_{\geq1},\ a_{s,2}\geq0,
2a_{s,1}+a_{s,3}-a_{s+1,2}\geq0,
a_{s,1}+a_{s,3}-a_{s+1,1}\geq0,\\
a_{s,1}+2a_{s+1,2}-a_{s+1,3}-a_{s+1,1}\geq0,\\
a_{s,1}+a_{s+1,2}+a_{s+1,1}-a_{s+2,2}\geq0,
a_{s,1}+a_{s+1,2}-a_{s+2,1}\geq0,\cdots\\
a_{s,3}\geq0,\ 
2a_{s+1,2}-a_{s+1,3}\geq0,
2a_{s+1,1}+a_{s+1,2}-a_{s+2,2}\geq0,\cdots\\
a_{s,1}\geq0,
a_{s+1,2}-a_{s+1,1}\geq0,
a_{s+1,3}+a_{s+1,1}-a_{s+2,2}\geq0,\\
a_{s+1,3}-a_{s+2,1}\geq0,
a_{s+2,2}+a_{s+1,1}-a_{s+2,3}\geq0,\cdots
\end{array} \right.
\right\}.
\]

\end{ex}

\subsection{Type ${\rm C}_{n-1}^{(1)}$-case}

\subsubsection{Assignment of inequalities to ${\rm REYD}_{{\rm D}^{(2)},k}$ $(1<k< n)$}

Let us fix an index $k\in I$ such that $1<k<n$.
For $t\in\mathbb{Z}$, let $P^k(t)\in\mathbb{Z}_{\geq0}$ be the non-negative integer defined as follows:
We set $P^k(k):=0$ and inductively define as
\[
P^k(t):=P^k(t-1)+p_{\pi_2(t),\pi_2(t-1)}\ (\text{for } t>k),
\]
\[
P^k(t):=P^k(t+1)+p_{\pi_2(t),\pi_2(t+1)}\ (\text{for } t<k),
\]
where we set $p_{1,1}=0$, $p_{n,n}=0$ and $\pi_2$ is defined in Definition \ref{pi1-def} (iv).
For $(i,j)\in \mathbb{Z}\times\mathbb{Z}$ and $s\in\mathbb{Z}_{\geq1}$, one defines
\[
L^2_{s,k,{\rm ad}}(i,j)=x_{s+P^k(i+k)+[i]_-+k-j,\pi_2(i+k)},\quad
L^2_{s,k,{\rm re}}(i,j)=x_{s+P^k(i+k-1)+[i-1]_-+k-j,\pi_2(i+k-1)}.
\]
Just as in (\ref{pos-A2-1}) and (\ref{pos-A2-2}), we see
that if $(i,j)$ is admissible in ${\rm REYD}_{{\rm D}^{(2)},k}$ then
\begin{equation}\label{pos-C1-1}
s+P^k(i+k)+[i]_-+k-j\geq s\geq 1
\end{equation}
by (\ref{pos-lem1}). If $(i,j)$ is removable then by (\ref{pos-lem2}),
\begin{equation}\label{pos-C1-2}
s+P^k(i+k-1)+[i-1]_-+k-j\geq s+1\geq 2.
\end{equation}

\begin{defn}\label{ad-rem-pt3}
Let $T=(y_t)_{t\in\mathbb{Z}}$ be a sequence in ${\rm REYD}_{{\rm D}^{(2)},k}$ of Definition \ref{DEYD}, $i\in\mathbb{Z}$
and $l\in\{0,n\}$.
\begin{enumerate}
\item We suppose that the point $(i,y_i)$ is admissible.
If $y_{i-1}<y_i=y_{i+1}$ and it holds either $i+k\equiv l+1$ and $i<0$ or $i+k\equiv l$ and $i>0$
then we say the point $(i,y_i)$ is a double $\pi_2(l)$-admissible point. 
\item We suppose that the point $(i,y_{i-1})$ is removable. If $y_{i-2}=y_{i-1}<y_{i}$ and it holds
either
$i+k-1\equiv l+1$ and $i>1$ or $i+k-1\equiv l$ and $i<1$
then we say the point $(i,y_{i-1})$ is a double $\pi_2(l)$-removable point. 
\item Other admissible (resp. removable) points $(i,y_i)$ (resp. $(i,y_{i-1})$) other than (i) (resp. (ii)) are said to be single $\pi_2(i+k)$-admissible (resp. $\pi_2(i+k-1)$-removable) points.
\end{enumerate}
Here, in (i) and (ii), the notation $a\equiv b$ means $a\equiv b$ (mod $2n$). 
\end{defn}
\nd
Note that if $l=0$ (resp. $l=n$) then $\pi_2(l)=1$ (resp. $\pi_2(l)=n$).
For $T\in{\rm REYD}_{{\rm D}^{(2)},k}$, we set
\begin{eqnarray}
L^2_{s,k,\iota}(T)&:=&
\sum_{t\in I} \left(
\sum_{P:\text{single }t\text{-admissible point of }T} L^2_{s,k,{\rm ad}}(P)
-\sum_{P:\text{single }t\text{-removable point of }T} L^2_{s,k,{\rm re}}(P)\right)\label{L2kdef}\\
& & + 
\sum_{P:\text{double }1\text{-admissible point of }T} 2L^2_{s,k,{\rm ad}}(P)
-\sum_{P:\text{double }1\text{-removable point of }T} 2L^2_{s,k,{\rm re}}(P) \nonumber
\\
& & + 
\sum_{P:\text{double }n\text{-admissible point of }T} 2L^2_{s,k,{\rm ad}}(P)
-\sum_{P:\text{double }n\text{-removable point of }T} 2L^2_{s,k,{\rm re}}(P)\in (\mathbb{Q}^{\infty})^*.\nonumber
\end{eqnarray}

\subsubsection{Assignment of inequalities to Young walls}

Let $k$ be $k=1$ or $k=n$ and we fix an integer $s\in\mathbb{Z}_{\geq 1}$.
We draw Young walls on $\mathbb{R}_{\leq0}\times \mathbb{R}_{\geq k}$. For example, in the case $n=3$, Young walls are drawn
as follows:
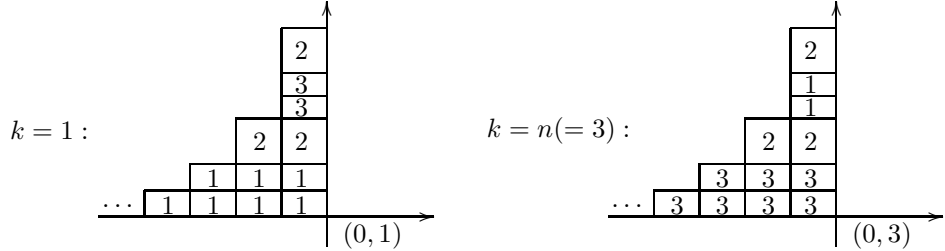
\begin{equation*}
\begin{xy}
(-30.5,8) *{k=1:}="000-1-a",
(-15.5,-2) *{\ 1}="000-1",
(-9.5,-2) *{\ 1}="00-1",
(-3.5,-2) *{\ 1}="0-1",
(1.5,-2) *{\ \ 1}="0-1",
(-21.5,-2) *{\dots}="00000",
(-9.5,1.5) *{\ 1}="000",
(-3.5,1.5) *{\ 1}="00",
(1.5,1.5) *{\ \ 1}="0",
(1.5,6.5) *{\ \ 2}="02",
(-4,6.5) *{\ \ 2}="002",
(1.5,11) *{\ \ 3}="03a",
(1.5,14) *{\ \ 3}="03b",
(1.5,18.5) *{\ \ 2}="04",
(12,-6) *{(0,1)}="origin",
(0,0) *{}="1",
(0,-3.5) *{}="1-u",
(6,0)*{}="2",
(6,-3.5)*{}="2-u",
(6,-7.5)*{}="-1-u",
(6,3.5)*{}="3",
(0,3.5)*{}="4",
(-6,3.5)*{}="5",
(-6,0)*{}="6",
(-6,-3.5)*{}="6-u",
(-12,3.5)*{}="7",
(-12,0)*{}="8",
(-12,-3.6)*{}="8-u",
(-18,3.5)*{}="9",
(-18,0)*{}="10-a",
(-18,0)*{}="10",
(-18,-3.5)*{}="10-u",
(-24,-3.5)*{}="11-u",
(20,-3.5)*{}="0-u",
(6,9.5)*{}="3-1",
(6,15.5)*{}="3-2",
(6,21.5)*{}="3-3",
(6,25)*{}="3-4",
(0,9.5)*{}="4-1",
(0,12.5)*{}="4-a",
(6,12.5)*{}="4-b",
(0,15.5)*{}="4-2",
(0,21.5)*{}="4-3",
(-6,9.5)*{}="5-1",
\ar@{-} "-1-u";"2-u"^{}
\ar@{->} "2-u";"0-u"^{}
\ar@{-} "10-u";"11-u"^{}
\ar@{-} "8";"10-a"^{}
\ar@{-} "10-u";"10-a"^{}
\ar@{-} "10-u";"2-u"^{}
\ar@{-} "8";"8-u"^{}
\ar@{-} "6";"6-u"^{}
\ar@{-} "2";"2-u"^{}
\ar@{-} "1";"1-u"^{}
\ar@{-} "1";"2"^{}
\ar@{-} "1";"4"^{}
\ar@{-} "2";"3"^{}
\ar@{-} "3";"4"^{}
\ar@{-} "5";"6"^{}
\ar@{-} "5";"4"^{}
\ar@{-} "1";"6"^{}
\ar@{-} "7";"8"^{}
\ar@{-} "7";"5"^{}
\ar@{-} "6";"8"^{}
\ar@{-} "3";"3-1"^{}
\ar@{-} "4-1";"3-1"^{}
\ar@{-} "4-1";"4"^{}
\ar@{-} "5-1";"5"^{}
\ar@{-} "4-1";"5-1"^{}
\ar@{-} "4-1";"4-2"^{}
\ar@{-} "3-2";"3-1"^{}
\ar@{-} "3-2";"4-2"^{}
\ar@{-} "3-2";"3-3"^{}
\ar@{-} "4-a";"4-b"^{}
\ar@{-} "4-3";"4-2"^{}
\ar@{-} "4-3";"3-3"^{}
\ar@{->} "3-3";"3-4"^{}
\end{xy}\qquad
\begin{xy}
(-30.5,8) *{k=n(=3):}="000-1-a",
(-15.5,-2) *{\ 3}="000-1",
(-9.5,-2) *{\ 3}="00-1",
(-3.5,-2) *{\ 3}="0-1",
(1.5,-2) *{\ \ 3}="0-1",
(-21.5,-2) *{\dots}="00000",
(-9.5,1.5) *{\ 3}="000",
(-3.5,1.5) *{\ 3}="00",
(1.5,1.5) *{\ \ 3}="0",
(1.5,6.5) *{\ \ 2}="02",
(-4,6.5) *{\ \ 2}="002",
(1.5,11) *{\ \ 1}="03a",
(1.5,14) *{\ \ 1}="03b",
(1.5,18.5) *{\ \ 2}="04",
(12,-6) *{(0,3)}="origin",
(0,0) *{}="1",
(0,-3.5) *{}="1-u",
(6,0)*{}="2",
(6,-3.5)*{}="2-u",
(6,-7.5)*{}="-1-u",
(6,3.5)*{}="3",
(0,3.5)*{}="4",
(-6,3.5)*{}="5",
(-6,0)*{}="6",
(-6,-3.5)*{}="6-u",
(-12,3.5)*{}="7",
(-12,0)*{}="8",
(-12,-3.6)*{}="8-u",
(-18,3.5)*{}="9",
(-18,0)*{}="10-a",
(-18,0)*{}="10",
(-18,-3.5)*{}="10-u",
(-24,-3.5)*{}="11-u",
(20,-3.5)*{}="0-u",
(6,9.5)*{}="3-1",
(6,15.5)*{}="3-2",
(6,21.5)*{}="3-3",
(6,25)*{}="3-4",
(0,9.5)*{}="4-1",
(0,12.5)*{}="4-a",
(6,12.5)*{}="4-b",
(0,15.5)*{}="4-2",
(0,21.5)*{}="4-3",
(-6,9.5)*{}="5-1",
\ar@{-} "-1-u";"2-u"^{}
\ar@{->} "2-u";"0-u"^{}
\ar@{-} "10-u";"11-u"^{}
\ar@{-} "8";"10-a"^{}
\ar@{-} "10-u";"10-a"^{}
\ar@{-} "10-u";"2-u"^{}
\ar@{-} "8";"8-u"^{}
\ar@{-} "6";"6-u"^{}
\ar@{-} "2";"2-u"^{}
\ar@{-} "1";"1-u"^{}
\ar@{-} "1";"2"^{}
\ar@{-} "1";"4"^{}
\ar@{-} "2";"3"^{}
\ar@{-} "3";"4"^{}
\ar@{-} "5";"6"^{}
\ar@{-} "5";"4"^{}
\ar@{-} "1";"6"^{}
\ar@{-} "7";"8"^{}
\ar@{-} "7";"5"^{}
\ar@{-} "6";"8"^{}
\ar@{-} "3";"3-1"^{}
\ar@{-} "4-1";"3-1"^{}
\ar@{-} "4-1";"4"^{}
\ar@{-} "5-1";"5"^{}
\ar@{-} "4-1";"5-1"^{}
\ar@{-} "4-1";"4-2"^{}
\ar@{-} "3-2";"3-1"^{}
\ar@{-} "3-2";"4-2"^{}
\ar@{-} "3-2";"3-3"^{}
\ar@{-} "4-a";"4-b"^{}
\ar@{-} "4-3";"4-2"^{}
\ar@{-} "4-3";"3-3"^{}
\ar@{->} "3-3";"3-4"^{}
\end{xy}
\end{equation*}
We inductively define integers $P^k(t)$ ($t\in\mathbb{Z}_{\geq k}$) as
\[
P^k(k):=0,\ \ P^k(t)=P^k(t-1)+p_{\pi'(t),\pi'(t-1)},
\]
where $\pi':\mathbb{Z}_{\geq1}\rightarrow\{1,2,\cdots,n\}$ was defined in (\ref{piprime}).
Let $i\in\mathbb{Z}_{\geq0}$, $l\in\mathbb{Z}_{\geq k}$ and $S$ be a slot or block
\[
S=
\begin{xy}
(-8,15) *{(-i-1,l+1)}="0000",
(17,15) *{(-i,l+1)}="000",
(15,-3) *{(-i,l)}="00",
(-7,-3) *{(-i-1,l)}="0",
(0,0) *{}="1",
(12,0)*{}="2",
(12,12)*{}="3",
(0,12)*{}="4",
\ar@{-} "1";"2"^{}
\ar@{-} "1";"4"^{}
\ar@{-} "2";"3"^{}
\ar@{-} "3";"4"^{}
\end{xy}
\]
in $\mathbb{R}_{\leq0}\times \mathbb{R}_{\geq k}$. 
If $S$ is colored by $t\in I\setminus \{1,n\}$ in the pattern of Definition \ref{def-YW} (ii)
then we set
\begin{equation}\label{l2def1}
L^2_{s,k,{\rm ad}}(S):=x_{s+P^k(l)+i,t},\quad L^2_{s,k,{\rm re}}(S):=x_{s+P^k(l)+i+1,t},
\end{equation}
which are similar assignments to (\ref{l1def1}).

Let $i\in\mathbb{Z}_{\geq0}$, $l\in\mathbb{Z}_{\geq k}$ and $S'$ be a slot or block colored by $t\in\{1,n\}$
in $\mathbb{R}_{\leq0}\times \mathbb{R}_{\geq k}$ such that the place is one of the following two:
\[
S'=
\begin{xy}
(-8,10) *{(-i-1,l+\frac{1}{2})}="0000",
(17,10) *{(-i,l+\frac{1}{2})}="000",
(15,-3) *{(-i,l)}="00",
(-7,-3) *{(-i-1,l)}="0",
(0,0) *{}="1",
(12,0)*{}="2",
(12,6)*{}="3",
(0,6)*{}="4",
(32,3) *{{\rm or}}="or",
(52,10) *{(-i-1,l+1)}="a0000",
(77,10) *{(-i,l+1)}="a000",
(75,-3) *{(-i,l+\frac{1}{2})}="a00",
(54,-3) *{(-i-1,l+\frac{1}{2})}="a0",
(60,0) *{}="a1",
(72,0)*{}="a2",
(72,6)*{}="a3",
(60,6)*{}="a4",
\ar@{-} "1";"2"^{}
\ar@{-} "1";"4"^{}
\ar@{-} "2";"3"^{}
\ar@{-} "3";"4"^{}
\ar@{-} "a1";"a2"^{}
\ar@{-} "a1";"a4"^{}
\ar@{-} "a2";"a3"^{}
\ar@{-} "a3";"a4"^{}
\end{xy}
\]
Then we set
\[
L^2_{s,k,{\rm ad}}(S'):=x_{s+P^k(l)+i,t},\quad L^2_{s,k,{\rm re}}(S'):=x_{s+P^k(l)+i+1,t}
\]
just as in (\ref{l1def2}).
It is easy to see 
\begin{equation}\label{C1YW-pr}
s+P^k(l)+i\geq s\geq1,\quad s+P^k(l)+i+1\geq s+1\geq2.
\end{equation}
In these cases, it holds $t=\pi'(l)$.
For a proper Young wall $Y$ of type ${\rm D}^{(2)}_{n}$ of ground state $\Lambda_k$, we define
\begin{eqnarray}
L^2_{s,k,\iota}(Y)&:=&
\sum_{t\in I} \left(
\sum_{P:\text{single }t\text{-admissible slot}} L^2_{s,k,{\rm ad}}(P)
-\sum_{P:\text{single removable }t\text{-block}} L^2_{s,k,{\rm re}}(P)\right)\nonumber\\
& & + 
\sum_{P:\text{double }1\text{-admissible slot}} 2L^2_{s,k,{\rm ad}}(P)
-\sum_{P:\text{double removable }1\text{-block}} 2L^2_{s,k,{\rm re}}(P)\label{L22-def}\\
& & + 
\sum_{P:\text{double }n\text{-admissible slot}} 2L^2_{s,k,{\rm ad}}(P)
-\sum_{P:\text{double removable }n\text{-block}} 2L^2_{s,k,{\rm re}}(P).\nonumber
\end{eqnarray}

Let ${\rm YW}_{{\rm D}^{(2)},k}$ be the set of all proper Young walls of type ${\rm D}^{(2)}_{n}$ of ground state $\Lambda_k$.

\subsubsection{Combinatorial description of ${\rm Im}(\Psi_{\iota})$ of type ${\rm C}_{n-1}^{(1)}$}

\begin{thm}\label{thmC1}
If $\mathfrak{g}$ is of type ${\rm C}_{n-1}^{(1)}$ $(n\geq3)$ and $\iota$ is adapted
then $\iota$ satisfies the $\Xi'$-positivity condition and
\[
{\rm Im}(\Psi_{\iota})
=\left\{ 
\textbf{a}\in\mathbb{Z}^{\infty} \left|
\begin{array}{l}
\text{for any }s\in\mathbb{Z}_{\geq1},\ k\in I\setminus\{1,n\},\ \\
\text{and}\ T\in{\rm REYD}_{{\rm D}^{(2)},k},\ \text{it holds}\ 
L^2_{s,k,\iota}(T)(\textbf{a})\geq0 \\
\text{and for any }k\in\{1,n\}\ \text{and }Y\in{\rm YW}_{{\rm D}^{(2)},k},\ 
\text{it holds}\ L^2_{s,k,\iota}(Y)(\textbf{a})\geq0
\end{array} \right.
\right\}.
\]

\end{thm}
In this way, inequalities of type ${\rm C}^{(1)}_{n-1}$ are expressed by 
combinatorial objects (revised extended Young diagrams, Young walls)
of type ${\rm D}^{(2)}_n$ as remarked in the end of subsection \ref{seno}.

\section{Action of $S'$}

\subsection{Type ${\rm A}_{n-1}^{(1)}$-case and ${\rm D}_{n}^{(2)}$-case}

In this subsection, we consider the following replacement of a concave corner in an extended Young diagram $T$
with a convex corner
and calculate how the values $\overline{L}_{s,k,\iota}(T)$, $L_{s,k,\iota}(T)$ (defined in (\ref{ovlL2}), (\ref{LL2})) are changed.
\begin{equation}\label{change}
\begin{xy}
(-6,13) *{(i,j)}="0000",
(4,-3) *{(i,j-1)}="000",
(10,15) *{(i+1,j)}="00",
(72,0) *{(i+1,j-1)}="a0000",
(54,-3) *{(i,j-1)}="a000",
(60,15) *{(i+1,j)}="a00",
(30,8) *{\rightarrow}="a00",
(0,0) *{\bullet}="1",
(12,12)*{\bullet}="3",
(0,12)*{\bullet}="4",
(50,0) *{\bullet}="a1",
(62,12)*{\bullet}="a3",
(62,0)*{\bullet}="a4",
\ar@{-} "a1";"a4"^{}
\ar@{-} "a3";"a4"^{}
\ar@{-} "1";"4"^{}
\ar@{-} "3";"4"^{}
\end{xy}
\end{equation}

\begin{prop}\label{prop-closednessAD}
We suppose that $T=(y_i)_{i\in\mathbb{Z}_{\geq0}}\in {\rm EYD}_{k}$ has an $(i+y_i)$-diagonal concave corner at a point $(i,y_i)$
and put $j:=y_i$.
Let $s\in\mathbb{Z}_{\geq1}$ and $T'\in{\rm EYD}_{k}$ be the extended Young diagram obtained from $T$
by replacing the $(i+j)$-diagonal concave corner by a convex corner. 
\begin{enumerate}
\item In the case $\mathfrak{g}$ is of type ${\rm A}_{n-1}^{(1)}$, it holds
\[
\overline{L}_{s,k,\iota}(T')=\overline{L}_{s,k,\iota}(T)-\beta_{s+P^k(i+j)+{\rm min}\{i,k-j\},\ovl{i+j}}.
\]
\item
In the case $\mathfrak{g}$ is of type ${\rm D}_{n}^{(2)}$, it holds
\[
L_{s,k,\iota}(T')=L_{s,k,\iota}(T)-\beta_{s+P^k(i+j)+{\rm min}\{i,k-j\},\pi(i+j)}.
\]
\end{enumerate}
\end{prop}

\nd
{\it Proof.} By the assumption $T$ has a concave corner at $(i,j)$, it holds $i>0$ and $y_{i-1}<y_{i}$ or $i=0$ so that
lines around the point $(i,j)$ in $T$ and $T'$ are as follows:
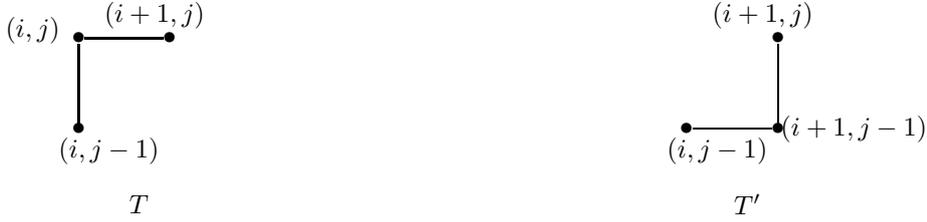
\[
\begin{xy}
(-6,13) *{(i,j)}="0000",
(4,-3) *{(i,j-1)}="000",
(10,15) *{(i+1,j)}="00",
(8,-10) *{T}="T",
(102,0) *{(i+1,j-1)}="a0000",
(84,-3) *{(i,j-1)}="a000",
(90,15) *{(i+1,j)}="a00",
(88,-10) *{T'}="T'",
(0,0) *{\bullet}="1",
(12,12)*{\bullet}="3",
(0,12)*{\bullet}="4",
(80,0) *{\bullet}="a1",
(92,12)*{\bullet}="a3",
(92,0)*{\bullet}="a4",
\ar@{-} "a1";"a4"^{}
\ar@{-} "a3";"a4"^{}
\ar@{-} "1";"4"^{}
\ar@{-} "3";"4"^{}
\end{xy}
\]
Our assumption means the point $(i,j)$ (resp. $(i+1,j-1)$) is a concave corner (resp. is not a corner) in $T$ and
is not a corner (resp. is a convex corner) in $T'$. Note that if $(i+1,j)$ is a convex corner in $T$, that is,
$y_{i}<y_{i+1}$ then $(i+1,j)$ is not a corner in $T'$. If $(i+1,j)$ is not a corner in $T$, that is,
$y_{i}=y_{i+1}$ then $(i+1,j)$ is a concave corner in $T'$.
Similarly,
If $(i,j-1)$ is a convex corner (resp. is not a corner) in $T$
then $(i,j-1)$ is not a corner (resp. is a concave corner) in $T'$. 
Other corners of $T$ are same as those of $T'$. 

(i) We suppose that $\mathfrak{g}$ is of type ${\rm A}_{n-1}^{(1)}$.
In the case $i>k-j$, it follows
\[
\overline{L}_{s,k,\iota}(i,j)=x_{s+P^k(i+j)+k-j,\overline{i+j}},\quad
\overline{L}_{s,k,\iota}(i+1,j)=x_{s+P^k(i+j+1)+k-j,\overline{i+j+1}},
\]
\[
\overline{L}_{s,k,\iota}(i,j-1)=x_{s+P^k(i+j-1)+k-j+1,\overline{i+j-1}},\quad
\overline{L}_{s,k,\iota}(i+1,j-1)=x_{s+P^k(i+j)+k-j+1,\overline{i+j}}.
\]
Hence, taking (\ref{pij2}), (\ref{pij3}), (\ref{A1pk1}) and the above argument into account,
\begin{eqnarray*}
\overline{L}_{s,k,\iota}(T')-\overline{L}_{s,k,\iota}(T)
&=&-x_{s+P^k(i+j)+k-j,\overline{i+j}}+x_{s+P^k(i+j+1)+k-j,\overline{i+j+1}}\\
& &+x_{s+P^k(i+j-1)+k-j+1,\overline{i+j-1}}-x_{s+P^k(i+j)+k-j+1,\overline{i+j}}\\
&=&
-x_{s+P^k(i+j)+k-j,\overline{i+j}}+x_{s+P^k(i+j)+p_{\ovl{i+j+1},\ovl{i+j}}+k-j,\overline{i+j+1}}\\
& &+x_{s+P^k(i+j)+p_{\ovl{i+j-1},\ovl{i+j}}+k-j,\overline{i+j-1}}-x_{s+P^k(i+j)+k-j+1,\overline{i+j}}\\
&=&-\beta_{s+P^k(i+j)+k-j,\ovl{i+j}}.
\end{eqnarray*}
Similarly, in the case $i=k-j$ so that $k=i+j$, it holds
\begin{eqnarray*}
\overline{L}_{s,k,\iota}(T')-\overline{L}_{s,k,\iota}(T)
&=&-x_{s+i,\overline{k}}+x_{s+P^k(k+1)+i,\overline{k+1}}
+x_{s+P^k(k-1)+i,\overline{k-1}}-x_{s+i+1,\overline{k}}\\
&=&
-x_{s+i,\overline{k}}+x_{s+p_{\ovl{k+1},\ovl{k}}+i,\overline{k+1}}
+x_{s+p_{\ovl{k-1},\ovl{k}}+i,\overline{k-1}}-x_{s+i+1,\overline{k}}\\
&=&-\beta_{s+i,\ovl{k}}.
\end{eqnarray*}
In the case $i<k-j$, it holds
\[
\overline{L}_{s,k,\iota}(i,j)=x_{s+P^k(i+j)+i,\overline{i+j}},\quad
\overline{L}_{s,k,\iota}(i+1,j)=x_{s+P^k(i+j+1)+i+1,\overline{i+j+1}},
\]
\[
\overline{L}_{s,k,\iota}(i,j-1)=x_{s+P^k(i+j-1)+i,\overline{i+j-1}},\quad
\overline{L}_{s,k,\iota}(i+1,j-1)=x_{s+P^k(i+j)+i+1,\overline{i+j}}
\]
and by (\ref{pij2}), (\ref{pij3}) and (\ref{A1pk2}),
\begin{eqnarray*}
\overline{L}_{s,k,\iota}(T')-\overline{L}_{s,k,\iota}(T)
&=&-x_{s+P^k(i+j)+i,\overline{i+j}}+x_{s+P^k(i+j+1)+i+1,\overline{i+j+1}}\\
& &+x_{s+P^k(i+j-1)+i,\overline{i+j-1}}-x_{s+P^k(i+j)+i+1,\overline{i+j}}\\
&=&
-x_{s+P^k(i+j)+i,\overline{i+j}}+x_{s+P^k(i+j)+p_{\ovl{i+j+1},\ovl{i+j}}+i,\overline{i+j+1}}\\
& &+x_{s+P^k(i+j)+p_{\ovl{i+j-1},\ovl{i+j}}+i,\overline{i+j-1}}-x_{s+P^k(i+j)+i+1,\overline{i+j}}\\
&=&-\beta_{s+P^k(i+j)+i,\ovl{i+j}}.
\end{eqnarray*}
Therefore, we get $\overline{L}_{s,k,\iota}(T')=\overline{L}_{s,k,\iota}(T)-\beta_{s+P^k(i+j)+{\rm min}\{i,k-j\},\ovl{i+j}}$.

(ii) Next, we assume $\mathfrak{g}$ is of type ${\rm D}_{n}^{(2)}$.
In the case $i>k-j$, the definition of $L_{s,k,\iota}$ means
\[
L_{s,k,\iota}(i,j)=x_{s+P^k(i+j)+k-j,\pi(i+j)},\quad
L_{s,k,\iota}(i+1,j)=x_{s+P^k(i+j+1)+k-j,\pi(i+j+1)},
\]
\[
L_{s,k,\iota}(i,j-1)=x_{s+P^k(i+j-1)+k-j+1,\pi(i+j-1)},\quad
L_{s,k,\iota}(i+1,j-1)=x_{s+P^k(i+j)+k-j+1,\pi(i+j)}.
\]
Note that if $\pi(i+j)=1$ then $\pi(i+j-1)=\pi(i+j+1)=2$ and $P^k(i+j-1)+1=P^k(i+j-1)+p_{1,2}+p_{2,1}=P^k(i+j)+p_{2,1}=P^k(i+j+1)$ so that
\[
L_{s,k,\iota}(i+1,j)=L_{s,k,\iota}(i,j-1)=x_{s+P^k(i+j)+k-j+p_{2,1},2}.
\]
Similarly, if $\pi(i+j)=n$ then
\[
L_{s,k,\iota}(i+1,j)=L_{s,k,\iota}(i,j-1)=x_{s+P^k(i+j)+k-j+p_{n-1,n},n-1}.
\]
Thus, a similar argument to (i) yields
$L_{s,k,\iota}(T')-L_{s,k,\iota}(T)=-\beta_{s+P^k(i+j)+k-j,\pi(i+j)}$.
Considering the cases $i=k-j$ and $i<k-j$, it follows 
\[
L_{s,k,\iota}(T')-L_{s,k,\iota}(T)=-\beta_{s+P^k(i+j)+{\rm min}\{i,k-j\},\pi(i+j)}.
\]
\qed

\subsection{Type ${\rm A}_{2n-2}^{(2)}$-case (action on extended Young diagrams)}

In this subsection, we assume $\mathfrak{g}$ is of type ${\rm A}_{2n-2}^{(2)}$ and take $k\in I\setminus\{1\}$ and prove the following proposition.
Recall that we defined the notion of admissible and removable points in Definition \ref{ad-rem-pt} and the notation $L^1_{s,k,\iota}$ in (\ref{L1kdef}).

\begin{prop}\label{A2closed}
Let $T=(y_t)_{t\in\mathbb{Z}}$ be a sequence in ${\rm REYD}_{{\rm A}^{(2)},k}$ of Definition \ref{AEYD} and $i\in\mathbb{Z}$.
\begin{enumerate}
\item We suppose that the point $(i,y_i)$ is single or double admissible and let $T'=(y_t')_{t\in\mathbb{Z}}$ be the sequence in ${\rm REYD}_{{\rm A}^{(2)},k}$
such that $y_{i}'=y_i-1$ and $y_t'=y_t$ $(t\neq i)$. Then for $s\in \mathbb{Z}_{\geq1}$, putting $j:=y_i$, it holds
\[
L^1_{s,k,\iota}(T')-L^1_{s,k,\iota}(T)=-\beta_{s+P^k(i+k)+[i]_-+k-j,\pi_1(i+k)}.
\]
\item We suppose that the point $(i,y_{i-1})$ is single or double removable and let $T''=(y_t'')_{t\in\mathbb{Z}}$ be the sequence in ${\rm REYD}_{{\rm A}^{(2)},k}$
such that $y_{i-1}''=y_{i-1}+1$ and $y_t''=y_t$ $(t\neq i-1)$. Then for $s\in \mathbb{Z}_{\geq1}$, putting $j:=y_{i-1}''=y_{i-1}+1$, it holds
\[
L^1_{s,k,\iota}(T'')-L^1_{s,k,\iota}(T)=\beta_{s+P^k(i+k-1)+[i-1]_-+k-j,\pi_1(i+k-1)}.
\]
\end{enumerate}
\end{prop}

\nd
In this proposition, (ii) follows from (i) since $(i-1,y_{i-1}'')$ is admissible in $T''$ and $T=(y_t)_{t\in\mathbb{Z}}$ is the sequence in ${\rm REYD}_{{\rm A}^{(2)},k}$
such that $y_{i-1}=y_{i-1}''-1$ and $y_t=y_t''$ $(t\neq i-1)$. Thus, we prove (i) by dividing our claim into four lemmas.
In the proofs of lemmas, if we say $m$-admissible/removable then it means $m$-single admissible/removable for $m\in I$.

\begin{lem}\label{lemA2-1}
In the case $\pi_1(i+k)>2$, Proposition \ref{A2closed} (i) holds.
\end{lem}

\nd
{\it Proof.}

It follows from $\pi_1(i+k)>2$
that $i+k\not\equiv0$ and $i-1+k\not\equiv0$ (mod $2n-1$).
By Definition \ref{AEYD} (3), it holds $y_{i-1}\in\{y_i,y_i-1\}$ and $y_i\in\{y_{i+1},y_{i+1}-1\}$.
The assumption of Proposition \ref{A2closed} (i) means $y_{i-1}=y_{i}-1$ and $y_{i}=y_{i+1}$ (if not, $T'\notin{\rm REYD}_{{\rm A}^{(2)},k}$)
 so that lines around the point $(i,j)$ in $T$ are as follows:
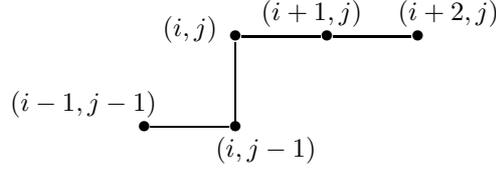
\[
\begin{xy}
(-20,3) *{(i-1,j-1)}="00000",
(-6,13) *{(i,j)}="0000",
(4,-3) *{(i,j-1)}="000",
(10,15) *{(i+1,j)}="00",
(28,15) *{(i+2,j)}="0",
(0,0) *{\bullet}="1",
(-12,0) *{\bullet}="2",
(12,12)*{\bullet}="3",
(0,12)*{\bullet}="4",
(24,12)*{\bullet}="5",
\ar@{-} "1";"4"^{}
\ar@{-} "1";"2"^{}
\ar@{-} "3";"4"^{}
\ar@{-} "3";"5"^{}
\end{xy}
\]
Since $y_{i-1}'=y_i'$ and $y_{i+1}'=y_i'+1$, lines
around the point $(i,j)$ in $T'$ are as follows:
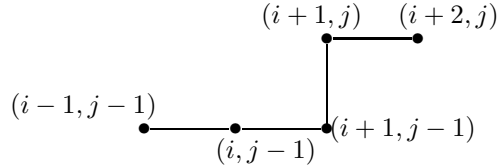
\[
\begin{xy}
(-20,3) *{(i-1,j-1)}="00000",
(22,0) *{(i+1,j-1)}="0000",
(4,-3) *{(i,j-1)}="000",
(10,15) *{(i+1,j)}="00",
(28,15) *{(i+2,j)}="0",
(0,0) *{\bullet}="1",
(-12,0) *{\bullet}="2",
(12,12)*{\bullet}="3",
(12,0)*{\bullet}="4",
(24,12)*{\bullet}="5",
\ar@{-} "1";"4"^{}
\ar@{-} "1";"2"^{}
\ar@{-} "3";"4"^{}
\ar@{-} "3";"5"^{}
\end{xy}
\]
The point $(i,j)$ (resp. $(i+1,j-1)$) is a single $\pi_1(i+k)$-admissible (resp. $\pi_1(i+k)$-removable) point in $T$ (resp. in $T'$).
By $i+1+k\not\equiv0$, $i-2+k\not\equiv0$ (mod $2n-1$) and
the rule in Definition \ref{AEYD} (3), the point $(i+2,j)$ (resp. $(i-1,j-1)$) is either a single removable (resp. admissible) corner or not a corner in $T$ (resp. in $T'$). 
Note that the point $(i+2,j)$ is a single removable corner in $T$ if and only if the point $(i+1,j)$ is not an admissible corner in $T'$.
Similarly, the point $(i-1,j-1)$ is a single admissible corner in $T'$ if and only if the point $(i,j-1)$ is not a single removable corner in $T$.
One can summarize them as follows:

\begin{table}[htb]
  \begin{tabular}{|c|c|c|} \hline
    points & in $T$ & in $T'$ \\ \hline
    $(i,j)$ & $\pi_1(i+k)$-admissible & normal \\
    $(i+1,j-1)$ & normal & $\pi_1(i+k)$-removable \\
    $(i+2,j)$ & $\pi_1(i+k+1)$-removable (resp. normal) & normal \\
    $(i+1,j)$ & normal & normal (resp. $\pi_1(i+k+1)$-admissible) \\
    \hline
  \end{tabular}
\end{table}

\begin{table}[H]
  \begin{tabular}{|c|c|c|} \hline
    points & in $T$ & in $T'$ \\ \hline
    $(i-1,j-1)$ & normal & $\pi_1(i+k-1)$-admissible (resp. normal) \\
    $(i,j-1)$ & normal (resp. $\pi_1(i+k-1)$-removable corner) & normal \\
    \hline
  \end{tabular}
\end{table}
\nd
Here, `normal' means the point is neither admissible nor removable.
Other points in $T$ are same as in $T'$. Taking (\ref{L1kdef}) into account, we obtain
\begin{eqnarray*}
L^1_{s,k,\iota}(T')-L^1_{s,k,\iota}(T)&=&
-x_{s+P^k(i+k)+[i]_-+k-j,\pi_1(i+k)}
-x_{s+P^k(i+k)+[i]_-+k-j+1,\pi_1(i+k)}\\
& &+x_{s+P^k(i+k+1)+[i+1]_-+k-j,\pi_1(i+k+1)}
+x_{s+P^k(i+k-1)+[i-1]_-+k-j+1,\pi_1(i+k-1)}.
\end{eqnarray*}
If $i>0$ then we get
\[
P^k(i+k+1)=P^k(i+k)+p_{\pi_1(i+k+1),\pi_1(i+k)},
\]
\[
P^k(i+k-1)=P^k(i+k)-p_{\pi_1(i+k),\pi_1(i+k-1)}
=P^k(i+k)+p_{\pi_1(i+k-1),\pi_1(i+k)}-1.
\]
Thus,
\begin{eqnarray*}
L^1_{s,k,\iota}(T')-L^1_{s,k,\iota}(T)&=&
-x_{s+P^k(i+k)+k-j,\pi_1(i+k)}
-x_{s+P^k(i+k)+k-j+1,\pi_1(i+k)}\\
& &+x_{s+P^k(i+k)+p_{\pi_1(i+k+1),\pi_1(i+k)}+k-j,\pi_1(i+k+1)}
+x_{s+P^k(i+k)+p_{\pi_1(i+k-1),\pi_1(i+k)}+k-j,\pi_1(i+k-1)}\\
&=&-\beta_{s+P^k(i+k)+k-j,\pi_1(i+k)}.
\end{eqnarray*}
If $i=0$ then $P^k(i+k)=P^k(k)=0$ and
\begin{eqnarray*}
L^1_{s,k,\iota}(T')-L^1_{s,k,\iota}(T)&=&
-x_{s+k-j,\pi_1(k)}
-x_{s+k-j+1,\pi_1(k)}\\
& &+x_{s+p_{\pi_1(k+1),\pi_1(k)}+k-j,\pi_1(k+1)}
+x_{s+p_{\pi_1(k-1),\pi_1(k)}+k-j,\pi_1(k-1)}\\
&=&
-\beta_{s+k-j,\pi_1(k)}.
\end{eqnarray*}
If $i<0$ then one gets $[i]_-=i$, $[i+1]_-=i+1$, $[i-1]_-=i-1$,
\[
P^k(i+k+1)=P^k(i+k)-p_{\pi_1(i+k),\pi_1(i+k+1)}
=P^k(i+k)-1+p_{\pi_1(i+k+1),\pi_1(i+k)},
\]
\[
P^k(i+k-1)=P^k(i+k)+p_{\pi_1(i+k-1),\pi_1(i+k)},
\]
and
\begin{eqnarray*}
L^1_{s,k,\iota}(T')-L^1_{s,k,\iota}(T)&=&
-x_{s+P^k(i+k)+i+k-j,\pi_1(i+k)}
-x_{s+P^k(i+k)+i+k-j+1,\pi_1(i+k)}\\
& &+x_{s+P^k(i+k+1)+i+1+k-j,\pi_1(i+k+1)}
+x_{s+P^k(i+k-1)+i-1+k-j+1,\pi_1(i+k-1)}\\
&=&
-x_{s+P^k(i+k)+i+k-j,\pi_1(i+k)}
-x_{s+P^k(i+k)+i+k-j+1,\pi_1(i+k)}\\
& &+x_{s+P^k(i+k)+p_{\pi_1(i+k+1),\pi_1(i+k)}+i+k-j,\pi_1(i+k+1)}\\
& &
+x_{s+P^k(i+k)+p_{\pi_1(i+k-1),\pi_1(i+k)}+i+k-j,\pi_1(i+k-1)}\\
&=&-\beta_{s+P^k(i+k)+i+k-j,\pi_1(i+k)}.
\end{eqnarray*}
Therefore, our claim $L^1_{s,k,\iota}(T')-L^1_{s,k,\iota}(T)=-\beta_{s+P^k(i+k)+[i]_-+k-j,\pi_1(i+k)}$ follows. 
\qed

\begin{lem}\label{lemA2-2}
In the case $\pi_1(i+k)=2$, Proposition \ref{A2closed} (i) holds.
\end{lem}

\nd
{\it Proof.} By a similar argument to the proof of Lemma \ref{lemA2-1}, putting $j:=y_i$,
we get
$y_{i-1}=y_{i}-1$ and $y_{i}=y_{i+1}$ so that
lines around $(i,j)$ in $T$ are as follows:
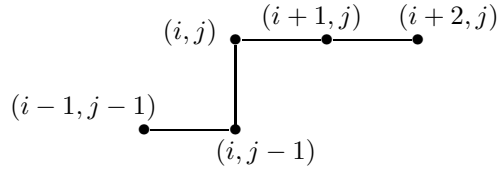
\[
\begin{xy}
(-20,3) *{(i-1,j-1)}="00000",
(-6,13) *{(i,j)}="0000",
(4,-3) *{(i,j-1)}="000",
(10,15) *{(i+1,j)}="00",
(28,15) *{(i+2,j)}="0",
(0,0) *{\bullet}="1",
(-12,0) *{\bullet}="2",
(12,12)*{\bullet}="3",
(0,12)*{\bullet}="4",
(24,12)*{\bullet}="5",
\ar@{-} "1";"4"^{}
\ar@{-} "1";"2"^{}
\ar@{-} "3";"4"^{}
\ar@{-} "3";"5"^{}
\end{xy}
\]
Lines
around the point $(i,j)$ in $T'$ are as follows:
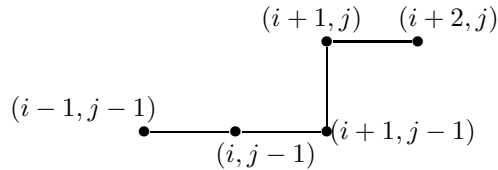
\[
\begin{xy}
(-20,3) *{(i-1,j-1)}="00000",
(22,0) *{(i+1,j-1)}="0000",
(4,-3) *{(i,j-1)}="000",
(10,15) *{(i+1,j)}="00",
(28,15) *{(i+2,j)}="0",
(0,0) *{\bullet}="1",
(-12,0) *{\bullet}="2",
(12,12)*{\bullet}="3",
(12,0)*{\bullet}="4",
(24,12)*{\bullet}="5",
\ar@{-} "1";"4"^{}
\ar@{-} "1";"2"^{}
\ar@{-} "3";"4"^{}
\ar@{-} "3";"5"^{}
\end{xy}
\]
First, we assume $\pi_1(i+k+1)=1$ and $\pi_1(i+k-1)=3$, which means 
$i+k+1\equiv 0$ (${\rm mod}\ 2n-1$) and $\pi_1(i+k+2)=1$.
Since $k+i+1\equiv 0$ (${\rm mod}\ 2n-1$), it holds either $y_{i+2}=y_{i+1}$ or $y_{i+2}>y_{i+1}$ or $y_{i+2}<y_{i+1}$.
Now, we supposed $1<k\leq n$ so that $k \not\equiv 0,1$ (${\rm mod}\ 2n-1$) and $i\neq -1,-2$.

\vspace{2mm}

\nd
\underline{Case 1. $y_{i+2}=y_{i+1}$ and $y_{i+2}<y_{i+3}$}

\vspace{2mm}

\nd
In this case, the lines around $(i,j)$ in $T$ are as follows:
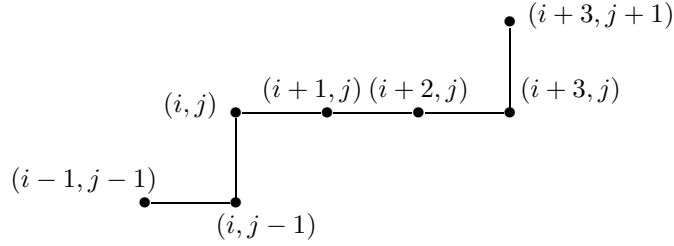
\[
\begin{xy}
(-20,3) *{(i-1,j-1)}="00000",
(-6,13) *{(i,j)}="0000",
(4,-3) *{(i,j-1)}="000",
(10,15) *{(i+1,j)}="00",
(24,15) *{(i+2,j)}="0",
(44,15) *{(i+3,j)}="0-a",
(48,25) *{(i+3,j+1)}="0-ab",
(0,0) *{\bullet}="1",
(-12,0) *{\bullet}="2",
(12,12)*{\bullet}="3",
(0,12)*{\bullet}="4",
(24,12)*{\bullet}="5",
(36,12)*{\bullet}="6",
(36,24)*{\bullet}="7",
\ar@{-} "1";"4"^{}
\ar@{-} "1";"2"^{}
\ar@{-} "3";"4"^{}
\ar@{-} "3";"5"^{}
\ar@{-} "6";"5"^{}
\ar@{-} "6";"7"^{}
\end{xy}
\]
Thus, the table of admissible and removable points are as follows:
\begin{table}[H]
  \begin{tabular}{|c|c|c|} \hline
    points & in $T$ & in $T'$ \\ \hline
    $(i,j)$ & $2$-admissible & normal  \\ \hline
    $(i+1,j-1)$ & normal & 2-removable \\ \hline
    $(i+2,j)$ & $1$-removable if $i+1<0$, & normal \\ 
    & normal if $i+1>0$ & \\ \hline
    $(i+1,j)$ & normal & $1$-double admissible if $i+1>0$, \\ 
    & & $1$-single admissible if $i+1<0$ \\\hline
    $(i-1,j-1)$ & normal & normal (resp. $3$-admissible) \\ \hline
    $(i,j-1)$ & $3$-removable (resp. normal) & normal \\
    \hline
  \end{tabular}
\end{table}

\nd
Since other points in $T$ are same as in $T'$, it follows by (\ref{L1kdef}) that
\begin{eqnarray*}
L^1_{s,k,\iota}(T')-L^1_{s,k,\iota}(T)&=&
-x_{s+P^k(i+k)+[i]_-+k-j,2}
-x_{s+P^k(i+k)+[i]_-+k-j+1,2}\\
& &+2x_{s+P^k(i+k+1)+[i+1]_-+k-j,1}
+x_{s+P^k(i+k-1)+[i-1]_-+k-j+1,3}.
\end{eqnarray*}
By the same argument as in the proof of Lemma \ref{lemA2-1}, one gets
\[
L^1_{s,k,\iota}(T')-L^1_{s,k,\iota}(T)=-\beta_{s+P^k(i+k)+[i]_-+k-j,2}.
\]

\vspace{2mm}

\nd
\underline{Case 2. $y_{i+2}=y_{i+1}$ and $y_{i+2}=y_{i+3}$}

\vspace{2mm}

\nd
In this case, the lines around $(i,j)$ in $T$ are as follows:
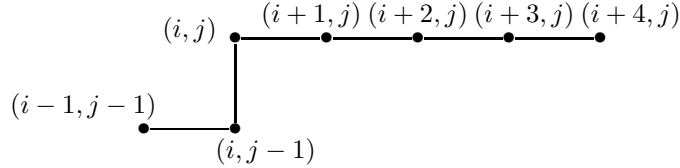
\[
\begin{xy}
(-20,3) *{(i-1,j-1)}="00000",
(-6,13) *{(i,j)}="0000",
(4,-3) *{(i,j-1)}="000",
(10,15) *{(i+1,j)}="00",
(24,15) *{(i+2,j)}="0",
(38,15) *{(i+3,j)}="0-a",
(52,15) *{(i+4,j)}="0-ab",
(0,0) *{\bullet}="1",
(-12,0) *{\bullet}="2",
(12,12)*{\bullet}="3",
(0,12)*{\bullet}="4",
(24,12)*{\bullet}="5",
(36,12)*{\bullet}="6",
(48,12)*{\bullet}="7",
\ar@{-} "1";"4"^{}
\ar@{-} "1";"2"^{}
\ar@{-} "3";"4"^{}
\ar@{-} "3";"5"^{}
\ar@{-} "6";"5"^{}
\ar@{-} "6";"7"^{}
\end{xy}
\]
The table is same as in Case 1.
except for the point $(i+2,j)$:

\begin{table}[H]
  \begin{tabular}{|c|c|c|} \hline
    points & in $T$ & in $T'$ \\ \hline
    $(i,j)$ & $2$-admissible & normal  \\ \hline
    $(i+1,j-1)$ & normal & 2-removable \\ \hline
     & $1$-removable and $1$-admissible & $1$-admissible if $i+1<0$ \\ 
    $(i+2,j)$ & if $i+1<0$, & \\ 
    & normal if $i+1>0$ & normal if  $i+1>0$ \\ \hline
    $(i+1,j)$ & normal & $1$-double admissible if $i+1>0$, \\ 
    & & $1$-single admissible if $i+1<0$ \\\hline
    $(i-1,j-1)$ & normal & normal (resp. $3$-admissible) \\ \hline
    $(i,j-1)$ & $3$-removable (resp. normal) & normal \\
    \hline
  \end{tabular}
\end{table}
\nd
Thus, just as in Case 1.,
one gets
\[
L^1_{s,k,\iota}(T')-L^1_{s,k,\iota}(T)=-\beta_{s+P^k(i+k)+[i]_-+k-j,2}.
\]

\vspace{2mm}

\nd
\underline{Case 3. $y_{i+2}>y_{i+1}$}

\vspace{2mm}

\nd
In this case, the lines around $(i,j)$ in $T$ are as follows:
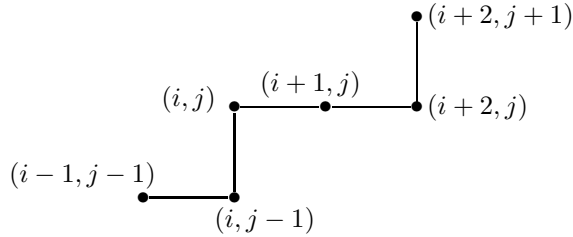
\[
\begin{xy}
(-20,3) *{(i-1,j-1)}="00000",
(-6,13) *{(i,j)}="0000",
(4,-3) *{(i,j-1)}="000",
(10,15) *{(i+1,j)}="00",
(32,12) *{(i+2,j)}="0",
(35,24) *{(i+2,j+1)}="0-a",
(0,0) *{\bullet}="1",
(-12,0) *{\bullet}="2",
(12,12)*{\bullet}="3",
(0,12)*{\bullet}="4",
(24,12)*{\bullet}="5",
(24,24)*{\bullet}="6",
\ar@{-} "1";"4"^{}
\ar@{-} "1";"2"^{}
\ar@{-} "3";"4"^{}
\ar@{-} "3";"5"^{}
\ar@{-} "6";"5"^{}
\end{xy}
\]
The table of admissible and removable points is as follows:

\begin{table}[H]
  \begin{tabular}{|c|c|c|} \hline
    points & in $T$ & in $T'$ \\ \hline
    $(i,j)$ & $2$-admissible & normal  \\ \hline
    $(i+1,j-1)$ & normal & 2-removable \\ \hline
    $(i+2,j)$ & $1$-double removable if $i+2<1$, & normal  \\ 
    & $1$-single removable if $i+2>1$ &  \\ \hline
    $(i+1,j)$ & normal & $1$-single admissible if $i+1>0,$ \\ 
    & & normal if $i+1<0$ \\\hline
    $(i-1,j-1)$ & normal & normal (resp. $3$-admissible) \\ \hline
    $(i,j-1)$ & $3$-removable (resp. normal) & normal \\
    \hline
  \end{tabular}
\end{table}
\nd
Hence, one can verify
\begin{eqnarray*}
L^1_{s,k,\iota}(T')-L^1_{s,k,\iota}(T)&=&
-x_{s+P^k(i+k)+[i]_-+k-j,2}
-x_{s+P^k(i+k)+[i]_-+k-j+1,2}\\
& &+2x_{s+P^k(i+k+1)+[i+1]_-+k-j,1}
+x_{s+P^k(i+k-1)+[i-1]_-+k-j+1,3}\\
&=&-\beta_{s+P^k(i+k)+[i]_-+k-j,2}.
\end{eqnarray*}

\vspace{2mm}

\nd
\underline{Case 4. $y_{i+1}>y_{i+2}$}

\vspace{2mm}

\nd
In this case, the lines around $(i,j)$ in $T$ are as follows:
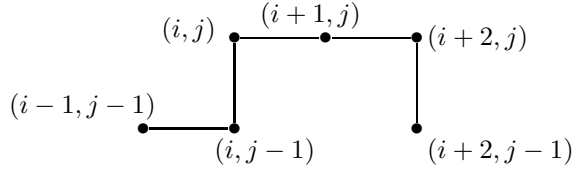
\[
\begin{xy}
(-20,3) *{(i-1,j-1)}="00000",
(-6,13) *{(i,j)}="0000",
(4,-3) *{(i,j-1)}="000",
(10,15) *{(i+1,j)}="00",
(32,12) *{(i+2,j)}="0",
(35,-3) *{(i+2,j-1)}="0-a",
(0,0) *{\bullet}="1",
(-12,0) *{\bullet}="2",
(12,12)*{\bullet}="3",
(0,12)*{\bullet}="4",
(24,12)*{\bullet}="5",
(24,0)*{\bullet}="6",
\ar@{-} "1";"4"^{}
\ar@{-} "1";"2"^{}
\ar@{-} "3";"4"^{}
\ar@{-} "3";"5"^{}
\ar@{-} "6";"5"^{}
\end{xy}
\]
Only if $i+1<0$, this Case 4. happens.  
The table is as follows:

\begin{table}[H]
  \begin{tabular}{|c|c|c|} \hline
    points & in $T$ & in $T'$ \\ \hline
    $(i,j)$ & $2$-admissible & normal  \\ \hline
    $(i+1,j-1)$ & normal & $2$-removable \\ \hline
    $(i+2,j)$ & $1$-single removable & normal  \\ \hline
    $(i+1,j)$ & normal & $1$-single admissible \\ \hline
    $(i-1,j-1)$ & normal & normal (resp. $3$-admissible) \\ \hline
    $(i,j-1)$ & $3$-removable (resp. normal) & normal \\
    \hline
  \end{tabular}
\end{table}
\nd
Hence, one can verify
\begin{eqnarray*}
L^1_{s,k,\iota}(T')-L^1_{s,k,\iota}(T)&=&
-x_{s+P^k(i+k)+[i]_-+k-j,2}
-x_{s+P^k(i+k)+[i]_-+k-j+1,2}\\
& &+2x_{s+P^k(i+k+1)+[i+1]_-+k-j,1}
+x_{s+P^k(i+k-1)+[i-1]_-+k-j+1,3}\\
&=&-\beta_{s+P^k(i+k)+[i]_-+k-j,2}.
\end{eqnarray*}

Next, we assume $\pi_1(i+k+1)=3$ and $\pi_1(i+k-1)=1$, which means 
$i+k-1\equiv 1$, $i+k-2\equiv 0$ (${\rm mod}\ 2n-1$). By $1<k\leq n$, it holds $i\neq 1,2$.
Since one can prove our claim by a similar argument to the case $\pi_1(i+k+1)=1$ and $\pi_1(i+k-1)=3$,
we write only the cases one should consider and the table of admissible and removable points.

\vspace{2mm}

\nd
\underline{Case 1. $y_{i-2}=y_{i-1}$ and $y_{i-2}>y_{i-3}$}

\vspace{2mm}

\begin{table}[H]
  \begin{tabular}{|c|c|c|} \hline
    points & in $T$ & in $T'$ \\ \hline
    $(i,j)$ & $2$-admissible & normal  \\ \hline
    $(i+1,j-1)$ & normal & $2$-removable  \\ \hline
    $(i-1,j-1)$ & normal & $1$-single admissible if $i-2<0$, \\
    & & normal if $i-2>0$
    \\ \hline
    $(i,j-1)$ & $1$-double removable if $i>1$, & normal \\ 
    & $1$-single removable if $i<1$ & \\ \hline
    $(i+2,j)$ & $3$-removable (resp. normal) & normal \\ \hline
    $(i+1,j)$ & normal & normal (resp. $3$-admissible) \\
    \hline
  \end{tabular}
\end{table}
\nd
We remark that $i>1$ (resp. $i<1$) if and only if $i>2$ (resp. $i<2$) by $i\neq 1,2$.

\vspace{2mm}

\nd
\underline{Case 2. $y_{i-2}=y_{i-1}$ and $y_{i-2}=y_{i-3}$}

\vspace{2mm}

\begin{table}[H]
  \begin{tabular}{|c|c|c|} \hline
    points & in $T$ & in $T'$ \\ \hline
    $(i,j)$ & $2$-admissible & normal  \\ \hline
    $(i+1,j-1)$ & normal & $2$-removable  \\ \hline
    $(i-1,j-1)$ & $1$-removable if $i-2<0$, & $1$-removable and $1$-admissible if $i-2<0$, \\
    & normal if $i-2>0$ & normal if $i-2>0$ \\ \hline
    $(i,j-1)$ & $1$-double removable if $i>1$, & normal \\
     & $1$-single removable if $i<1$ & \\ \hline
    $(i+2,j)$ & $3$-removable (resp. normal) & normal \\ \hline
    $(i+1,j)$ & normal & normal (resp. $3$-admissible) \\
    \hline
  \end{tabular}
\end{table}

\vspace{2mm}

\nd
\underline{Case 3. $y_{i-2}<y_{i-1}$}

\vspace{2mm}

\begin{table}[H]
  \begin{tabular}{|c|c|c|} \hline
    points & in $T$ & in $T'$ \\ \hline
    $(i,j)$ & $2$-admissible & normal  \\ \hline
    $(i+1,j-1)$ & normal & $2$-removable  \\ \hline
    $(i-1,j-1)$ & normal & $1$-double admissible if $i-1<0$, \\
    & & $1$-single admissible if $i-1>0$
    \\ \hline
    $(i,j-1)$ & $1$-removable if $i-2>0$, & normal \\ 
    & normal if $i-2<0$ & \\ \hline
    $(i+2,j)$ & $3$-removable (resp. normal) & normal \\ \hline
    $(i+1,j)$ & normal & normal (resp. $3$-admissible) \\
    \hline
  \end{tabular}
\end{table}

\vspace{2mm}

\nd
\underline{Case 4. $y_{i-2}>y_{i-1}$}

\vspace{2mm}

\nd
In this case, it holds $i-2<0$.

\begin{table}[H]
  \begin{tabular}{|c|c|c|} \hline
    points & in $T$ & in $T'$ \\ \hline
    $(i,j)$ & $2$-admissible & normal  \\ \hline
    $(i+1,j-1)$ & normal & $2$-removable  \\ \hline
    $(i-1,j-1)$ & normal & $1$-admissible \\ \hline
    $(i,j-1)$ & $1$-removable & normal \\ \hline
    $(i+2,j)$ & $3$-removable (resp. normal) & normal \\ \hline
    $(i+1,j)$ & normal & normal (resp. $3$-admissible) \\
    \hline
  \end{tabular}
\end{table}
\qed

\begin{lem}\label{lemA2-3}
In the case $i+k\equiv 1$ $({\rm mod}\ 2n-1)$ so that $\pi_1(i+k)=1$, Proposition \ref{A2closed} (i) holds.
\end{lem}

\nd
{\it Proof.} Since $i-1+k\equiv 0$ $({\rm mod}\ 2n-1)$, there are three patterns (1) $y_{i-1}<y_{i}$, (2) $y_{i-1}=y_{i}$, (3) $y_{i-1}>y_{i}$ (Definition \ref{AEYD}).
We get $y_{i}=y_{i+1}$.
Putting $j:=y_i$,
lines around the point $(i,j)$ in $T$ and $T'$ of each pattern are as follows:
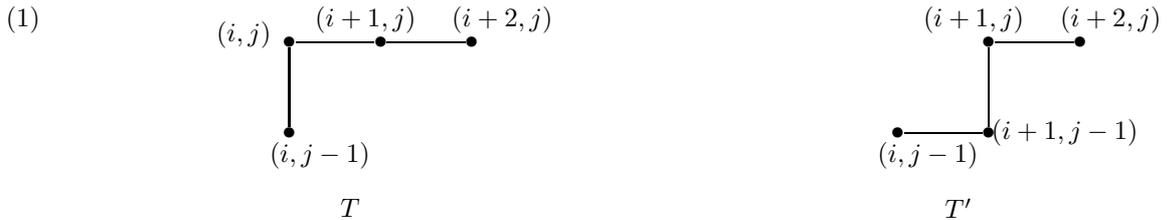
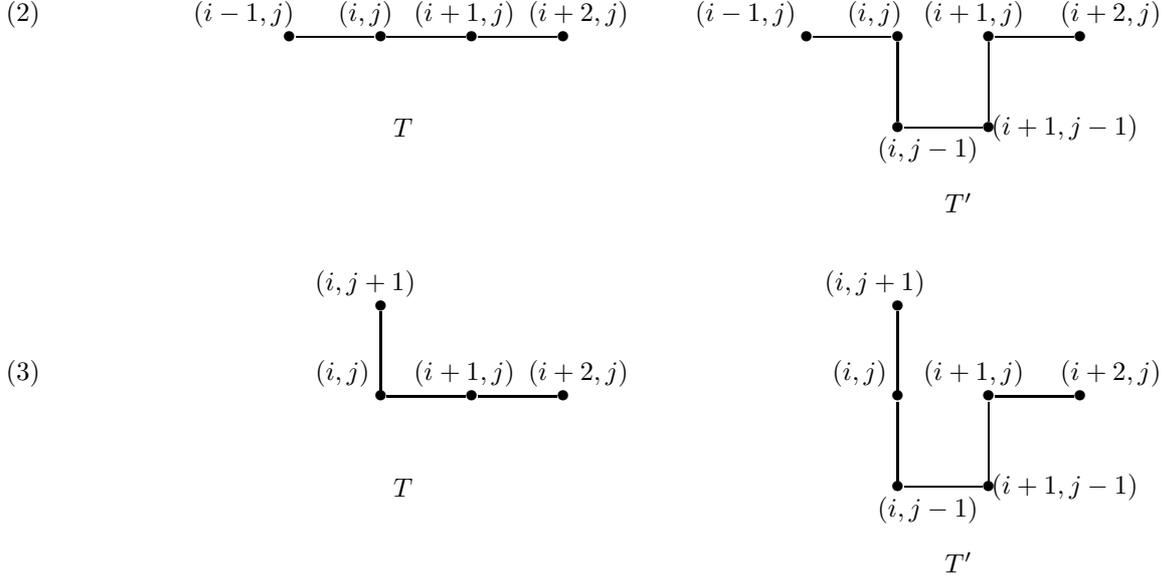
\[
\begin{xy}
(-35,15) *{(1)}="i",
(-6,13) *{(i,j)}="0000",
(4,-3) *{(i,j-1)}="000",
(10,15) *{(i+1,j)}="00",
(28,15) *{(i+2,j)}="0",
(8,-10) *{T}="T",
(102,0) *{(i+1,j-1)}="a0000",
(84,-3) *{(i,j-1)}="a000",
(90,15) *{(i+1,j)}="a00",
(108,15) *{(i+2,j)}="a0",
(88,-10) *{T'}="T'",
(0,0) *{\bullet}="1",
(12,12)*{\bullet}="3",
(0,12)*{\bullet}="4",
(24,12)*{\bullet}="5",
(80,0) *{\bullet}="a1",
(92,12)*{\bullet}="a3",
(92,0)*{\bullet}="a4",
(104,12)*{\bullet}="a5",
\ar@{-} "a1";"a4"^{}
\ar@{-} "a3";"a4"^{}
\ar@{-} "a3";"a5"^{}
\ar@{-} "1";"4"^{}
\ar@{-} "3";"4"^{}
\ar@{-} "3";"5"^{}
\end{xy}
\]

\[
\begin{xy}
(-35,15) *{(2)}="i",
(-6,15) *{(i-1,j)}="0000",
(60,15) *{(i-1,j)}="000",
(10,15) *{(i,j)}="00",
(23,15) *{(i+1,j)}="0",
(38,15) *{(i+2,j)}="00000",
(15,0) *{T}="T",
(77,15) *{(i,j)}="a00000",
(102,0) *{(i+1,j-1)}="a0000",
(84,-3) *{(i,j-1)}="a000",
(90,15) *{(i+1,j)}="a00",
(108,15) *{(i+2,j)}="a0",
(88,-10) *{T'}="T'",
(12,12)*{\bullet}="3",
(0,12)*{\bullet}="4",
(24,12)*{\bullet}="5",
(36,12)*{\bullet}="6",
(80,0) *{\bullet}="a1",
(80,12) *{\bullet}="a2",
(92,12)*{\bullet}="a3",
(92,0)*{\bullet}="a4",
(104,12)*{\bullet}="a5",
(68,12) *{\bullet}="a6",
\ar@{-} "a1";"a4"^{}
\ar@{-} "a1";"a2"^{}
\ar@{-} "a3";"a4"^{}
\ar@{-} "a3";"a5"^{}
\ar@{-} "a6";"a2"^{}
\ar@{-} "3";"4"^{}
\ar@{-} "3";"5"^{}
\ar@{-} "6";"5"^{}
\end{xy}
\]

\[
\begin{xy}
(-35,15) *{(3)}="i",
(10,27) *{(i,j+1)}="0000",
(77,27) *{(i,j+1)}="000",
(7,15) *{(i,j)}="00",
(23,15) *{(i+1,j)}="0",
(38,15) *{(i+2,j)}="00000",
(15,0) *{T}="T",
(75,15) *{(i,j)}="a00000",
(102,0) *{(i+1,j-1)}="a0000",
(84,-3) *{(i,j-1)}="a000",
(90,15) *{(i+1,j)}="a00",
(108,15) *{(i+2,j)}="a0",
(88,-10) *{T'}="T'",
(12,24)*{\bullet}="2",
(12,12)*{\bullet}="3",
(24,12)*{\bullet}="5",
(36,12)*{\bullet}="6",
(80,0) *{\bullet}="a1",
(80,12) *{\bullet}="a2",
(92,12)*{\bullet}="a3",
(92,0)*{\bullet}="a4",
(104,12)*{\bullet}="a5",
(80,24) *{\bullet}="a6",
\ar@{-} "a1";"a4"^{}
\ar@{-} "a1";"a2"^{}
\ar@{-} "a3";"a4"^{}
\ar@{-} "a3";"a5"^{}
\ar@{-} "a6";"a2"^{}
\ar@{-} "3";"2"^{}
\ar@{-} "3";"5"^{}
\ar@{-} "6";"5"^{}
\end{xy}
\]
Here, the patterns (2), (3) are happen only in the case $i-1<0$.

\vspace{2mm}

\nd
(1) First, we consider the pattern (1).

\vspace{2mm}

\nd
\underline{Case 1. $y_{i-1}=y_i-1$ and $y_{i-2}=y_{i-1}$}

\vspace{2mm}

\nd
The table of admissible and removable points are as follows:

\begin{table}[H]
  \begin{tabular}{|c|c|c|} \hline
    points & in $T$ & in $T'$ \\ \hline
    $(i,j)$ & $1$-double admissible if $i<0$, & normal  \\ 
    & $1$-single admissible if $i>0$ &
    \\ \hline
    $(i+1,j-1)$ & normal & $1$-double removable if $i+1>1$,  \\ 
    & & $1$-single removable if $i+1<1$ \\ \hline
    $(i-1,j-1)$ & normal & normal \\ \hline
    $(i,j-1)$ & $1$-double removable if $i<1$, & $1$-single removable if $i-1<0$, \\ 
    & $1$-single removable if $i>1$ & normal if $i-1>0$ \\ \hline
    $(i+2,j)$ & normal (resp. $2$-removable) & normal \\ \hline
    $(i+1,j)$ & normal & $2$-admissible (resp. normal) \\
    \hline
  \end{tabular}
\end{table}
\nd
By $i+k\equiv1$ and $\pi_1(k)\neq 1$, it holds $i\neq0,1$. Hence $i>0$ (resp. $i<0$) if and only if $i>1$ (resp. $i<1$).
Since other points in $T$ are same as in $T'$,
it follows by (\ref{L1kdef}) and $P^k(i+k)=P^k(i+k-1)$ that
\begin{eqnarray*}
L^1_{s,k,\iota}(T')-L^1_{s,k,\iota}(T)&=&
-x_{s+P^k(i+k)+[i]_-+k-j,1}
-x_{s+P^k(i+k)+[i]_-+k-j+1,1}\\
& &+x_{s+P^k(i+k+1)+[i+1]_-+k-j,2}\\
&=&-\beta_{s+P^k(i+k)+[i]_-+k-j,1}.
\end{eqnarray*}
In other cases and patterns, one can similarly show
$L^1_{s,k,\iota}(T')-L^1_{s,k,\iota}(T)=-\beta_{s+P^k(i+k)+[i]_-+k-j,1}$.
We write only the table for each case.

\vspace{2mm}

\nd
\underline{Case 2. $y_{i-1}=y_i-1$ and $y_{i-2}<y_{i-1}$}

\vspace{2mm}

\begin{table}[H]
  \begin{tabular}{|c|c|c|} \hline
    points & in $T$ & in $T'$ \\ \hline
    $(i,j)$ & $1$-double admissible if $i<0$, & normal  \\ 
    & $1$-single admissible if $i>0$ &
    \\ \hline
    $(i+1,j-1)$ & normal & $1$-double removable if $i+1>1$,  \\ 
    & & $1$-single removable if $i+1<1$ \\ \hline
    $(i-1,j-1)$ & $1$-single admissible if $i-1>0$, & $1$-double admissible if $i-1>0$,  \\ 
    & normal if $i-1<0$ & $1$-single admissible if $i-1<0$ \\ \hline
    $(i,j-1)$ & normal & normal \\ \hline
    $(i+2,j)$ & normal (resp. $2$-removable) & normal \\ \hline
    $(i+1,j)$ & normal & $2$-admissible (resp. normal) \\
    \hline
  \end{tabular}
\end{table}

\vspace{2mm}

\nd
\underline{Case 3. $y_{i-1}<y_i-1$}

\vspace{2mm}

\nd
In this case, we have $i-1>0$.
The table of admissible and removable points are as follows:

\begin{table}[H]
  \begin{tabular}{|c|c|c|} \hline
    points & in $T$ & in $T'$ \\ \hline
    $(i,j)$ & $1$-single admissible & normal  \\ \hline
    $(i+1,j-1)$ & normal & $1$-single removable \\ \hline
    $(i-1,j-1)$ & normal & normal \\ \hline
    $(i,j-1)$ & normal & normal \\ \hline
    $(i+2,j)$ & normal (resp. $2$-removable) & normal \\ \hline
    $(i+1,j)$ & normal & $2$-admissible (resp. normal) \\
    \hline
  \end{tabular}
\end{table}

(2) Next, we consider the pattern (2). It holds $i<1$.

\vspace{2mm}

\nd
\underline{Case 1. $y_{i-2}=y_{i-1}$}

\vspace{2mm}

\begin{table}[H]
  \begin{tabular}{|c|c|c|} \hline
    points & in $T$ & in $T'$ \\ \hline
    $(i,j)$ & $1$-admissible and $1$-removable & $1$-removable  \\ \hline
    $(i+1,j-1)$ & normal & $1$-removable \\ \hline
    $(i-1,j)$ & normal & normal  \\ \hline
    $(i,j-1)$ & normal & normal \\ \hline
    $(i+2,j)$ & normal (resp. $2$-removable) & normal \\ \hline
    $(i+1,j)$ & normal & $2$-admissible (resp. normal) \\
    \hline
  \end{tabular}
\end{table}

\vspace{2mm}

\nd
\underline{Case 2. $y_{i-2}<y_{i-1}$}

\vspace{2mm}

\begin{table}[H]
  \begin{tabular}{|c|c|c|} \hline
    points & in $T$ & in $T'$ \\ \hline
    $(i,j)$ & $1$-admissible & normal  \\ \hline
    $(i+1,j-1)$ & normal & $1$-removable \\ \hline
    $(i-1,j)$ & $1$-admissible & $1$-admissible \\ \hline
    $(i,j-1)$ & normal & normal \\ \hline
    $(i+2,j)$ & normal (resp. $2$-removable) & normal \\ \hline
    $(i+1,j)$ & normal & $2$-admissible (resp. normal) \\
    \hline
  \end{tabular}
\end{table}

(3) Finally, we consider the pattern (3). It holds $i<1$. The table is as follows:

\begin{table}[H]
  \begin{tabular}{|c|c|c|} \hline
    points & in $T$ & in $T'$ \\ \hline
    $(i,j)$ & $1$-admissible & normal  \\ \hline
    $(i+1,j-1)$ & normal & $1$-removable \\ \hline
    $(i+2,j)$ & normal (resp. $2$-removable) & normal \\ \hline
    $(i+1,j)$ & normal & $2$-admissible (resp. normal) \\
    \hline
  \end{tabular}
\end{table}

\qed

\begin{lem}\label{lemA2-4}
In the case $i+k\equiv 0$ $({\rm mod}\ 2n-1)$ so that $\pi_1(i+k)=1$, Proposition \ref{A2closed} (i) holds.
\end{lem}

\nd
{\it Proof.}
Since $i+k\equiv 0$ (${\rm mod}\ 2n-1$), there are three patterns (1) $y_{i}=y_{i+1}$, (2) $y_{i}<y_{i+1}$, (3) $y_{i}>y_{i+1}$ (Definition \ref{AEYD}).
We get $y_{i}=y_{i-1}+1$.
Putting $j:=y_i$,
lines around the point $(i,j)$ in $T$ and $T'$ of each pattern are as follows:
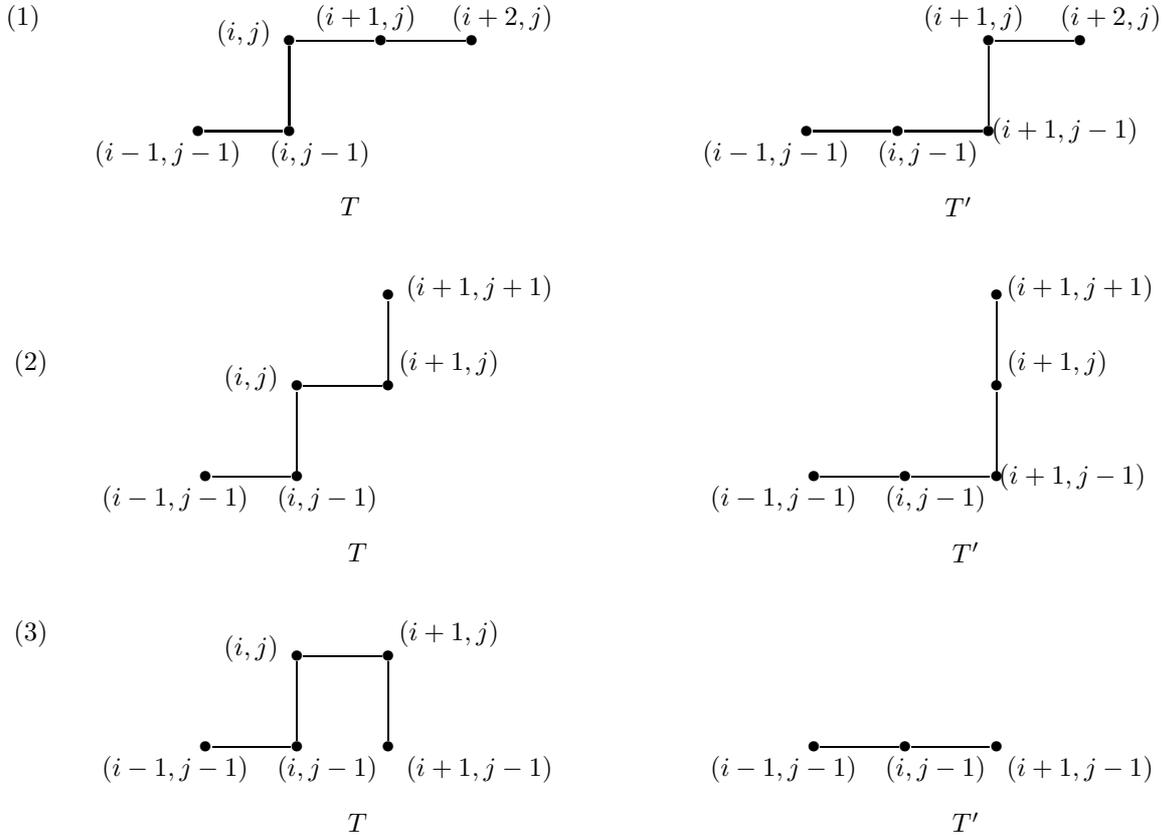
\[
\begin{xy}
(-35,15) *{(1)}="i",
(-16,-3) *{(i-1,j-1)}="00000",
(-6,13) *{(i,j)}="0000",
(4,-3) *{(i,j-1)}="000",
(10,15) *{(i+1,j)}="00",
(28,15) *{(i+2,j)}="0",
(8,-10) *{T}="T",
(64,-3) *{(i-1,j-1)}="a0000",
(102,0) *{(i+1,j-1)}="a0000",
(84,-3) *{(i,j-1)}="a000",
(90,15) *{(i+1,j)}="a00",
(108,15) *{(i+2,j)}="a0",
(88,-10) *{T'}="T'",
(-12,0) *{\bullet}="11",
(0,0) *{\bullet}="1",
(12,12)*{\bullet}="3",
(0,12)*{\bullet}="4",
(24,12)*{\bullet}="5",
(68,0) *{\bullet}="a11",
(80,0) *{\bullet}="a1",
(92,12)*{\bullet}="a3",
(92,0)*{\bullet}="a4",
(104,12)*{\bullet}="a5",
\ar@{-} "a1";"a11"^{}
\ar@{-} "a1";"a4"^{}
\ar@{-} "a3";"a4"^{}
\ar@{-} "a3";"a5"^{}
\ar@{-} "1";"11"^{}
\ar@{-} "1";"4"^{}
\ar@{-} "3";"4"^{}
\ar@{-} "3";"5"^{}
\end{xy}
\]

\[
\begin{xy}
(-35,15) *{(2)}="i",
(-16,-3) *{(i-1,j-1)}="00000",
(-6,13) *{(i,j)}="0000",
(4,-3) *{(i,j-1)}="000",
(20,15) *{(i+1,j)}="00",
(24,25) *{(i+1,j+1)}="0",
(8,-10) *{T}="T",
(64,-3) *{(i-1,j-1)}="a0000",
(102,0) *{(i+1,j-1)}="a0000",
(84,-3) *{(i,j-1)}="a000",
(100,15) *{(i+1,j)}="a00",
(103,25) *{(i+1,j+1)}="a0",
(88,-10) *{T'}="T'",
(-12,0) *{\bullet}="11",
(0,0) *{\bullet}="1",
(12,12)*{\bullet}="3",
(0,12)*{\bullet}="4",
(12,24)*{\bullet}="5",
(68,0) *{\bullet}="a11",
(80,0) *{\bullet}="a1",
(92,12)*{\bullet}="a3",
(92,0)*{\bullet}="a4",
(92,24)*{\bullet}="a5",
\ar@{-} "a1";"a11"^{}
\ar@{-} "a1";"a4"^{}
\ar@{-} "a3";"a4"^{}
\ar@{-} "a3";"a5"^{}
\ar@{-} "1";"11"^{}
\ar@{-} "1";"4"^{}
\ar@{-} "3";"4"^{}
\ar@{-} "3";"5"^{}
\end{xy}
\]

\[
\begin{xy}
(-35,15) *{(3)}="i",
(-16,-3) *{(i-1,j-1)}="00000",
(-6,13) *{(i,j)}="0000",
(4,-3) *{(i,j-1)}="000",
(20,15) *{(i+1,j)}="00",
(24,-3) *{(i+1,j-1)}="0",
(8,-10) *{T}="T",
(64,-3) *{(i-1,j-1)}="a0000",
(84,-3) *{(i,j-1)}="a000",
(103,-3) *{(i+1,j-1)}="a0",
(88,-10) *{T'}="T'",
(-12,0) *{\bullet}="11",
(0,0) *{\bullet}="1",
(12,12)*{\bullet}="3",
(0,12)*{\bullet}="4",
(12,0)*{\bullet}="5",
(68,0) *{\bullet}="a11",
(80,0) *{\bullet}="a1",
(92,0)*{\bullet}="a4",
\ar@{-} "a1";"a11"^{}
\ar@{-} "a1";"a4"^{}
\ar@{-} "1";"11"^{}
\ar@{-} "1";"4"^{}
\ar@{-} "3";"4"^{}
\ar@{-} "3";"5"^{}
\end{xy}
\]
Here, the pattern (2) (resp. (3)) happens only in the case $i>0$ (resp. $i<0$).

\nd
(1) First, we consider the pattern (1).

\vspace{2mm}

\nd
\underline{Case 1. $y_{i+1}=y_{i+2}$}

\vspace{2mm}

\begin{table}[H]
  \begin{tabular}{|c|c|c|} \hline
    points & in $T$ & in $T'$ \\ \hline
    $(i,j)$ & $1$-double admissible if $i>0$, & normal \\ 
    & $1$-single admissible if $i<0$ & 
    \\ \hline
    $(i+1,j-1)$ & normal & $1$-double removable if $i+1<1$, \\ 
    & & $1$-single removable if $i+1>1$ \\ \hline
    $(i+1,j)$ & $1$-single admissible if $i<0$, & $1$-double admissible if $i+1<0$, \\
    & normal if $i>0$ & $1$-single admissible if $i+1>0$ \\ \hline
    $(i+2,j)$ & normal & normal \\ \hline
    $(i,j-1)$ & normal (resp. $2$-removable) & normal \\ \hline
   $(i-1,j-1)$ & normal & $2$-admissible (resp. normal)  \\ \hline
   \end{tabular}
\end{table}
The conditions $i+k\equiv 0$ (${\rm mod}\ 2n-1$) and $\pi_1(k)\neq1$ mean $i\neq 0,-1$ so that
$i>0$ (resp. $i<0$) if and only if $i+1>0$ (resp. $i+1<0$).
Considering the definition (\ref{L1kdef}) and $P^k(i+k)=P^k(i+k+1)$, we get
\begin{eqnarray*}
L^1_{s,k,\iota}(T')-L^1_{s,k,\iota}(T)&=&
-x_{s+P^k(i+k)+[i]_-+k-j,1}
-x_{s+P^k(i+k)+[i]_-+k-j+1,1}\\
& &+x_{s+P^k(i+k-1)+[i-1]_-+k-j+1,2}\\
&=&-\beta_{s+P^k(i+k)+[i]_-+k-j,1}.
\end{eqnarray*}
In other cases and patterns, one can similarly verify
$L^1_{s,k,\iota}(T')-L^1_{s,k,\iota}(T)=-\beta_{s+P^k(i+k)+[i]_-+k-j,1}$.
We write only the table for each case.

\vspace{2mm}

\nd
\underline{Case 2. $y_{i+1}<y_{i+2}$}

\vspace{2mm}

\begin{table}[H]
  \begin{tabular}{|c|c|c|} \hline
    points & in $T$ & in $T'$ \\ \hline
    $(i,j)$ & $1$-double admissible if $i>0$, & normal \\ 
    & $1$-single admissible if $i<0$ & 
    \\ \hline
    $(i+1,j-1)$ & normal & $1$-double removable if $i+1<1$, \\ 
    & & $1$-single removable if $i+1>1$ \\ \hline
    $(i+1,j)$ & normal & normal \\ \hline
    $(i+2,j)$ & $1$-double removable if $i+2>1$, & $1$-single removable if $i>0$, \\ 
    & $1$-single removable if $i+2<1$ & normal if $i<0$ \\ \hline
    $(i,j-1)$ & normal (resp. $2$-removable) & normal \\ \hline
   $(i-1,j-1)$ & normal & $2$-admissible (resp. normal)  \\ \hline
   \end{tabular}
\end{table}

(2) Next, we consider the pattern (2). It holds $i>0$.

\begin{table}[H]
  \begin{tabular}{|c|c|c|} \hline
    points & in $T$ & in $T'$ \\ \hline
    $(i,j)$ & $1$-single admissible & normal \\ \hline
    $(i+1,j-1)$ & normal & $1$-single removable \\ \hline
    $(i+1,j)$ & normal & normal \\ \hline
    $(i,j-1)$ & normal (resp. $2$-removable) & normal \\ \hline
   $(i-1,j-1)$ & normal & $2$-admissible (resp. normal)  \\ \hline
   \end{tabular}
\end{table}

(3) Next, we consider the pattern (3). It holds $i<0$.

\vspace{2mm}

\nd
\underline{Case 1. $y_{i+1}<y_i-1$}

\vspace{2mm}

\begin{table}[H]
  \begin{tabular}{|c|c|c|} \hline
    points & in $T$ & in $T'$ \\ \hline
    $(i,j)$ & $1$-single admissible & normal \\ \hline
    $(i+1,j-1)$ & normal & $1$-single removable \\ \hline
    $(i+1,j)$ & normal & normal \\ \hline
    $(i,j-1)$ & normal (resp. $2$-removable) & normal \\ \hline
   $(i-1,j-1)$ & normal & $2$-admissible (resp. normal)  \\ \hline
   \end{tabular}
\end{table}

\vspace{2mm}

\nd
\underline{Case 2. $y_{i+1}=y_i-1$ and $y_{i+1}=y_{i+2}$}

\vspace{2mm}

\begin{table}[H]
  \begin{tabular}{|c|c|c|} \hline
    points & in $T$ & in $T'$ \\ \hline
    $(i,j)$ & $1$-single admissible & normal \\ \hline
    $(i+1,j-1)$ & $1$-single admissible & $1$-single admissible and $1$-single removable \\ \hline
    $(i+1,j)$ & normal & normal \\ \hline
    $(i,j-1)$ & normal (resp. $2$-removable) & normal \\ \hline
   $(i-1,j-1)$ & normal & $2$-admissible (resp. normal)  \\ \hline
   \end{tabular}
\end{table}

\vspace{2mm}

\nd
\underline{Case 3. $y_{i+1}=y_i-1$ and $y_{i+1}<y_{i+2}$}

\vspace{2mm}

\begin{table}[H]
  \begin{tabular}{|c|c|c|} \hline
    points & in $T$ & in $T'$ \\ \hline
    $(i,j)$ & $1$-single admissible & normal \\ \hline
    $(i+1,j-1)$ & normal & $1$-single removable \\ \hline
    $(i+1,j)$ & normal & normal \\ \hline
    $(i,j-1)$ & normal (resp. $2$-removable) & normal \\ \hline
   $(i-1,j-1)$ & normal & $2$-admissible (resp. normal)  \\ \hline
   \end{tabular}
\end{table}

\qed

Hence, Proposition \ref{A2closed} follows from Lemma \ref{lemA2-1}-\ref{lemA2-4}.

\subsection{Type ${\rm A}_{2n-2}^{(2)}$-case (action on Young walls)}

In this subsection, we assume $\mathfrak{g}$ is of type ${\rm A}_{2n-2}^{(2)}$. Recall that we defined
the notation $L^1_{s,1,\iota}$
in (\ref{L11-def}).

\begin{prop}\label{prop-closednessAw-YW}
\begin{enumerate}
\item
Let $t\in I\setminus\{1\}$ and
we suppose that $Y\in {\rm YW}_{{\rm A}^{(2)},1}$ has a $t$-admissible slot
\[
\begin{xy}
(-8,15) *{(-i-1,l+1)}="0000",
(17,15) *{(-i,l+1)}="000",
(15,-3) *{(-i,l)}="00",
(-7,-3) *{(-i-1,l)}="0",
(0,0) *{}="1",
(12,0)*{}="2",
(12,12)*{}="3",
(0,12)*{}="4",
\ar@{--} "1";"2"^{}
\ar@{--} "1";"4"^{}
\ar@{--} "2";"3"^{}
\ar@{--} "3";"4"^{}
\end{xy}
\]
Let $Y'\in{\rm YW}_{{\rm A}^{(2)},1}$ be the proper Young wall obtained from $Y$
by adding the $t$-block to the slot.
Then for $s\in\mathbb{Z}_{\geq1}$ it follows
\[
L^1_{s,1,\iota}(Y')=L^1_{s,1,\iota}(Y)-\beta_{s+P^1(l)+i,t}.
\]
\item We suppose that $Y\in {\rm YW}_{{\rm A}^{(2)},1}$ has a $1$-admissible slot
\[
\begin{xy}
(-8,9) *{(-i-1,l+\frac{1}{2})}="0000",
(17,9) *{(-i,l+\frac{1}{2})}="000",
(15,-3) *{(-i,l)}="00",
(-7,-3) *{(-i-1,l)}="0",
(0,0) *{}="1",
(12,0)*{}="2",
(12,6)*{}="3",
(0,6)*{}="4",
(38,3)*{{\rm or}}="or",
(60,9) *{(-i-1,l+1)}="a0000",
(85,9) *{(-i,l+1)}="a000",
(83,-3) *{(-i,l+\frac{1}{2})}="a00",
(61,-3) *{(-i-1,l+\frac{1}{2})}="a0",
(68,0) *{}="a1",
(80,0)*{}="a2",
(80,6)*{}="a3",
(68,6)*{}="a4",
\ar@{--} "1";"2"^{}
\ar@{--} "1";"4"^{}
\ar@{--} "2";"3"^{}
\ar@{--} "3";"4"^{}
\ar@{--} "a1";"a2"^{}
\ar@{--} "a1";"a4"^{}
\ar@{--} "a2";"a3"^{}
\ar@{--} "a3";"a4"^{}
\end{xy}
\]
Let $Y'\in{\rm YW}_{{\rm A}^{(2)},1}$ be the proper Young wall obtained from $Y$
by adding the $1$-block to the slot.
Then for $s\in\mathbb{Z}_{\geq1}$ it follows
\[
L^1_{s,1,\iota}(Y')=L^1_{s,1,\iota}(Y)-\beta_{s+P^1(l)+i,1}.
\]
\end{enumerate}
\end{prop}

\nd
{\it Proof.} (i) Let $A$ be the $t$-admissible slot in our claim. By (\ref{piprime}), (\ref{t-pi}), we get $t=\pi'(l)$. 
It follows from $t>1$ that $l>1$.
First, we suppose that $2<t\leq n$. Since the slot $A$ is $t$-admissible, by Definition \ref{def-YW} and \ref{def-YW2}, blocks and slots around $A$ are as follows:
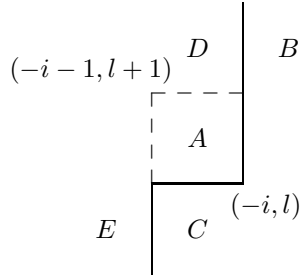
\[
\begin{xy}
(-8,15) *{(-i-1,l+1)}="0000",
(15,-3) *{(-i,l)}="00",
(6,6) *{A}="A",
(18,18) *{B}="B",
(6,-6) *{C}="C",
(6,18) *{D}="D",
(-6,-6) *{E}="E",
(0,0) *{}="1",
(12,0)*{}="2",
(12,12)*{}="3",
(0,12)*{}="4",
(12,24)*{}="5",
(0,-12)*{}="6",
\ar@{-} "1";"2"^{}
\ar@{--} "1";"4"^{}
\ar@{-} "2";"3"^{}
\ar@{--} "3";"4"^{}
\ar@{-} "3";"5"^{}
\ar@{-} "1";"6"^{}
\end{xy}
\]
Here, the blocks and slots around the slot $A$ are named as $B$, $C$, $D$ and $E$ as above. If $i=0$ then we identify
$B$ as a non-removable block.
Note that $B$ is removable in $Y$ if and only if $D$ is not admissible in $Y'$.
We also see that $C$ is not removable in $Y$ if and only if $E$ is admissible in $Y'$. 
One can summarize them as follows:

\begin{table}[htb]
  \begin{tabular}{|c|c|c|} \hline
    slot or block & in $Y$ & in $Y'$ \\ \hline
    $A$ & $t$-admissible & removable $t$ \\
    $B$ & removable $\pi'(l+1)$ (resp. normal) & normal \\
    $D$ & normal & normal (resp. $\pi'(l+1)$-admissible) \\
    \hline
  \end{tabular}
\end{table}

\begin{table}[htb]
  \begin{tabular}{|c|c|c|} \hline
    slot or block & in $Y$ & in $Y'$ \\ \hline
    $C$ & normal (resp. removable $\pi'(l-1)$) & normal \\
    $E$ & normal & $\pi'(l-1)$-admissible (resp. normal) \\
    \hline
  \end{tabular}
\end{table}
\nd
Here `normal' means it is neither admissible nor removable. Note that in the case $t<n$, the relation $t=\pi'(l)$ implies
$\pi'(l-1)=t-1$ (resp. $\pi'(l-1)=t+1$) if and only if $\pi'(l+1)=t+1$ (resp. $\pi'(l+1)=t-1$).
In the case $t=n$, it follows $\pi'(l-1)=\pi'(l+1)=n-1$.
Since other blocks and slots do not change, it holds
\begin{eqnarray*}
L^1_{s,1,\iota}(Y')-L^1_{s,1,\iota}(Y)&=& 
-x_{s+P^1(l)+i,t}-x_{s+P^1(l)+i+1,t}
+x_{s+P^1(l+1)+i,\pi'(l+1)}+x_{s+P^1(l-1)+i+1,\pi'(l-1)}\\
&=& 
-x_{s+P^1(l)+i,\pi'(l)}-x_{s+P^1(l)+i+1,\pi'(l)}\\
& &+x_{s+P^1(l)+p_{\pi'(l+1),\pi'(l)}+i,\pi'(l+1)}
+x_{s+P^1(l)+p_{\pi'(l-1),\pi'(l)}+i,\pi'(l-1)}\\
&=&-\beta_{s+P^1(l)+i,\pi'(l)}=-\beta_{s+P^1(l)+i,t}.
\end{eqnarray*}

Next, we suppose that $t=2$. If $\pi'(l+1)=1$ so that $\pi'(l-1)=3$ then
blocks and slots around $A$ are as follows:
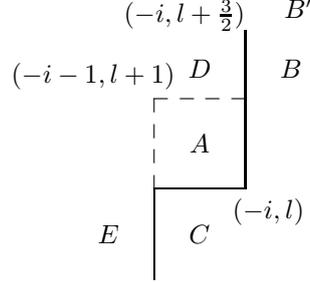
\[
\begin{xy}
(-8,15) *{(-i-1,l+1)}="0000",
(15,-3) *{(-i,l)}="00",
(4,23) *{(-i,l+\frac{3}{2})}="000",
(6,6) *{A}="A",
(18,16) *{B}="B",
(19,24) *{B'}="B'",
(6,-6) *{C}="C",
(6,16) *{D}="D",
(-6,-6) *{E}="E",
(0,0) *{}="1",
(12,0)*{}="2",
(12,12)*{}="3",
(0,12)*{}="4",
(12,21)*{}="5",
(0,-12)*{}="6",
\ar@{-} "1";"2"^{}
\ar@{--} "1";"4"^{}
\ar@{-} "2";"3"^{}
\ar@{--} "3";"4"^{}
\ar@{-} "3";"5"^{}
\ar@{-} "1";"6"^{}
\end{xy}
\]
\nd
We see that $C$ is not a removable $3$-block in $Y$ if and only if $E$ is a $3$-admissible slot in $Y'$.
As for $B$, $D$ and $B'$, there are three patterns:
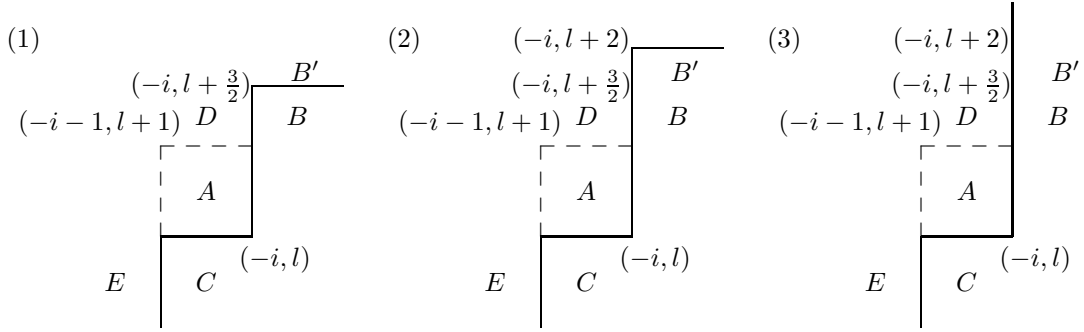
\[
\begin{xy}
(-68,26) *{(1)}="l",
(-18,26) *{(2)}="lr",
(32,26) *{(3)}="r",
(-58,15) *{(-i-1,l+1)}="0000l",
(-35,-3) *{(-i,l)}="00l",
(-46,20) *{(-i,l+\frac{3}{2})}="000l",
(-44,6) *{A}="Al",
(-32,16) *{B}="Bl",
(-31,22) *{B'}="B'l",
(-44,-6) *{C}="Cl",
(-44,16) *{D}="Dl",
(-56,-6) *{E}="El",
(-50,0) *{}="1l",
(-38,0)*{}="2l",
(-38,12)*{}="3l",
(-50,12)*{}="4l",
(-38,20)*{}="5l",
(-26,20)*{}="55l",
(-50,-12)*{}="6l",
(-8,15) *{(-i-1,l+1)}="0000",
(15,-3) *{(-i,l)}="00",
(4,20) *{(-i,l+\frac{3}{2})}="000",
(4,26) *{(-i,l+2)}="00000",
(6,6) *{A}="A",
(18,16) *{B}="B",
(19,22) *{B'}="B'",
(6,-6) *{C}="C",
(6,16) *{D}="D",
(-6,-6) *{E}="E",
(0,0) *{}="1",
(12,0)*{}="2",
(12,12)*{}="3",
(0,12)*{}="4",
(12,25)*{}="5",
(24,25)*{}="55",
(0,-12)*{}="6",
(42,15) *{(-i-1,l+1)}="0000r",
(65,-3) *{(-i,l)}="00r",
(54,20) *{(-i,l+\frac{3}{2})}="000r",
(54,26) *{(-i,l+2)}="00000r",
(56,6) *{A}="Ar",
(68,16) *{B}="Br",
(69,22) *{B'}="B'r",
(56,-6) *{C}="Cr",
(56,16) *{D}="Dr",
(44,-6) *{E}="Er",
(50,0) *{}="1r",
(62,0)*{}="2r",
(62,12)*{}="3r",
(50,12)*{}="4r",
(62,25)*{}="5r",
(62,31)*{}="55r",
(50,-12)*{}="6r",
\ar@{-} "1l";"2l"^{}
\ar@{--} "1l";"4l"^{}
\ar@{-} "2l";"3l"^{}
\ar@{--} "3l";"4l"^{}
\ar@{-} "3l";"5l"^{}
\ar@{-} "5l";"55l"^{}
\ar@{-} "1l";"6l"^{}
\ar@{-} "1";"2"^{}
\ar@{--} "1";"4"^{}
\ar@{-} "2";"3"^{}
\ar@{--} "3";"4"^{}
\ar@{-} "3";"5"^{}
\ar@{-} "5";"55"^{}
\ar@{-} "1";"6"^{}
\ar@{-} "1r";"2r"^{}
\ar@{--} "1r";"4r"^{}
\ar@{-} "2r";"3r"^{}
\ar@{--} "3r";"4r"^{}
\ar@{-} "3r";"5r"^{}
\ar@{-} "5r";"55r"^{}
\ar@{-} "1r";"6r"^{}
\end{xy}
\]
\nd
In the pattern (1),
$B$ is a single removable $1$-block in $Y$ and is a non-removable block in $Y'$,
$D$ is a non-admissible slot in $Y$ and is a single $1$-admissible slot in $Y'$.
$B'$ is a single $1$-admissible slot or non-admissible slot but, the admissibility in $Y$ is same as in $Y'$.
Considering (2), (3) similarly,
one can summarize removable blocks and admissible slots for each pattern as follows:
\begin{table}[H]
  \begin{tabular}{|c|l|l|} \hline
    slot or block & in $Y$ & in $Y'$ \\ \hline
    $A$ & $2$-admissible & removable $2$ \\
        $C$ & normal (resp. removable $3$) & normal \\
    $E$ & normal & $3$-admissible (resp. normal) \\
    $B$ & (1) single removable $1$, & (1), (2), (3) normal \\
     & (2), (3) normal & \\
    $D$ & (1), (2), (3) normal & (1), (2) single $1$-admissible, \\
     & & (3) double $1$-admissible \\
    $B'$ & (1) single $1$-admissible or normal, & (1) same as in $Y$,  \\
     & (2) double removable $1$, & (2) single removable $1$,  \\
     & (3) normal & (3) normal  \\
    \hline
  \end{tabular}
\end{table}
Thus, 
\begin{eqnarray*}
L^1_{s,1,\iota}(Y')-L^1_{s,1,\iota}(Y)&=& 
-x_{s+P^1(l)+i,2}-x_{s+P^1(l)+i+1,2}
+2x_{s+P^1(l+1)+i,1}+x_{s+P^1(l-1)+i+1,3}\\
&=& 
-x_{s+P^1(l)+i,2}-x_{s+P^1(l)+i+1,2}\\
& &+2x_{s+P^1(l)+p_{1,2}+i,1}
+x_{s+P^1(l)+p_{3,2}+i,3}\\
&=&-\beta_{s+P^1(l)+i,2}.
\end{eqnarray*}

If $\pi'(l+1)=3$ so that $\pi'(l-1)=1$ then 
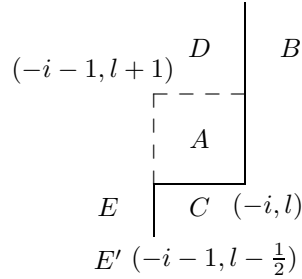
\[
\begin{xy}
(-8,15) *{(-i-1,l+1)}="0000",
(8,-10) *{(-i-1,l-\frac{1}{2})}="000",
(15,-3) *{(-i,l)}="00",
(6,6) *{A}="A",
(18,18) *{B}="B",
(6,-3) *{C}="C",
(6,18) *{D}="D",
(-6,-3) *{E}="E",
(-6,-10) *{E'}="E'",
(0,0) *{}="1",
(12,0)*{}="2",
(12,12)*{}="3",
(0,12)*{}="4",
(12,24)*{}="5",
(0,-7)*{}="6",
\ar@{-} "1";"2"^{}
\ar@{--} "1";"4"^{}
\ar@{-} "2";"3"^{}
\ar@{--} "3";"4"^{}
\ar@{-} "3";"5"^{}
\ar@{-} "1";"6"^{}
\end{xy}
\]
We can similarly verify
\[
L^1_{s,1,\iota}(Y')-L^1_{s,1,\iota}(Y)=-\beta_{s+P^1(l)+i,2}.
\]

(ii) In this case, it holds $\pi'(l)=1$. First, we consider the case the $1$-admissible slot is
\begin{equation}\label{t1dia}
\begin{xy}
(-8,9) *{(-i-1,l+\frac{1}{2})}="0000",
(17,9) *{(-i,l+\frac{1}{2})}="000",
(15,-3) *{(-i,l)}="00",
(-7,-3) *{(-i-1,l)}="0",
(0,0) *{}="1",
(12,0)*{}="2",
(12,6)*{}="3",
(0,6)*{}="4",
\ar@{--} "1";"2"^{}
\ar@{--} "1";"4"^{}
\ar@{--} "2";"3"^{}
\ar@{--} "3";"4"^{}
\end{xy}
\end{equation}
Let $A$ be this slot.
It holds $\pi'(l-1)=2$. Then blocks and slots around $A$ are 
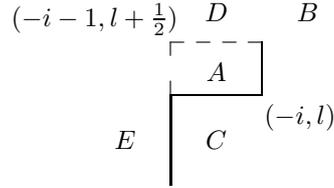
\[
\begin{xy}
(-10,13) *{(-i-1,l+\frac{1}{2})}="0000",
(17,0) *{(-i,l)}="00",
(6,6) *{A}="A",
(18,14) *{B}="B",
(6,-3) *{C}="C",
(6,14) *{D}="D",
(-6,-3) *{E}="E",
(0,3) *{}="1",
(12,3)*{}="2",
(12,10)*{}="3",
(0,10)*{}="4",
(0,-9)*{}="6",
\ar@{-} "1";"2"^{}
\ar@{--} "1";"4"^{}
\ar@{-} "2";"3"^{}
\ar@{--} "3";"4"^{}
\ar@{-} "1";"6"^{}
\end{xy}
\]
We see that $C$ is not a removable $2$-block in $Y$ if and only if $E$ is a $2$-admissible slot in $Y'$.
As for $A$, $B$ and $D$, there are three patterns:
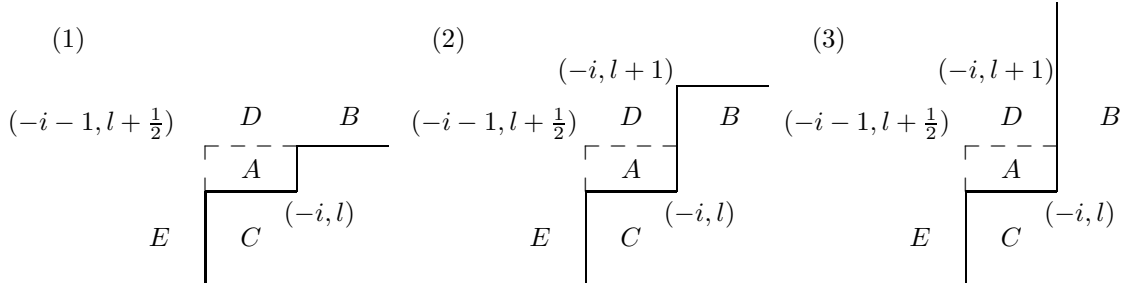
\[
\begin{xy}
(-68,26) *{(1)}="l",
(-18,26) *{(2)}="lr",
(32,26) *{(3)}="r",
(-65,15) *{(-i-1,l+\frac{1}{2})}="0000l",
(-35,3) *{(-i,l)}="00l",
(-44,9) *{A}="Al",
(-31,16) *{B}="B'l",
(-44,0) *{C}="Cl",
(-44,16) *{D}="Dl",
(-56,0) *{E}="El",
(-50,6) *{}="1l",
(-38,6)*{}="2l",
(-38,12)*{}="3l",
(-50,12)*{}="4l",
(-38,12)*{}="5l",
(-26,12)*{}="55l",
(-50,-6)*{}="6l",
(-12,15) *{(-i-1,l+\frac{1}{2})}="0000",
(15,3) *{(-i,l)}="00",
(4,22) *{(-i,l+1)}="00000",
(6,9) *{A}="A",
(19,16) *{B}="B'",
(6,0) *{C}="C",
(6,16) *{D}="D",
(-6,0) *{E}="E",
(0,6) *{}="1",
(12,6)*{}="2",
(12,12)*{}="3",
(0,12)*{}="4",
(12,20)*{}="5",
(24,20)*{}="55",
(0,-6)*{}="6",
(37,15) *{(-i-1,l+\frac{1}{2})}="0000r",
(65,3) *{(-i,l)}="00r",
(54,22) *{(-i,l+1)}="00000r",
(56,9) *{A}="Ar",
(69,16) *{B}="B'r",
(56,0) *{C}="Cr",
(56,16) *{D}="Dr",
(44,0) *{E}="Er",
(50,6) *{}="1r",
(62,6)*{}="2r",
(62,12)*{}="3r",
(50,12)*{}="4r",
(62,25)*{}="5r",
(62,31)*{}="55r",
(50,-6)*{}="6r",
\ar@{-} "1l";"2l"^{}
\ar@{--} "1l";"4l"^{}
\ar@{-} "2l";"3l"^{}
\ar@{--} "3l";"4l"^{}
\ar@{-} "3l";"5l"^{}
\ar@{-} "5l";"55l"^{}
\ar@{-} "1l";"6l"^{}
\ar@{-} "1";"2"^{}
\ar@{--} "1";"4"^{}
\ar@{-} "2";"3"^{}
\ar@{--} "3";"4"^{}
\ar@{-} "3";"5"^{}
\ar@{-} "5";"55"^{}
\ar@{-} "1";"6"^{}
\ar@{-} "1r";"2r"^{}
\ar@{--} "1r";"4r"^{}
\ar@{-} "2r";"3r"^{}
\ar@{--} "3r";"4r"^{}
\ar@{-} "3r";"5r"^{}
\ar@{-} "5r";"55r"^{}
\ar@{-} "1r";"6r"^{}
\end{xy}
\]
\nd
Removable blocks and admissible slots in each pattern are as follows:
\begin{table}[H]
  \begin{tabular}{|c|l|l|} \hline
    slot or block & in $Y$ & in $Y'$ \\ \hline
    $A$ & (1), (2) single $1$-admissible, & (1), (2), (3) single removable $1$ \\
     & (3) double $1$-admissible & \\
    $B$ & (1) single $1$-admissible or normal & (1) same as in $Y$, \\
    & (2) single removable $1$, & (2) single removable $1$, \\
    & (3) normal & (3) normal \\
    $D$ & (1), (2), (3) normal & (1), (2) normal, \\
     &  & (3) single $1$-admissible \\
    \hline
  \end{tabular}
\end{table}

\begin{table}[htb]
  \begin{tabular}{|c|c|c|} \hline
    slot or block & in $Y$ & in $Y'$ \\ \hline
    $C$ & normal (resp. removable $2$) & normal \\
    $E$ & normal & $2$-admissible (resp. normal) \\
    \hline
  \end{tabular}
\end{table}
Thus, 
\begin{eqnarray*}
L^1_{s,1,\iota}(Y')-L^1_{s,1,\iota}(Y)&=& 
-x_{s+P^1(l)+i,1}-x_{s+P^1(l)+i+1,1}
+x_{s+P^1(l-1)+i+1,2}\\
&=& 
-x_{s+P^1(l)+i,1}-x_{s+P^1(l)+i+1,1}
+x_{s+P^1(l)+p_{2,1}+i,2}\\
&=&-\beta_{s+P^1(l)+i,1}.
\end{eqnarray*}

Next, we consider the case the slot is
\[
\begin{xy}
(60,9) *{(-i-1,l+1)}="a0000",
(85,9) *{(-i,l+1)}="a000",
(83,-3) *{(-i,l+\frac{1}{2})}="a00",
(61,-3) *{(-i-1,l+\frac{1}{2})}="a0",
(68,0) *{}="a1",
(80,0)*{}="a2",
(80,6)*{}="a3",
(68,6)*{}="a4",
\ar@{--} "a1";"a2"^{}
\ar@{--} "a1";"a4"^{}
\ar@{--} "a2";"a3"^{}
\ar@{--} "a3";"a4"^{}
\end{xy}
\]
We obtain $\pi'(l+1)=2$.
Let $A$ be this slot.
Then blocks and slots around $A$ are 
\[
\begin{xy}
(-10,13) *{(-i-1,l+1)}="0000",
(20,0) *{(-i,l+\frac{1}{2})}="00",
(6,6) *{A}="A",
(18,16) *{B}="B",
(6,0) *{C}="C",
(6,16) *{D}="D",
(-6,0) *{E}="E",
(-6,-6) *{F}="F",
(0,3) *{}="1",
(12,3)*{}="2",
(12,10)*{}="3",
(0,10)*{}="4",
(12,22)*{}="5",
(0,-3)*{}="6",
\ar@{-} "1";"2"^{}
\ar@{--} "1";"4"^{}
\ar@{-} "2";"5"^{}
\ar@{--} "3";"4"^{}
\end{xy}
\]
There are three patterns : In $Y$, (1) $E$ is a $1$-block, (2) $E$ is a single $1$-admissible slot,
(3) $F$ is a slot so that $E$ is a non-admissible slot. As with the case the slot is (\ref{t1dia}),
we get the following tables:
\begin{table}[H]
  \begin{tabular}{|c|l|l|} \hline
    slot or block & in $Y$ & in $Y'$ \\ \hline
    $A$ & (1), (2), (3) single $1$-admissible & (1), (2) single removable $1$ \\
    &  & (3) double removable $1$ \\
    $C$ & (1), (2) normal  & (1), (2), (3) normal \\
     & (3) single removable $1$ & \\
    $E$ & (1) single removable $1$ or normal & (1) same as in $Y$ \\
    & (2) single $1$-admissible & (2) single $1$-admissible \\
    & (3) normal & (3) normal \\
    \hline
  \end{tabular}
\end{table}
\begin{table}[htb]
  \begin{tabular}{|c|c|c|} \hline
    slot or block & in $Y$ & in $Y'$ \\ \hline
    $B$ & removable $2$ (resp. normal) & normal \\
    $D$ & normal & normal (resp. $2$-admissible) \\
    \hline
  \end{tabular}
\end{table}
Thus, 
\begin{eqnarray*}
L^1_{s,1,\iota}(Y')-L^1_{s,1,\iota}(Y)&=& 
-x_{s+P^1(l)+i,1}-x_{s+P^1(l)+i+1,1}
+x_{s+P^1(l+1)+i,2}\\
&=& 
-x_{s+P^1(l)+i,1}-x_{s+P^1(l)+i+1,1}
+x_{s+P^1(l)+p_{2,1}+i,2}\\
&=&-\beta_{s+P^1(l)+i,1}.
\end{eqnarray*}

\qed

\subsection{Type ${\rm C}_{n-1}^{(1)}$-case (action on extended Young diagrams)}

In this subsection, let $\mathfrak{g}$ be of type ${\rm C}_{n-1}^{(1)}$ and $k\in I\setminus\{1,n\}$.

\begin{prop}\label{D2closed}
Let $T=(y_t)_{t\in\mathbb{Z}}$ be a sequence in ${\rm REYD}_{{\rm D}^{(2)},k}$ of Definition \ref{DEYD} and $i\in\mathbb{Z}$.
\begin{enumerate}
\item We suppose that the point $(i,y_i)$ is single or double admissible and let $T'=(y_t')_{t\in\mathbb{Z}}$ be the sequence in ${\rm REYD}_{{\rm D}^{(2)},k}$
such that $y_{i}'=y_i-1$ and $y_t'=y_t$ $(t\neq i)$. Then for $s\in \mathbb{Z}_{\geq1}$, putting $j:=y_i$, it holds
\[
L^2_{s,k,\iota}(T')-L^2_{s,k,\iota}(T)=-\beta_{s+P^k(i+k)+[i]_-+k-j,\pi_2(i+k)}.
\]
\item We suppose that the point $(i,y_{i-1})$ is single or double removable and let $T''=(y_t'')_{t\in\mathbb{Z}}$ be the sequence in ${\rm REYD}_{{\rm D}^{(2)},k}$
such that $y_{i-1}''=y_{i-1}+1$ and $y_t''=y_t$ $(t\neq i-1)$. Then for $s\in \mathbb{Z}_{\geq1}$, putting $j:=y_{i-1}''=y_{i-1}+1$, it holds
\[
L^2_{s,k,\iota}(T'')-L^2_{s,k,\iota}(T)=\beta_{s+P^k(i+k-1)+[i-1]_-+k-j,\pi_2(i+k-1)}.
\]
\end{enumerate}
\end{prop}

\nd
{\it Proof.} The claim (ii) follows from (i) just as in the proof of Proposition \ref{A2closed}.
In the case $1\leq \pi_2(i+k)\leq n-2$, we can prove (i) by a similar argument to Lemma \ref{lemA2-1}-\ref{lemA2-4}.
If $\pi_2(i+k)=n-1$ then considering the case $\pi_2(i+k+1)=n$, $\pi_2(i+k-1)=n-2$ and the case $\pi_2(i+k+1)=n-2$, $\pi_2(i+k-1)=n$,
one can prove (i) by
replacing $1$, $2$, $3$-admissible and removable points with $n$, $n-1$, $n-2$-admissible and removable points in the proof of
Lemma \ref{lemA2-2} respectively.
In the case $i+k\equiv n+1$ (mod $2n$) so that $\pi_2(i+k)=n$, by replacing $1$, $2$ 
-admissible and removable points with $n$, $n-1$-admissible and removable points 
in the proof of
Lemma \ref{lemA2-3} respectively, we can prove (i).
In the case $i+k\equiv n$ (mod $2n$) so that $\pi_2(i+k)=n$, one can similarly show (i) just as in Lemma \ref{lemA2-4}. \qed

\vspace{3mm}

We remark that in the above proof,
when $n=3$ and $\pi_2(i+k)=n-1=2$, the tables of admissibility/removability slightly become complicated
since one need to consider the cases $n-2(=1)$-admissible/removable points are double or single.
Dividing into cases properly, however, one can prove our claim even if $n=3$ by a similar way to
 Lemma \ref{lemA2-2}.

\subsection{Type ${\rm C}_{n-1}^{(1)}$-case (action on Young walls)}

In this subsection, let $\mathfrak{g}$ be of type ${\rm C}_{n-1}^{(1)}$ and $k\in\{1,n\}$.

\begin{prop}\label{prop-closednessCw-YW}
\begin{enumerate}
\item
Let $t\in I\setminus\{1,n\}$ and
we suppose that $Y\in {\rm YW}_{{\rm D}^{(2)},k}$ has a $t$-admissible slot
\[
\begin{xy}
(-8,15) *{(-i-1,l+1)}="0000",
(17,15) *{(-i,l+1)}="000",
(15,-3) *{(-i,l)}="00",
(-7,-3) *{(-i-1,l)}="0",
(0,0) *{}="1",
(12,0)*{}="2",
(12,12)*{}="3",
(0,12)*{}="4",
\ar@{--} "1";"2"^{}
\ar@{--} "1";"4"^{}
\ar@{--} "2";"3"^{}
\ar@{--} "3";"4"^{}
\end{xy}
\]
Let $Y'\in{\rm YW}_{{\rm D}^{(2)},k}$ be the Young wall obtained from $Y$
by adding the $t$-block to the slot.
Then for $s\in\mathbb{Z}_{\geq1}$ it follows
\[
L^2_{s,k,\iota}(Y')=L^2_{s,k,\iota}(Y)-\beta_{s+P^k(l)+i,t}.
\]
\item Let $t$ be $t=1$ or $t=n$. We suppose that $Y\in {\rm YW}_{{\rm D}^{(2)},k}$ has a $t$-admissible slot
\[
\begin{xy}
(-8,9) *{(-i-1,l+\frac{1}{2})}="0000",
(17,9) *{(-i,l+\frac{1}{2})}="000",
(15,-3) *{(-i,l)}="00",
(-7,-3) *{(-i-1,l)}="0",
(0,0) *{}="1",
(12,0)*{}="2",
(12,6)*{}="3",
(0,6)*{}="4",
(38,3)*{{\rm or}}="or",
(60,9) *{(-i-1,l+1)}="a0000",
(85,9) *{(-i,l+1)}="a000",
(83,-3) *{(-i,l+\frac{1}{2})}="a00",
(61,-3) *{(-i-1,l+\frac{1}{2})}="a0",
(68,0) *{}="a1",
(80,0)*{}="a2",
(80,6)*{}="a3",
(68,6)*{}="a4",
\ar@{--} "1";"2"^{}
\ar@{--} "1";"4"^{}
\ar@{--} "2";"3"^{}
\ar@{--} "3";"4"^{}
\ar@{--} "a1";"a2"^{}
\ar@{--} "a1";"a4"^{}
\ar@{--} "a2";"a3"^{}
\ar@{--} "a3";"a4"^{}
\end{xy}
\]
Let $Y'\in{\rm YW}_{{\rm D}^{(2)},k}$ be the Young wall obtained from $Y$
by adding the $t$-block to the slot.
Then for $s\in\mathbb{Z}_{\geq1}$ it follows
\[
L^2_{s,k,\iota}(Y')=L^2_{s,k,\iota}(Y)-\beta_{s+P^k(l)+i,t}.
\]
\end{enumerate}
\end{prop}

\nd
{\it Proof.} (i) and (ii) except for $t\in\{n-1,n\}$ can be proved by the same way as in  
the proof of Proposition \ref{prop-closednessAw-YW}.
As for the case $t=n-1$, replacing $L^1_{s,1,\iota}$, $P^1$, and $1$, $2$, $3\in I$ in the proof 
of Proposition \ref{prop-closednessAw-YW} (i) of the case $t=2$ by
$L^2_{s,k,\iota}$, $P^k$ and $n$, $n-1$, $n-2\in I$, respectively, one can similarly show our claim.
As for the case $t=n$, replacing $L^1_{s,1,\iota}$, $P^1$ and $1$, $2\in I$ in the proof 
of Proposition \ref{prop-closednessAw-YW} (ii) by
$L^2_{s,k,\iota}$, $P^k$ and $n$, $n-1\in I$, respectively, one can also similarly prove our claim. \qed

\section{Proof}

\subsection{Type ${\rm A}_{n-1}^{(1)}$-case and ${\rm D}_{n}^{(2)}$-case}

In this subsection, let $\mathfrak{g}$ be of type ${\rm A}_{n-1}^{(1)}$ or 
${\rm D}_{n}^{(2)}$.

\nd
{\it Proof of Theorem \ref{thmA1} and \ref{thmD2}.} 
For each $s\in\mathbb{Z}_{\geq1}$ and $k\in I$, we put
\[
\Xi_{s,k,\io}' :=  \{S_{j_l}'\cd S_{j_2}'S_{j_1}'x_{s,k}\,|\,
l\geq0,j_1,\cd,j_l\geq1\},
\]
\[
l_{s,k,\iota}:=
\begin{cases}
\ovl{L}_{s,k,\iota} & {\rm if}\ \mathfrak{g}\ \text{is of type }{\rm A}_{n-1}^{(1)},\\
L_{s,k,\iota} & {\rm if}\ \mathfrak{g}\ \text{is of type }{\rm D}_{n}^{(2)}
\end{cases}
\]
and let us prove $\Xi'_{s,k,\io}=l_{s,k,\iota}({\rm EYD}_k)$. First, we show $\Xi'_{s,k,\io}\subset l_{s,k,\iota}({\rm EYD}_k)$.
The extended diagram $\phi=(k,k,k,k,\cdots)\in {\rm EYD}_k$ is described as follows:
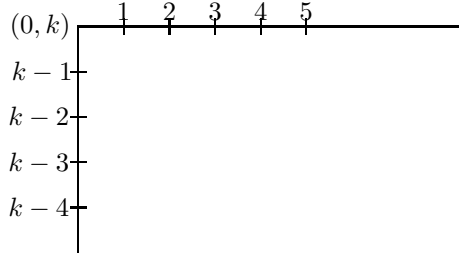
\[
\begin{xy}
(0,0) *{}="1",
(50,0)*{}="2",
(0,-30)*{}="3",
(-5,0)*{(0,k)}="4",
(6,2) *{1}="10",
(6,1) *{}="1010a",
(6,-1) *{}="1010",
(12,2) *{2}="11",
(12,1) *{}="1111a",
(12,-1) *{}="1111",
(18,2) *{3}="12",
(18,-1) *{}="1212",
(18,1) *{}="1212a",
(24,2) *{4}="13",
(24,-1) *{}="1313",
(24,1) *{}="1313a",
(30,2) *{5}="14",
(30,-1) *{}="1414",
(30,1) *{}="1414a",
(-3,-6)*{k-1\ \ \ }="5",
(-1,-6)*{}="5a",
(1,-6)*{}="5aa",
(-1,-12)*{}="6a",
(1,-12)*{}="6aa",
(-1,-18)*{}="7a",
(1,-18)*{}="7aa",
(-1,-24)*{}="8a",
(1,-24)*{}="8aa",
(1,-6)*{}="55",
(-4,-12)*{k-2\ \ }="6",
(1,-12)*{}="66",
(-4,-18)*{k-3\ \ }="7",
(1,-18)*{}="77",
(-4,-24)*{k-4\ \ }="8",
\ar@{-} "1";"2"^{}
\ar@{-} "1";"3"^{}
\ar@{-} "5a";"5aa"^{}
\ar@{-} "6a";"6aa"^{}
\ar@{-} "7a";"7aa"^{}
\ar@{-} "8a";"8aa"^{}
\ar@{-} "1010";"1010a"^{}
\ar@{-} "1111";"1111a"^{}
\ar@{-} "1212";"1212a"^{}
\ar@{-} "1313";"1313a"^{}
\ar@{-} "1414";"1414a"^{}
\end{xy}
\]
Thus, it has one concave corner $(0,k)$ and no convex corner, which implies by (\ref{ovlL2}) and (\ref{LL2}) that
\[
l_{s,k,\iota}(\phi)=x_{s,k}\in l_{s,k,\iota}({\rm EYD}_k).
\]
To prove the inclusion $\Xi'_{s,k,\io}\subset l_{s,k,\iota}({\rm EYD}_k)$, 
we need to show $l_{s,k,\iota}({\rm EYD}_k)$ is closed under the action of $S'_{t,d}$ for any $(t,d)\in\mathbb{Z}_{\geq1}\times I$.
For any $T\in {\rm EYD}_k$ and $(t,d)\in\mathbb{Z}_{\geq1}\times I$, if $x_{t,d}$ has a positive coefficient in $l_{s,k,\iota}(T)$
then there is a concave corner $(i,j)$ in $T$ and $l_{s,k,\iota}(i,j)=x_{t,d}$. 
We get
\[
t=s+P^k(i+j)+{\rm min}\{k-j,i\},\quad d=\ovl{i+j}\ (\mathfrak{g}:{\rm A}_{n-1}^{(1)}),\ d=\pi(i+j)\ (\mathfrak{g}:{\rm D}_{n}^{(2)}). 
\]
Let $T'\in{\rm EYD}_k$
be the extended Young diagram obtained from $T$ by replacing
the concave corner $(i,j)$ by a convex corner, which is the replacement in (\ref{change}).
It follows from (\ref{Sk}) and
Proposition \ref{prop-closednessAD} that 
\[
l_{s,k,\iota}(T')=l_{s,k,\iota}(T)-\beta_{t,d}=S'_{t,d}l_{s,k,\iota}(T).
\]
If $x_{t,d}$ has a negative coefficient in $l_{s,k,\iota}(T)$ then
there is a convex corner $P=(p_1,p_2)$ in $T$ such that $l_{s,k,\iota}(P)=x_{t,d}$.
Since $P$ is a convex corner, it holds $p_1>0$ and $p_2<k$ so that
there exist $i\in\mathbb{Z}_{\geq0}$ and $j\in\mathbb{Z}_{\leq k}$ such that
$p_1=i+1$, $p_2=j-1$. By $l_{s,k,\iota}(P)=l_{s,k,\iota}(i+1,j-1)=x_{t,d}$, we have
\[
t=s+P^k(i+j)+{\rm min}\{k-j,i\}+1,\quad d=\ovl{i+j}\ (\mathfrak{g}:{\rm A}_{n-1}^{(1)}),\ d=\pi(i+j)\ (\mathfrak{g}:{\rm D}_{n}^{(2)}). 
\]
By $i\in\mathbb{Z}_{\geq0}$, $j\in\mathbb{Z}_{\leq k}$ and $s\geq1$,
we obtain
\begin{equation}\label{pr-3}
t-1=s+P^k(i+j)+{\rm min}\{k-j,i\}\geq1.
\end{equation}
Let $T'\in{\rm EYD}_k$
be the extended Young diagram obtained from $T$ by replacing
the convex corner $(i+1,j-1)$ by a concave corner, which is the opposite replacement in (\ref{change}).
Hence, $T'$ has the concave corner $(i,j)$ and 
Proposition \ref{prop-closednessAD} says the following:
\begin{equation}\label{pr-2}
l_{s,k,\iota}(T)=l_{s,k,\iota}(T')-\beta_{t-1,d}.
\end{equation}
Considering (\ref{Sk}), it follows
\[
l_{s,k,\iota}(T')=l_{s,k,\iota}(T)+\beta_{t-1,d}=
S'_{t,d}l_{s,k,\iota}(T).
\]
Hence, we proved the closedness of $l_{s,k,\iota}({\rm EYD}_k)$, that is,
if $T\in {\rm EYD}_k$ then $S'_{t,d}l_{s,k,\iota}(T)\in l_{s,k,\iota}({\rm EYD}_k)$ for any $(t,d)\in\mathbb{Z}_{\geq1}\times I$, which implies
the inclusion $\Xi'_{s,k,\io}\subset l_{s,k,\iota}({\rm EYD}_k)$. 

Note that if $(i,j)$ is a convex corner in an extended Young diagram $T$ then $i\geq1$ and $k>j$. 
The non-negativity of $P^k(i+j)$ and positivity of $s$ imply
\[
s+P^k(i+j)+{\rm min}\{k-j,i\}\geq2,
\]
where $s+P^k(i+j)+{\rm min}\{k-j,i\}$ is the left index in (\ref{ovlL1}) and (\ref{LL1}).
By (\ref{ovlL1}), (\ref{ovlL2}), (\ref{LL1}), (\ref{LL2}), $\Xi'_{s,k,\io}\subset l_{s,k,\iota}({\rm EYD}_k)$ and
$\Xi'_{\iota}=\bigcup_{(s,k)\in\mathbb{Z}_{\geq1}\times I}\Xi'_{s,k,\io}$, the sequence $\iota$
satisfies the $\Xi'$-positivity condition. 

Next, we prove $\Xi'_{s,k,\io}\supset l_{s,k,\iota}({\rm EYD}_k)$.
We identify each $T\in {\rm EYD}_k$ as a Young diagram consisting of several boxes, where each box is the square
whose length of sides are $1$. For $T\in {\rm EYD}_k$, we show $l_{s,k,\iota}(T)\in\Xi'_{s,k,\io}$ by induction on the 
number of boxes in $T$. 
In the case the number of box is $0$ so that $T=\phi$, 
it holds $l_{s,k,\iota}(\phi)=x_{s,k}\in\Xi'_{s,k,\io}$. So we assume $T$ has at least one box.
Note that $T$ is obtained from $\phi$ by a sequence of replacements in (\ref{change}).
Considering Proposition \ref{prop-closednessAD}, one can write $l_{s,k,\iota}(T)$ as
\[
l_{s,k,\iota}(T)=x_{s,k}-\sum_{(t,d)\in\mathbb{Z}_{\geq1}\times I} c_{t,d}\beta_{t,d}
\]
with some non-negative integers $\{c_{t,d}\}$ such that $c_{t,d}=0$ except for finitely many $(t,d)$.
For each corner $(i,j)$ in $T$, it holds $s+P^k(i+j)+{\rm min}\{k-j,i\}\geq s$ so that if $c_{t,d}\neq0$ then $t\geq s$.
Since $T\neq\phi$, there exists $(t',d')\in\mathbb{Z}_{\geq1}\times I$ such that
\[
(t',d')={\rm max}\{(t,d)\in\mathbb{Z}_{\geq1}\times I | c_{t,d}\neq0 \}, 
\]
where the order on $\mathbb{Z}_{\geq1}\times I$ is defined as in the subsection \ref{seno}.
It follows $t'\geq s$ by $c_{t',d'}\neq0$.
The definition (\ref{betak}) of $\beta_{t,d}$ and Definition \ref{adapt} mean
the coefficient of $x_{t'+1,d'}$ is negative in $l_{s,k,\iota}(T)$. We define $(t'',d'')$ as
\[
(t'',d'')={\rm min}\{(t,d)\in\mathbb{Z}_{\geq1}\times I | \text{the coefficient of }x_{t,d}\ \text{in } l_{s,k,\iota}(T)\text{ is negative} \}.
\]
We see that $T$ has a convex corner $(i+1,j-1)$ such that $i\in\mathbb{Z}_{\geq0}$, $j\in\mathbb{Z}_{\leq k}$ and $ l_{s,k,\iota}(i+1,j-1)=x_{t'',d''}$. 
By the same way as in (\ref{pr-3}), it holds $t''-1\geq1$.
Let $T'\in {\rm EYD}_k$ be
the extended Young diagram obtained from $T$ by replacing the
convex corner $(i+1,j-1)$ by a concave corner. 
Just as in (\ref{pr-2}), we obtain
\[
l_{s,k,\iota}(T)+\beta_{t''-1,d''}=l_{s,k,\iota}(T').
\]
The minimality of $(t'',d'')$ means the coefficient of $x_{t''-1,d''}$ in 
$l_{s,k,\iota}(T')$ is positive and
\[
S'_{t''-1,d''}l_{s,k,\iota}(T')=l_{s,k,\iota}(T')-\beta_{t''-1,d''}=l_{s,k,\iota}(T).
\]
Note that the number of boxes in $T'$ is smaller than $T$. 
By the induction assumption, 
it holds $l_{s,k,\iota}(T')\in\Xi'_{s,k,\io}$ so that $l_{s,k,\iota}(T')=S_{j_l}'\cd S_{j_2}'S_{j_1}'x_{s,k}$
with some $l\in\mathbb{Z}_{\geq0}$ and $j_1,\cdots,j_l\in\mathbb{Z}_{\geq1}$. Therefore,
\[
l_{s,k,\iota}(T)=S'_{t''-1,d''}l_{s,k,\iota}(T')=S'_{t''-1,d''}S_{j_l}'\cd S_{j_2}'S_{j_1}'x_{s,k}\in \Xi'_{s,k,\io},
\]
which yields $\Xi'_{s,k,\io}= l_{s,k,\iota}({\rm EYD}_k)$.
Theorem \ref{thmA1} and \ref{thmD2}
follow by Theorem \ref{polyhthm} and $\Xi'_{\iota}=\bigcup_{(s,k)\in\mathbb{Z}_{\geq1}\times I}\Xi'_{s,k,\io}$. \qed

\subsection{Type ${\rm A}_{2n-2}^{(2)}$-case and ${\rm C}_{n-1}^{(1)}$-case}

In this subsection, let $\mathfrak{g}$ be of type ${\rm A}_{2n-2}^{(2)}$ or ${\rm C}_{n-1}^{(1)}$.

\nd
{\it Proof of Theorem \ref{thmA2} and \ref{thmC1}.} 
For each $s\in\mathbb{Z}_{\geq1}$ and $k\in I$, we put
\[
\Xi_{s,k,\io}' :=  \{S_{j_l}'\cd S_{j_2}'S_{j_1}'x_{s,k}\,|\,
l\geq0,j_1,\cd,j_l\geq1\}.
\]
Let us prove 
\[
\Xi'_{s,k,\io}=L^1_{s,k,\iota}({\rm REYD}_{{\rm A}^{(2)},k})\ \ ({\rm for}\ k\in I\setminus\{1\}),\quad
\Xi'_{s,1,\io}=L^1_{s,1,\iota}({\rm YW}_{{\rm A}^{(2)},1}),
\] 
if $\mathfrak{g}$ is of type ${\rm A}_{2n-2}^{(2)}$ and
\[
\Xi'_{s,k,\io}=L^2_{s,k,\iota}({\rm REYD}_{{\rm D}^{(2)},k})\ \ ({\rm for}\ k\in I\setminus\{1,n\}),\quad
\Xi'_{s,1,\io}=L^2_{s,1,\iota}({\rm YW}_{{\rm D}^{(2)},1}),\quad
\Xi'_{s,n,\io}=L^2_{s,n,\iota}({\rm YW}_{{\rm D}^{(2)},n})
\] 
if $\mathfrak{g}$ is of type ${\rm C}_{n-1}^{(1)}$.

\vspace{2mm}

\nd
\underline{Proof of $\Xi'_{s,k,\io}=L^1_{s,k,\iota}({\rm REYD}_{{\rm A}^{(2)},k})$,\ $\Xi'_{s,k,\io}=L^2_{s,k,\iota}({\rm REYD}_{{\rm D}^{(2)},k})$}

\vspace{2mm}

First, 
taking $k\in I\setminus\{1\}$ if $\mathfrak{g}$ is of type ${\rm A}_{2n-2}^{(2)}$, 
$k\in I\setminus\{1,n\}$ if $\mathfrak{g}$ is of type ${\rm C}_{n-1}^{(1)}$ and putting
\[
l_{s,k,\iota}:=
\begin{cases}
L^1_{s,k,\iota} & {\rm if}\ \mathfrak{g}\ \text{is of type }{\rm A}_{2n-2}^{(2)},\\
L^2_{s,k,\iota} & {\rm if}\ \mathfrak{g}\ \text{is of type }{\rm C}_{n-1}^{(1)},
\end{cases},\quad
{\rm REYD}_k:=
\begin{cases}
{\rm REYD}_{{\rm A}^{(2)},k} & {\rm if}\ \mathfrak{g}\ \text{is of type }{\rm A}_{2n-2}^{(2)},\\
{\rm REYD}_{{\rm D}^{(2)},k} & {\rm if}\ \mathfrak{g}\ \text{is of type }{\rm C}_{n-1}^{(1)},
\end{cases}
\]
we show
$\Xi'_{s,k,\io}=l_{s,k,\iota}({\rm REYD}_{k})$. We also set $l_{s,k,{\rm ad}}:=L^1_{s,k,{\rm ad}}$ (resp. $L^2_{s,k,{\rm ad}}$)
and $l_{s,k,{\rm re}}:=L^1_{s,k,{\rm re}}$ (resp. $L^2_{s,k,{\rm re}}$) if $\mathfrak{g}$ is of type ${\rm A}_{2n-2}^{(2)}$ (resp. ${\rm C}_{n-1}^{(1)}$). 
One takes
$\phi:=(\phi_l)_{l\in\mathbb{Z}}\in{\rm REYD}_k$ as $\phi_l=k+l$ for $l\in\mathbb{Z}_{<0}$ and $\phi_l=k$ for $l\in\mathbb{Z}_{\geq0}$. Then, it is described as
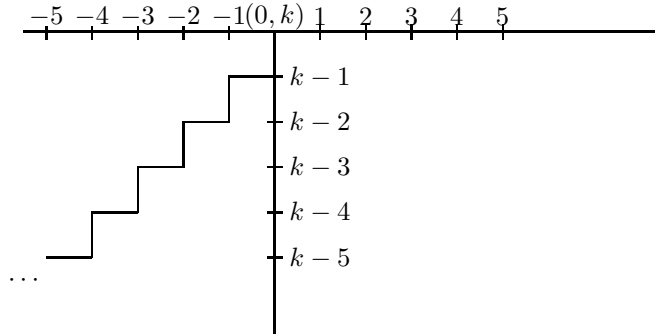
\[
\begin{xy}
(-33,0) *{}="-6",
(-6,2) *{-1}="-1",
(-12,2) *{-2}="-2",
(-18,2) *{-3}="-3",
(-24,2) *{-4}="-4",
(-30,2) *{-5}="-5",
(-6,-1) *{}="-1a",
(-12,-1) *{}="-2a",
(-18,-1) *{}="-3a",
(-24,-1) *{}="-4a",
(-30,-1) *{}="-5a",
(0,0) *{}="1",
(50,0)*{}="2",
(0,-40)*{}="3",
(0,2)*{(0,k)}="4",
(6,2) *{1}="10",
(6,-1) *{}="1010",
(12,2) *{2}="11",
(12,-1) *{}="1111",
(18,2) *{3}="12",
(18,-1) *{}="1212",
(24,2) *{4}="13",
(24,-1) *{}="1313",
(30,2) *{5}="14",
(30,-1) *{}="1414",
(6,-6)*{k-1}="k-1",
(6,-12)*{k-2}="k-2",
(6,-18)*{k-3}="k-3",
(6,-24)*{k-4}="k-4",
(6,-30)*{k-5}="k-5",
(-1,-6)*{}="5",
(1,-6)*{}="55",
(-1,-12)*{}="6",
(-6,-12)*{}="6a",
(-12,-12)*{}="6aa",
(-12,-18)*{}="6aaa",
(-18,-18)*{}="6aaaa",
(-18,-24)*{}="st1",
(-24,-24)*{}="st2",
(-24,-30)*{}="st3",
(-30,-30)*{}="st4",
(-33,-33)*{\cdots}="stdot",
(1,-12)*{}="66",
(-1,-18)*{}="t",
(1,-18)*{}="tt",
(-6,-6)*{}="7",
(1,-18)*{}="77",
(-1,-24)*{}="8",
(1,-24)*{}="88",
(-1,-30)*{}="9",
(1,-30)*{}="99",
\ar@{-} "1";"-6"^{}
\ar@{-} "1";"2"^{}
\ar@{-} "1";"3"^{}
\ar@{-} "-1";"-1a"^{}
\ar@{-} "-2";"-2a"^{}
\ar@{-} "-3";"-3a"^{}
\ar@{-} "-4";"-4a"^{}
\ar@{-} "-5";"-5a"^{}
\ar@{-} "5";"55"^{}
\ar@{-} "6";"66"^{}
\ar@{-} "t";"tt"^{}
\ar@{-} "8";"88"^{}
\ar@{-} "7";"5"^{}
\ar@{-} "7";"6a"^{}
\ar@{-} "6aa";"6a"^{}
\ar@{-} "6aaa";"6aa"^{}
\ar@{-} "6aaaa";"6aaa"^{}
\ar@{-} "st1";"6aaaa"^{}
\ar@{-} "st1";"st2"^{}
\ar@{-} "st2";"st3"^{}
\ar@{-} "st3";"st4"^{}
\ar@{-} "9";"99"^{}
\ar@{-} "10";"1010"^{}
\ar@{-} "11";"1111"^{}
\ar@{-} "12";"1212"^{}
\ar@{-} "13";"1313"^{}
\ar@{-} "14";"1414"^{}
\end{xy}
\]
All points are neither admissible nor removable
except for the $k$-admissible point $(0,k)$, which implies $l_{s,k,\iota}(\phi)=x_{s,k}$. Thus, it holds $x_{s,k}\in l_{s,k,\iota}({\rm REYD}_{k})$.
Let us show the inclusion $\Xi'_{s,k,\io}\subset l_{s,k,\iota}({\rm REYD}_{k})$.
We need to show $l_{s,k,\iota}({\rm REYD}_{k})$ is closed under the action of $S'_{t,d}$ for any $(t,d)\in\mathbb{Z}_{\geq1}\times I$.
Let $T=(y_m)_{m\in\mathbb{Z}}\in {\rm REYD}_{k}$ and $(t,d)\in\mathbb{Z}_{\geq1}\times I$. If $x_{t,d}$ has a positive coefficient in $l_{s,k,\iota}(T)$
then there is an admissible point $(i,y_i)$ in $T$ such that $l_{s,k,{\rm ad}}(i,y_i)=x_{t,d}$ by (\ref{L1kdef}) and (\ref{L2kdef}). Putting
$j:=y_i$, it follows
\[
t=s+P^k(i+k)+[i]_-+k-j,\quad d=\pi_1(i+k)\ (\mathfrak{g}:{\rm A}_{2n-2}^{(2)}),\ d=\pi_2(i+k)\ (\mathfrak{g}:{\rm C}_{n-1}^{(1)}).
\]
Let $T'=(y_m')_{m\in\mathbb{Z}}$ be the sequence such that
$y_i'=y_i-1$ and $y_m'=y_m$ $(m\neq i)$. Since $(i,j)=(i,y_i)$ is admissible, we obtain $T'\in {\rm REYD}_{k}$.
By Proposition \ref{A2closed} and \ref{D2closed} (i), we see that
\[
S'_{t,d}l_{s,k,\iota}(T)=l_{s,k,\iota}(T)-\beta_{t,d}=l_{s,k,\iota}(T')\in l_{s,k,\iota}({\rm REYD}_{k}).
\]
If $x_{t,d}$ has a negative coefficient in $l_{s,k,\iota}(T)$
then there is a removable point $(i,y_{i-1})$ in $T$ such that $l_{s,k,{\rm re}}(i,y_{i-1})=x_{t,d}$ by (\ref{L1kdef}) and (\ref{L2kdef}). 
It holds
\[
t=s+P^k(i+k-1)+[i-1]_-+k-y_{i-1},\quad d=\pi_1(i+k-1)\ (\mathfrak{g}:{\rm A}_{2n-2}^{(2)}),\ d=\pi_2(i+k-1)\ (\mathfrak{g}:{\rm C}_{n-1}^{(1)}).
\]
Let $T''=(y_t'')_{t\in\mathbb{Z}}$ be the sequence
such that $y_{i-1}''=y_{i-1}+1$ and $y_t''=y_t$ $(t\neq i-1)$. Since $(i,y_{i-1})$ is a removable point
we see that $T''\in{\rm REYD}_{k}$.
If $t\geq2$ so that $(t,d)^{(-)}=(t-1,d)\in \mathbb{Z}_{\geq1}\times I$ then Proposition \ref{A2closed} and \ref{D2closed} (ii) yield
\[
S'_{t,d}l_{s,k,\iota}(T)=l_{s,k,\iota}(T)+\beta_{t-1,d}=l_{s,k,\iota}(T')\in l_{s,k,\iota}({\rm REYD}_{k}).
\]
Therefore,
$l_{s,k,\iota}({\rm REYD}_{k})$ is closed under the action of $S'_{j}$ for any $j\in\mathbb{Z}_{\geq1}$.
Hence the inclusion $\Xi'_{s,k,\io}\subset l_{s,k,\iota}({\rm REYD}_{k})$ follows.

Let us prove the converse inclusion.
We identify each $T\in {\rm REYD}_k$ as a pile
of boxes in $\mathbb{R}\times \mathbb{R}_{\leq k}$, where the box is the square whose length of sides is $1$.
If $T$ is obtained from $\phi$ by adding $m$ boxes then
we say the number of boxes in $T$ is $m$. For instance, the number of boxes in the element (\ref{reydA-ex3}) is $13$.
For any $T=(y_l)_{l\in\mathbb{Z}}\in {\rm REYD}_{k}$,
we show $l_{s,k,\iota}(T)\in \Xi'_{s,k,\iota}$ by induction on the number of boxes. 
In the case
the number of boxes is $0$ so that $T=\phi$, it follows $l_{s,k,\iota}(\phi)=x_{s,k}\in\ \Xi'_{s,k,\iota}$.
Thus, we assume $T\neq \phi$.
Setting $m:={\rm min}\{l\in\mathbb{Z} | y_l < k+l \}$, one gets $m\leq 1$.
If $m=1$ then $y_l=k+l$ for $l\in\mathbb{Z}_{\leq0}$, in particular, $y_0=k$ so that $y_r=k$ for $r\in\mathbb{Z}_{>0}$ by $y_0\leq y_r\leq k$,
which implies $T=\phi$. Since we assumed $T\neq \phi$, it holds $m\leq 0$ and $y_m<k$.
Putting
\[
m_1:={\rm min}\{y_l | m\leq l\},
\]
one obtain $m_1<k$ and there exists $i\in\mathbb{Z}_{\geq m}$ such that $y_i=m_1$ and $y_{i+1}>y_i$. Then the point $(i+1,y_i)$ is removable.
Considering (\ref{pos-A2-2}) and (\ref{pos-C1-2}), it holds
\begin{equation}\label{prA2-1}
s+P^k(i+k)+[i]_-+k-y_i-1\geq s.
\end{equation}
Defining 
$T''=(y_l'')_{l\in\mathbb{Z}}\in {\rm REYD}_{k}$ as $y_i''=y_i+1$ and $y_t''=y_t$ $(t\neq i)$, we obtain
\[
l_{s,k,\iota}(T)=l_{s,k,\iota}(T'')-\beta_{t_1,d_1}
\] 
with some $(t_1,d_1)\in\mathbb{Z}_{\geq s}\times I$ by (\ref{prA2-1}), Proposition \ref{A2closed} and \ref{D2closed} (ii). 
Repeating this argument, we see
\[
l_{s,k,\iota}(T)=l_{s,k,\iota}(\phi)-\sum_{(t,d)\in\mathbb{Z}_{\geq s}\times I} c_{t,d}\beta_{t,d}
=x_{s,k}-\sum_{(t,d)\in\mathbb{Z}_{\geq s}\times I} c_{t,d}\beta_{t,d}
\]
with non-negative integers $\{c_{t,d}\}$ such that $c_{t,d}=0$ except for finitely many $(t,d)$.
By $T\neq\phi$, there exists $(t',d')\in\mathbb{Z}_{\geq s}\times I$ such that
\[
(t',d')={\rm max}\{(t,d)\in\mathbb{Z}_{\geq s}\times I | c_{t,d}>0 \}, 
\]
where the order on $\mathbb{Z}_{\geq1}\times I$ is defined as in the subsection \ref{seno}.
Considering the definition (\ref{pij3}) of $\beta_{t,d}$, 
the coefficient of $x_{t'+1,d'}$ is negative in $l_{s,k,\iota}(T)$. Let $(t'',d'')$ be
\[
(t'',d'')={\rm min}\{(t,d)\in\mathbb{Z}_{\geq s}\times I | \text{the coefficient of }x_{t,d}\ \text{in } l_{s,k,\iota}(T)\text{ is negative} \}.
\]
Since the coefficient of $x_{t'',d''}$ is negative, we see that there is a removable point $(\xi,y_{\xi-1})$ in $T$ such that
$l_{s,k,{\rm re}}(\xi,y_{\xi-1})=x_{t'',d''}$ and $t''=s+P^k(\xi+k-1)+[\xi-1]_-+k-y_{\xi-1}$. 
The relations (\ref{pos-A2-2}) and (\ref{pos-C1-2}) imply 
\begin{equation}\label{pos-cond-pr}
t''\geq s+1\geq2.
\end{equation}
Let $T''=(y_t'')_{t\in\mathbb{Z}}\in {\rm REYD}_{k}$ be
the element such that $y_{\xi-1}''=y_{\xi-1}+1$ and $y_{t}''=y_t$ ($t\neq \xi-1$).
Taking Proposition \ref{A2closed} and \ref{D2closed} (ii)
into account, we obtain
\[
l_{s,k,\iota}(T)=l_{s,k,\iota}(T'')-\beta_{t''-1,d''}.
\]
Combining this formula with the minimality
of $(t'',d'')$, 
the coefficient of $x_{t''-1,d''}$ in $l_{s,k,\iota}(T'')$ is positive. Thus,
\begin{equation}\label{prA2-2}
l_{s,k,\iota}(T)=S'_{t''-1,d''}l_{s,k,\iota}(T'').
\end{equation}
Note that the number of boxed in $T''$ is smaller than those of $T$.
Using the induction assumption, it holds $l_{s,k,\iota}(T'')\in \Xi_{s,k,\iota}'$.
In conjunction with (\ref{prA2-2}), the our claim 
$l_{s,k,\iota}(T)\in \Xi_{s,k,\iota}'$ follows.

\vspace{3mm}

\nd
\underline{Proof of $\Xi'_{s,1,\io}=L^1_{s,1,\iota}({\rm YW}_{{\rm A}^{(2)},1})$, $\Xi'_{s,1,\io}=L^2_{s,1,\iota}({\rm YW}_{{\rm D}^{(2)},1})$,
\ $\Xi'_{s,n,\io}=L^2_{s,n,\iota}({\rm YW}_{{\rm D}^{(2)},n})$}

\vspace{3mm}

We take $k=1$ in the case $\mathfrak{g}$ is of type ${\rm A}^{(2)}_{2n-2}$
and $k\in\{1,n\}$ in the case $\mathfrak{g}$ is of type ${\rm C}^{(1)}_{n-1}$. Putting
\[
l_{s,k,\iota}:=
\begin{cases}
L^1_{s,k,\iota} & {\rm if}\ \mathfrak{g}\ \text{is of type }{\rm A}_{2n-2}^{(2)},\\
L^2_{s,k,\iota} & {\rm if}\ \mathfrak{g}\ \text{is of type }{\rm C}_{n-1}^{(1)},
\end{cases},\quad
{\rm YW}_k:=
\begin{cases}
{\rm YW}_{{\rm A}^{(2)},k} & {\rm if}\ \mathfrak{g}\ \text{is of type }{\rm A}_{2n-2}^{(2)},\\
{\rm YW}_{{\rm D}^{(2)},k} & {\rm if}\ \mathfrak{g}\ \text{is of type }{\rm C}_{n-1}^{(1)},
\end{cases}
\]
let us prove $\Xi'_{s,k,\io}=l_{s,k,\iota}({\rm YW}_k)$. We also set $l_{s,k,{\rm ad}}=L^1_{s,k,{\rm ad}}$ (resp. $L^2_{s,k,{\rm ad}}$)
and $l_{s,k,{\rm re}}=L^1_{s,k,{\rm re}}$ (resp. $L^2_{s,k,{\rm re}}$) if $\mathfrak{g}$ is of type ${\rm A}_{2n-2}^{(2)}$ (resp. ${\rm C}_{n-1}^{(1)}$).
The ground state wall $Y_{\Lambda_k}\in {\rm YW}_k$ has neither admissible slot nor removable block
except for the $k$-single admissible slot
\[
\begin{xy}
(-8,8) *{(-1,k+1)}="0000",
(17,8) *{(0,k+1)}="000",
(15,-3) *{(0,k+\frac{1}{2})}="00",
(-7,-3) *{(-1,k+\frac{1}{2})}="0",
(0,0) *{}="1",
(12,0)*{}="2",
(12,6)*{}="3",
(0,6)*{}="4",
\ar@{--} "1";"2"^{}
\ar@{--} "1";"4"^{}
\ar@{--} "2";"3"^{}
\ar@{--} "3";"4"^{}
\end{xy}
\]
Thus, it follows $l_{s,k,\iota}(Y_{\Lambda_k})=x_{s,k}\in \Xi'_{s,k,\io}$.
To prove the inclusion $\Xi'_{s,k,\io}\subset l_{s,k,\iota}({\rm YW}_k)$,
we need to show $l_{s,k,\iota}({\rm YW}_k)$ is closed under the action of $S'_{\xi,d}$ for any $(\xi,d)\in\mathbb{Z}_{\geq1}\times I$.
For any $Y\in{\rm YW}_k$, if the coefficient of $x_{\xi,d}$ is positive in $l_{s,k,\iota}(Y)$ then
there is a $d$-admissible slot
\begin{equation}\label{YW-pr1}
P=
\begin{xy}
(-68,8) *{(-i-1,l+1)}="0000b",
(-36,8) *{(-i,l+1)}="000b",
(-42,-6) *{(-i,l)}="00b",
(-70,-6) *{(-i-1,l)}="0b",
(-8,8) *{(-i-1,l+1)}="0000",
(17,8) *{(-i,l+1)}="000",
(15,-3) *{(-i,l+\frac{1}{2})}="00",
(-7,-3) *{(-i-1,l+\frac{1}{2})}="0",
(40,8) *{(-i-1,l+\frac{1}{2})}="0000a",
(65,8) *{(-i,l+\frac{1}{2})}="000a",
(63,-3) *{(-i,l)}="00a",
(41,-3) *{(-i-1,l)}="0a",
(-24,3) *{\text{or}}="or2",
(-60,-6) *{}="1b",
(-48,-6)*{}="2b",
(-48,6)*{}="3b",
(-60,6)*{}="4b",
(30,3) *{\text{or}}="or1",
(0,0) *{}="1",
(12,0)*{}="2",
(12,6)*{}="3",
(0,6)*{}="4",
(48,0) *{}="1a",
(60,0)*{}="2a",
(60,6)*{}="3a",
(48,6)*{}="4a",
\ar@{--} "1b";"2b"^{}
\ar@{--} "1b";"4b"^{}
\ar@{--} "2b";"3b"^{}
\ar@{--} "3b";"4b"^{}
\ar@{--} "1";"2"^{}
\ar@{--} "1";"4"^{}
\ar@{--} "2";"3"^{}
\ar@{--} "3";"4"^{}
\ar@{--} "1a";"2a"^{}
\ar@{--} "1a";"4a"^{}
\ar@{--} "2a";"3a"^{}
\ar@{--} "3a";"4a"^{}
\end{xy}
\end{equation}
such that $l_{s,k,{\rm ad}}(P)=x_{\xi,d}$. Hence, it follows
$\xi=s+P^k(l)+i$.
Let $Y'\in {\rm YW}_k$ be the
proper Young wall obtained from $Y$ by adding the $d$-block to the slot $P$.
By Proposition \ref{prop-closednessAw-YW}, \ref{prop-closednessCw-YW}, it holds
\[
S'_{\xi,d}l_{s,k,\iota}(Y)=l_{s,k,\iota}(Y)-\beta_{s+P^k(l)+i,d}=l_{s,k,\iota}(Y')\in l_{s,k,\iota}({\rm YW}_k).
\]
If the coefficient of $x_{\xi,d}$ is negative in $l_{s,k,\iota}(Y)$ then there is a removable $d$-block
\[
B=
\begin{xy}
(-68,8) *{(-i-1,l+1)}="0000b",
(-36,8) *{(-i,l+1)}="000b",
(-42,-6) *{(-i,l)}="00b",
(-70,-6) *{(-i-1,l)}="0b",
(-8,8) *{(-i-1,l+1)}="0000",
(17,8) *{(-i,l+1)}="000",
(15,-3) *{(-i,l+\frac{1}{2})}="00",
(-7,-3) *{(-i-1,l+\frac{1}{2})}="0",
(40,8) *{(-i-1,l+\frac{1}{2})}="0000a",
(65,8) *{(-i,l+\frac{1}{2})}="000a",
(63,-3) *{(-i,l)}="00a",
(41,-3) *{(-i-1,l)}="0a",
(-24,3) *{\text{or}}="or2",
(-60,-6) *{}="1b",
(-48,-6)*{}="2b",
(-48,6)*{}="3b",
(-60,6)*{}="4b",
(30,3) *{\text{or}}="or1",
(0,0) *{}="1",
(12,0)*{}="2",
(12,6)*{}="3",
(0,6)*{}="4",
(48,0) *{}="1a",
(60,0)*{}="2a",
(60,6)*{}="3a",
(48,6)*{}="4a",
\ar@{-} "1b";"2b"^{}
\ar@{-} "1b";"4b"^{}
\ar@{-} "2b";"3b"^{}
\ar@{-} "3b";"4b"^{}
\ar@{-} "1";"2"^{}
\ar@{-} "1";"4"^{}
\ar@{-} "2";"3"^{}
\ar@{-} "3";"4"^{}
\ar@{-} "1a";"2a"^{}
\ar@{-} "1a";"4a"^{}
\ar@{-} "2a";"3a"^{}
\ar@{-} "3a";"4a"^{}
\end{xy}
\]
such that $l_{s,k,{\rm re}}(B)=x_{\xi,d}$, which implies $\xi=s+P^k(l)+i+1$.
Let $Y''\in {\rm YW}_k$ be the
proper Young wall obtained from $Y$ by removing the $d$-block.
Thus, $Y''$ has the $d$-admissible slot $P$ like as (\ref{YW-pr1}).
By Proposition \ref{prop-closednessAw-YW}, \ref{prop-closednessCw-YW}, it holds
$l_{s,k,\iota}(Y)=l_{s,k,\iota}(Y'')-\beta_{s+P^k(l)+i,d}$
so that 
\begin{equation}\label{YW-pr2}
S'_{\xi,d}l_{s,k,\iota}(Y)=l_{s,k,\iota}(Y)+\beta_{\xi-1,d}
=l_{s,k,\iota}(Y)+\beta_{s+P^k(l)+i,d}=l_{s,k,\iota}(Y'')\in l_{s,k,\iota}({\rm YW}_k).
\end{equation}
Therefore, the set $l_{s,k,\iota}({\rm YW}_k)$ is closed under the action of $S'$ and
$\Xi'_{s,k,\io}\subset l_{s,k,\iota}({\rm YW}_k)$.

Let us prove $l_{s,k,\iota}({\rm YW}_k)\subset \Xi'_{s,k,\io}$.
When we get a proper Young wall $Y$ by adding $m$ blocks to $Y_{\Lambda_k}$,
we say the number of blocks in $Y$ is $m$.
For any $Y\in {\rm YW}_k$, we prove $l_{s,k,\iota}(Y)\in \Xi'_{s,k,\iota}$
using the induction on the number of blocks in $Y$.
If the number of blocks is $0$ then it holds $Y=Y_{\Lambda_k}$ so that
$l_{s,k,\iota}(Y)=x_{s,k}\in \Xi'_{s,k,\iota}$. Hence, we may assume $Y\neq Y_{\Lambda_k}$.
Using Proposition \ref{prop-closednessAw-YW}, \ref{prop-closednessCw-YW},
one can describe $l_{s,k,\iota}(Y)$ as
\[
l_{s,k,\iota}(Y)=x_{s,k}-\sum_{(t,d)\in\mathbb{Z}_{\geq s}\times I} c_{t,d}\beta_{t,d}
\]
with non-negative integers $\{c_{t,d}\}$. Except for finitely many $(t,d)$, it holds $c_{t,d}=0$.
It follows by $Y\neq Y_{\Lambda_k}$ that one can take $(t',d')\in\mathbb{Z}_{\geq s}\times I$ as
$(t',d')
={\rm max}\{(t,d)\in\mathbb{Z}_{\geq s}\times I |c_{t,d}>0 \}$.
Considering the definition (\ref{pij3}) of $\beta_{t',d'}$, one can verify 
the coefficient of $x_{t'+1,d'}$ is negative in $l_{s,k,\iota}(Y)$. Thus,
we can take $(t'',d'')$ as
\[
(t'',d'')={\rm min}\{(t,d)\in\mathbb{Z}_{\geq 1}\times I | \text{the coefficient of }x_{t,d}\text{ in }l_{s,k,\iota}(Y)\text{ is negative}\}.
\]
There exists a removable block $B$ such that $l_{s,k,{\rm re}}(B)=x_{t'',d''}$ and $t''\geq s+1$ by (\ref{A2YW-pr}), (\ref{C1YW-pr}).
Let $Y''\in{\rm YW}_k$ be the Young wall which is obtained from $Y$ by removing the block $B$. 
Just as in (\ref{YW-pr2}),
the relation between $l_{s,k,\iota}(Y)$ and $l_{s,k,\iota}(Y'')$ is as follows:
\[
l_{s,k,\iota}(Y)=l_{s,k,\iota}(Y'')-\beta_{t''-1,d''}.
\]
This equation and the minimality of $(t'',d'')$ imply the coefficient of $x_{t''-1,d''}$ is positive in 
$l_{s,k,\iota}(Y'')$ so that
\begin{equation}\label{A2YW-pr2}
l_{s,k,\iota}(Y)=S'_{t''-1,d''}l_{s,k,\iota}(Y'').
\end{equation}
Note that the number of boxes in $Y''$ is smaller than $Y$. By induction assumption,
we obtain $l_{s,k,\iota}(Y'')\in \Xi'_{s,k,\iota}$, which yields $l_{s,k,\iota}(Y)\in \Xi'_{s,k,\iota}$ by (\ref{A2YW-pr2}).
Therefore, we get the inclusion $\Xi'_{s,k,\io}\supset l_{s,k,\iota}({\rm YW}_k)$.

\vspace{2mm}

\nd
\underline{Proof of $\Xi'$-positivity condition}

\vspace{2mm}

By the above argument and the definitions of $L_{s,k,\iota}^1$, $L_{s,k,\iota}^2$ in (\ref{L1kdef}), (\ref{L11-def}),
(\ref{L2kdef}), (\ref{L22-def}) and inequalities (\ref{pos-A2-2}), (\ref{A2YW-pr}), (\ref{pos-C1-2})
and (\ref{C1YW-pr}), the $\Xi'$-positivity condition holds.

Therefore, Theorem \ref{thmA2} and \ref{thmC1} follow from Theorem \ref{polyhthm} and the above argument. \qed


\begin{thebibliography}{9}







\bibitem{BZ}
A.Berenstein, A.Zelevinsky, 
Tensor product multiplicities, canonical bases and totally positive varieties,
Invent. Math. 143, no. 1, 77--128 (2001).


\bibitem{GKS16} V.Genz, G.Koshevoy and B.Schumann, Combinatorics of canonical bases revisited: type  A, 
Selecta Math. (N.S.) 27, no. 4, Paper No. 67, 45 pp, (2021).

\bibitem{Ha}
T.Hayashi,
$Q$-analogues of Clifford and Weyl algebras—spinor and oscillator representations of quantum enveloping algebras,
Comm. Math. Phys. 127, no. 1, 129--144 (1990).





\bibitem{H1}
A.Hoshino,
Polyhedral realizations of crystal bases for quantum algebras of finite types,
J. Math. Phys. 46, no. 11, 113514,  31 pp, (2005). 

\bibitem{H2}
A.Hoshino,
Polyhedral realizations of crystal bases for quantum algebras of classical affine types,
J. Math. Phys. 54, no. 5, 053511, 28 pp, (2013).


\bibitem{JMMO}
M.Jimbo, K. C. Misra, T.Miwa, M.Okado,
Combinatorics of representations of  
$U_q(\widehat{\mathfrak{s}\mathfrak{l}}(n))$ at $q=0$,
Comm. Math. Phys. 136, no. 3, 543-–566 (1991).


\bibitem{Kac} V. G. Kac,
{\it Infinite-dimensional Lie algebras}, 
third edition. Cambridge University Press, Cambridge, xxii+400 pp, (1990).  



\bibitem{KaN} Y.Kanakubo, T.Nakashima, Adapted sequence for polyhedral realization of crystal bases,
Communications in Algebra,
Volume 48, Issue 11, pp4732--4766. (2020)


\bibitem{Kang}
S.-J. Kang, Crystal bases for quantum affine algebras and combinatorics of Young walls,
Proc. London Math. Soc. (3) 86, no. 1, 29-–69 (2003).


\bibitem{KK}
S.-J. Kang, J.-H. Kwon, 
Crystal bases of the Fock space representations and string functions,
J. Algebra 280, no. 1, 313--349 (2004).



\bibitem{KMM}
S.-J. Kang, K. C. Misra, T.Miwa,
Fock space representations of the quantized universal enveloping algebras
$U_q(C_l^{(1)})$, $U_q(A_{2l}^{(2)})$, and $U_q(D_{l+1}^{(2)})$,
J. Algebra 155, no. 1, 238-–251 (1993).







\bibitem{K0} M.Kashiwara, Crystalling the $q$-analogue of universal 
              enveloping algebras, Comm. Math. Phys.,
            {\it 133}, 249--260 (1990).


\bibitem{K1} M.Kashiwara,
 On crystal bases of the $q$-analogue of universal enveloping algebras,
	Duke Math. J., {\it 63} (2), 465--516 (1991).

\bibitem{K3}M.Kashiwara, 
The crystal base and Littelmann's refined Demazure character formula,
Duke Math. J., 71, no 3,  839--858 (1993).


\bibitem{KN}M.Kashiwara, T.Nakashima,
Crystal graphs for representations of the $q$-analogue of classical Lie algebras,
J. Algebra 165, no. 2, 295--345 (1994).


\bibitem{KS}\label{KS}
Kim, J.-A., Shin, D.-U., 
Monomial realization of crystal bases $B(\infty)$ for the quantum finite algebras,
Algebr. Represent. Theory 11, no. 1, 93--105 (2008).


\bibitem{Lit}
P.Littelmann, 
Cones, crystals, and patterns,
Transform. Groups 3, no. 2, 145--179 (1998).

\bibitem{L}
G.Lusztig, Canonical bases arising from quantized enveloping algebras, 
J. Amer. Math. Soc. 3, no. 2, 447--498 (1990). 

\bibitem{MM}
K. C. Misra, T.Miwa,
Crystal base for the basic representation of $U_q(\widehat{\mathfrak{s}\mathfrak{l}}(n))$,
Comm. Math. Phys. 134, no. 1, 79--88 (1990).

\bibitem{NZ}
T.Nakashima,  A.Zelevinsky, Polyhedral realizations of 
crystal bases for quantized Kac-Moody algebras,
Adv. Math. {\bf 131}, no. 1, 253--278, (1997). 

\bibitem{P}
A.Premat,
Fock space representations and crystal bases for $C_n^{(1)}$,
Journal of Algebra,
Volume 278, Issue 1, 227--241 (2004).


\end{thebibliography}
\end{document}